\documentclass[11pt]{amsart}

\usepackage{fullpage}
\usepackage{amsmath,amssymb,mathtools}
\usepackage[foot]{amsaddr}
\usepackage{xcolor,graphicx,caption,float,placeins,booktabs,array,longtable,supertabular,enumitem,pgfplots}
\pgfplotsset{compat=1.18}
\graphicspath{{figures/}{../figures/}}
\usepackage[backend=biber,giveninits=true,sorting=nyt,style=alphabetic,natbib=true,sortcites=true,maxcitenames=2,maxbibnames=10,url=false,doi=true,backref=false]{biblatex}
\renewbibmacro{in:}{\ifentrytype{article}{}{\printtext{\bibstring{in}\intitlepunct}}}

\usepackage[colorlinks=true,linkcolor=red,citecolor=blue,urlcolor=darkgray]{hyperref}

\newtheorem{theorem}{Theorem}[section]
\newtheorem{proposition}[theorem]{Proposition}
\newtheorem{lemma}[theorem]{Lemma}
\newtheorem{corollary}[theorem]{Corollary}
\newtheorem{conjecture}[theorem]{Conjecture}
\theoremstyle{remark}
\newtheorem{remark}[theorem]{Remark}
\theoremstyle{plain}

\newcommand{\R}{\mathbb{R}}
\newcommand{\sphere}{\mathbb S^2}
\newcommand{\ball}{\mathbb B^3}

\newcolumntype{C}[1]{>{\raggedright\arraybackslash}p{#1}}   
\newlength{\panelw}                                           
\newcommand{\paneltitle}[1]{\parbox[b]{\panelw}{\centering\scriptsize #1}}
\newcommand{\zoo}[2]{\href{https://github.com/braxtonosting/genus0-fbms-zoo/tree/main/surfaces/#1/}{#2}}   
\allowdisplaybreaks
\AtBeginDocument{\setlength\abovedisplayskip{4pt plus 1pt minus 1pt}\setlength\belowdisplayskip{4pt plus 1pt minus 1pt}\setlength\abovedisplayshortskip{1pt}\setlength\belowdisplayshortskip{2pt}}
\begin{document}
\author{Haoyu Chen}
\address{Department of Mathematics, University of Utah, Salt Lake City, UT}
\email{Haoyu.Chen@utah.edu}
\thanks{The authors acknowledge partial support from NSF DMS-2136198 and DMS-2513175.}

\author{Braxton Osting}
\email{osting@math.utah.edu}

\title[Spectral zoo of free boundary minimal surfaces]{A spectral zoo of genus-zero free boundary minimal surfaces with many boundary components}



\begin{abstract}
We study critical points of the length-normalized first Steklov eigenvalue $\bar\sigma_1$ on genus-zero surfaces with $b$ boundary components. By a theorem of Fraser and Schoen, their first eigenfunctions give free boundary minimal surfaces in the unit ball. After uniformization, the domain is the round sphere minus $b$ round caps, with a weight on the boundary circles. A layer potential formulation gives the eigenvalues, their derivatives and a certificate for each critical point, and an exact multipole system for the caps. We compute a spectral zoo of 268 numerically certified critical points with up to $b=2252$ boundary components, including the $27$ embedded examples listed by Schulz. With their symmetry groups and Morse indices, the critical points satisfy a signed count for $3\le b\le10$. Most belong to families built from rings of boundary components. As the rings fill in, these appear to converge to the free curve minimal surfaces of Karpukhin and Noutchegueme, unions of flat disks and catenoidal annuli that meet along circles on the sphere. We construct these limits, among them the critical drum and the cake stand, a flat disk with a catenoidal horn. The families approach their limits like $j^{-1}$ for rings of $j$ caps, and the equatorial family approaches $4\pi$ with deficit $4\pi b^{-1}e^{-b}$. For large $b$, the leading terms of the multipole system reduce to the logarithmic energy of the cap centers and a formal evaluation gives $\bar\sigma_1^\star(b)=8\pi-16\pi\,b^{-1}\log b+O(b^{-1})$.
\end{abstract}

\maketitle

\begin{figure}[p]\centering
\makebox[\textwidth][c]{\includegraphics[width=1.15\textwidth,height=0.78\textheight,keepaspectratio]{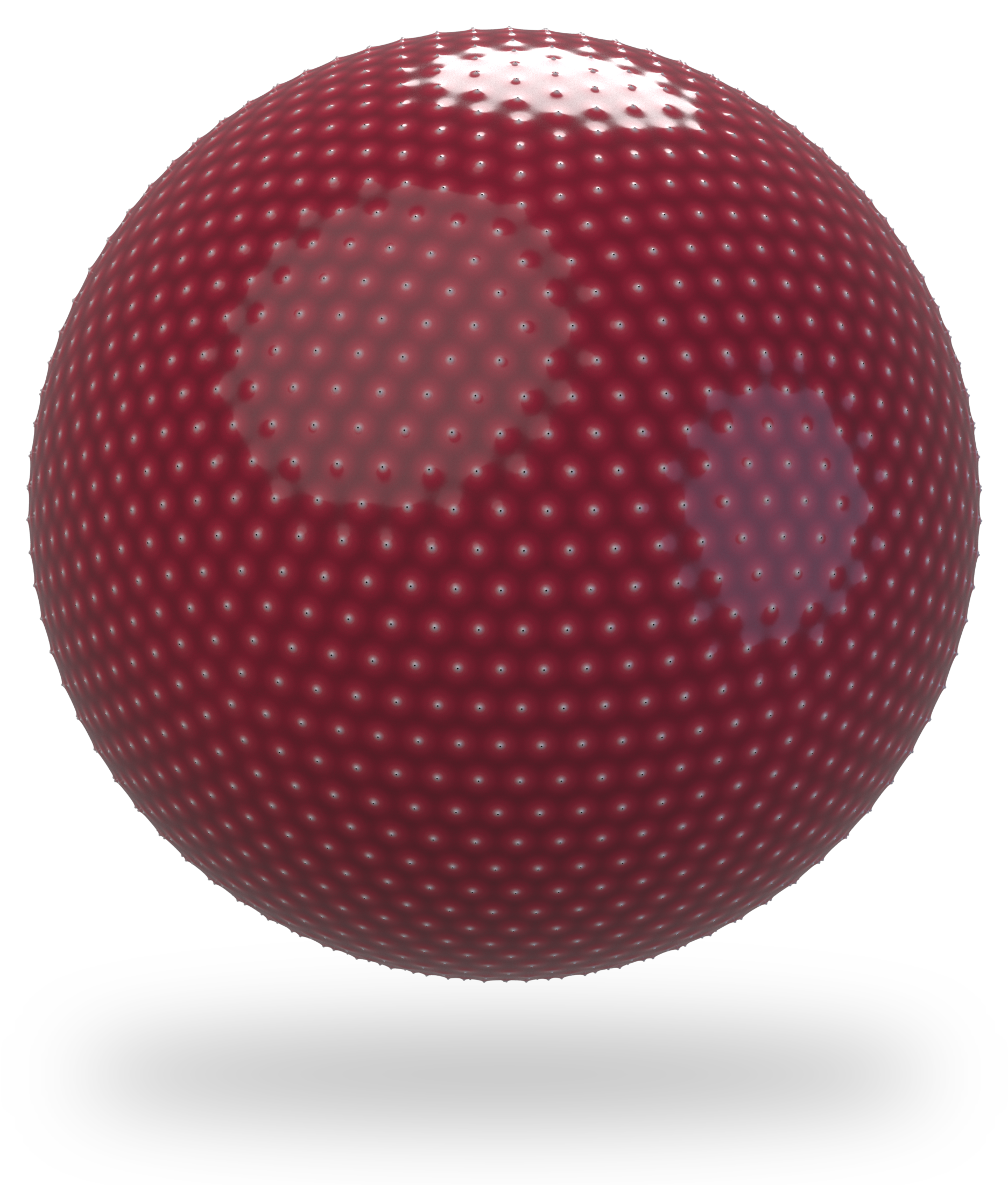}}
\caption{A free boundary minimal surface in the unit ball $\ball$ with $2252$ boundary components. It has genus zero and is a critical point of the first normalized Steklov eigenvalue $\bar\sigma_1$, computed and numerically certified by the method of this paper. The surface is a sheet just inside $\partial\ball$, joined to it by $2252$ small necks. Its color shows the distance $|X|$ from the origin, from red on the sheet to silver at the rim of each neck, and its boundary curves (silver) sit at the vertices of the geodesic icosahedron with $\nu=15$. It is the best configuration found for $b=2252$ (Table~\ref{tab:best}), the largest computed. The total length of its boundary curves is $\bar\sigma_1^\star=24.97$, and the area of the surface is $\tfrac12\bar\sigma_1^\star=12.49$. This value is close to the limit $8\pi\approx25.13$ that $\bar\sigma_1^\star(b)$ approaches as $b\to\infty$ (Section~\ref{ssec:large}).}\label{fig:headline}
\end{figure}

\section{Introduction}\label{sec:intro}
Let $(\Sigma,g)$ be a smooth, compact, connected, orientable Riemannian surface of genus $\gamma$ with $b\ge1$ boundary components. The Steklov eigenproblem is to find eigenpairs $(\sigma,u)$ with
\begin{subequations}\label{eq:steklov}
\begin{align}
&\Delta_gu=0 &&\text{ in }\Sigma, \label{eq:steklov-interior}\\
&\partial_{\nu_g}u=\sigma u &&\text{ on }\partial\Sigma. \label{eq:steklov-boundary}
\end{align}
\end{subequations}
The eigenvalues, counted with multiplicity, form a sequence $0=\sigma_0<\sigma_1\le\sigma_2\le\dots\to\infty$. They are the eigenvalues of the Dirichlet-to-Neumann operator \cite{Girouard_2017}. Since $\sigma_k$ scales like an inverse length, we normalize it by the boundary length $L_g(\partial\Sigma)$ and consider the extremal Steklov problem
\begin{equation}\label{eq:extremal}
\sigma^\star(\gamma,b):=\sup_g\bar\sigma_1(\Sigma,g),\qquad \bar\sigma_k(\Sigma,g):=\sigma_k(\Sigma,g)\,L_g(\partial\Sigma),
\end{equation}
the supremum taken over all smooth metrics on $\Sigma$. The problem goes back to Weinstock \cite{Weinstock_1954}, who proved $\bar\sigma_1\le2\pi$ for simply connected planar domains, with equality only for the disk. In general $\bar\sigma_k\le2\pi(\gamma+b+k-1)$ \cite{Fraser_2011,Girouard_2012,Karpukhin_2017}, so $\sigma^\star(\gamma,b)\le2\pi(\gamma+b)$. For genus zero, the supremum in \eqref{eq:extremal} is attained by a smooth metric \cite{Fraser_2015}, and in higher genus under a spectral gap condition \cite{Matthiesen_Petrides_2020}.

What sets the Steklov problem apart from other extremal eigenvalue problems is a link between its maximizers and minimal surfaces, discovered by Fraser and Schoen \cite{Fraser_2011,Fraser_2015}. A \emph{free boundary minimal surface} (FBMS) in the unit ball $\mathbb B^m$ is a critical point of area among surfaces with boundary on the sphere $\partial\mathbb B^m$. Equivalently, it is a minimal surface that meets the sphere orthogonally along its boundary. If $X\colon\Sigma\to\mathbb B^m$ is a free boundary minimal immersion and $g=X^*g_{\mathbb R^m}$, the coordinate functions satisfy $\Delta_gX=0$ in $\Sigma$ and $\partial_{\nu_g}X=X$ on $\partial\Sigma$. In other words, they solve \eqref{eq:steklov-interior}--\eqref{eq:steklov-boundary} with $\sigma=1$. Conversely, let $g$ maximize $\bar\sigma_1$. Then suitably scaled independent first eigenfunctions $u_1,\dots,u_m$, with $m$ the multiplicity of $\sigma_1$, define a branched conformal free boundary minimal immersion $X=(u_1,\dots,u_m)\colon\Sigma\to\mathbb B^m$ \cite[Theorem 1.1]{Fraser_2015}. We call $X$ the \emph{spectral map}. It maps $\partial\Sigma$ into the unit sphere and meets it orthogonally, since $|X|=1$ and $\partial_{\nu_g}X=\sigma_1X$ there. Because $X$ is conformal, arclength along $X(\partial\Sigma)$ is $\sigma_1$ times arclength for $g$, $ds_X=\sigma_1\,ds_g$. Up to the constant $\sigma_1$, the metric $g$ on the boundary is therefore the one induced by the surface, and
\[
L\bigl(X(\partial\Sigma)\bigr)=\sigma_1L_g(\partial\Sigma)=\bar\sigma_1.
\]
Since the coordinates are harmonic and $X$ is conformal, the area is half the boundary length:
\begin{equation}\label{eq:area}
\operatorname{Area}\bigl(X(\Sigma)\bigr)=\tfrac12\int_\Sigma|\nabla_gX|^2\,dA_g=\tfrac12\int_{\partial\Sigma}X\cdot\partial_{\nu_g}X\,ds_g=\tfrac12\sigma_1L_g(\partial\Sigma)=\tfrac12\bar\sigma_1 .
\end{equation}
The multiplicity $m_k$ of $\sigma_k$ on an orientable surface of genus $\gamma$ is at most $4\gamma+2k+1$, so $m_1\le3$ for a genus-zero surface \cite[Theorem 2.3]{Fraser_2015}. In this case, the maximizer of \eqref{eq:extremal} gives an FBMS in $\ball$, which for $b\ge2$ is not planar and is embedded \cite{Fraser_2015}. More generally, a genus-zero FBMS in $\ball$ with $b\ge2$ whose coordinate functions are first Steklov eigenfunctions is embedded, star-shaped with respect to the origin and free of branch points \cite[Proposition 8.1]{Fraser_2015}. So every critical point of $\bar\sigma_1$ in genus zero gives an embedded surface. In the other direction, Chodosh and Gianocca prove that every embedded FBMS of genus zero in $\ball$ has $\sigma_1=1$ \cite{Chodosh_Gianocca_2026}. Its coordinate functions are then first eigenfunctions, so every such surface comes from a critical metric of $\bar\sigma_1$. In particular, by \cite[Theorem 1.2]{Fraser_2015}, the critical catenoid is the only embedded free boundary minimal annulus in $\ball$. We use this correspondence to compute a \emph{spectral zoo} of genus-zero FBMS in $\ball$: surfaces whose coordinates are first Steklov eigenfunctions. Through the spectral map, every critical point of $\bar\sigma_1$ that we compute gives a genus-zero FBMS in $\ball$ of area $\tfrac12\bar\sigma_1$ (Theorem~\ref{thm:certificate} and \eqref{eq:area}), which is embedded by \cite[Proposition 8.1]{Fraser_2015}.

The values of $\sigma^\star(0,b)$ are known for small $b$ and in the limit $b\to\infty$. The equatorial disk realizes the supremum for $b=1$, with $\sigma^\star(0,1)=2\pi$, and the critical catenoid realizes it for $b=2$, with $\sigma^\star(0,2)\approx10.47$ \cite{Fraser_2015}. For $b\ge3$, no maximizer is known in closed form. The values increase strictly with $b$ \cite[Proposition 4.3]{Fraser_2015}, \cite{Matthiesen_Petrides_2018}. The maximizing surfaces converge to the sphere and the values converge to $8\pi$, the Hersch bound for the Laplacian on the sphere. Girouard and Lagac\'e proved this limit by homogenization \cite{Girouard_Lagace_2021}, and Karpukhin and Stern proved it independently \cite{Karpukhin_Stern_2021}. Karpukhin and Stern also showed that the deficit $8\pi-\sigma^\star(0,b)$ lies between two multiples of $b^{-1}\log b$ \cite[Theorem 1.8]{Karpukhin_Stern_2021}, and asked whether $\tfrac{b}{\log b}\bigl(8\pi-\sigma^\star(0,b)\bigr)$ has a limit \cite[Open Question 2]{Karpukhin_Stern_2021}. 

For finite $b\ge3$, there are two lines of previous work. On the computational side, Oudet, Kao and Osting \cite{Kao_2020} computed maximizers for $b=2,\dots,9$, $12$, $15$ and $20$ on the unit disk with $b-1$ circular holes. They observed that the boundary circles sit at the vertices of symmetric polyhedra resembling the minimizers of the Thomson problem. Girouard and Lagac\'e \cite{Girouard_Lagace_2021} computed the values of the configurations with Platonic symmetry. On the side of existence, Folha, Pacard and Zolotareva \cite{Folha_Pacard_Zolotareva_2017} and Ketover \cite{Ketover_2016} constructed genus-zero FBMS with $b$ boundary components, by gluing and by min-max. Their surfaces are two disks joined by $b$ equatorial bridges. Karpukhin, Kusner, McGrath and Stern \cite{karpukhin2024} proved that $\bar\sigma_1$ has maximizers among metrics invariant under a finite group, which gives FBMS with prescribed symmetry. Carlotto, Schulz and Wiygul \cite{Carlotto_Schulz_Wiygul_2025} and Franz and Schulz \cite{Franz_Schulz_2026} constructed many more by equivariant min-max with topological control, including pairs with the same topology and symmetry group. Schulz's online gallery \cite{Schulz_gallery} lists the numerically computed areas of the known embedded genus-zero FBMS with up to eight boundary components. These constructions give existence and symmetry, and by \cite{Chodosh_Gianocca_2026} each of these surfaces comes from a critical point of $\bar\sigma_1$. Which configuration is the maximizer for each $b$ has not been established. 

\emph{Results.} Most of our critical points belong to families with $j$ caps on each of a few latitudinal rings: pyramids, bipyramids, prisms, antiprisms and stacked rings, along with the equatorial family of $b$ caps on a great circle. We begin with their limits, which Section~\ref{sec:fcms} treats as free curve minimal surfaces of revolution in $\ball$. As the number $j$ of caps in a ring grows, each ring of necks shrinks to a circle on $\partial\ball$. The surfaces then appear to converge to unions of coaxial catenoidal annuli and flat disks glued along these circles (Figure~\ref{fig:fcms}). These unions are free curve minimal surfaces in the sense of \cite{Karpukhin_Noutchegueme_2025}. Lemma~\ref{lem:profile} reduces their definition to conditions on the profile curve, and Propositions~\ref{prop:disk}--\ref{prop:c3} construct five of them. Twice their areas are the limit values of the families (Sections~\ref{sec:fcms} and \ref{ssec:ringlimits}). The equatorial surfaces approach the doubly covered disk, with deficit $\tfrac{4\pi}{b}e^{-b}$. The prisms and antiprisms approach the critical drum of \cite{Karpukhin_Noutchegueme_2025}, whose value is $\bar\sigma_1=\tfrac{8\pi}{1+W(1/e)}=19.658536$, with $W$ the Lambert function (Proposition~\ref{prop:drum}). The antiprisms thus run from the critical catenoid at $j=1$ to the critical drum. A polar cap adds a new piece, a catenoidal horn that meets $\partial\ball$ orthogonally. The pyramids approach the cake stand, a flat disk with a horn, which is new to our knowledge (Proposition~\ref{prop:cake}). The bipyramids approach the wheel, the union of two catenoidal annuli of \cite{Carlotto_Schulz_Wiygul_2025} (Proposition~\ref{prop:wheel}). Kapouleas and Zou desingularize the unions of two flat disks and many catenoidal annuli into free boundary minimal surfaces \cite{Kapouleas_Zou_2021}, and Carlotto, Schulz and Wiygul desingularize the wheel \cite[Theorem 6.1]{Carlotto_Schulz_Wiygul_2025}. As far as we know, the drum and the cake stand have not been desingularized.

Our second result concerns the sphere with $b$ round caps removed. Proposition~\ref{prop:multipole} writes the Steklov problem there as an exact system for the charges and multipoles of the caps (Section~\ref{ssec:multipole}). When the weight is constant on each boundary circle, the flux of a cap is proportional to the mean of the eigenfunction on its circle, and each multipole is proportional to a Taylor coefficient of the field of the other caps (Corollary~\ref{cor:const}). For large $b$, keeping only the charges of the caps, or the charges and the dipoles, gives the monopole and dipole approximations (Section~\ref{ssec:large}). In the monopole approximation, the centers of the boundary curves enter only through their logarithmic energy, and a smaller energy gives a larger value. The dipole approximation suggests $\bar\sigma_1^\star(b)=8\pi-16\pi\,b^{-1}\log b+O(b^{-1})$, so that the limit in the question of Karpukhin and Stern is $16\pi$.

Our computations give $268$ critical points of $\bar\sigma_1$, with $b$ from $2$ to $2252$ (Section~\ref{sec:experiments}, Table~\ref{tab:catalog}). We search the polyhedral classes, these families, classes with the symmetry of the Fekete configurations, the minimizers of the logarithmic energy on the sphere, and an unrestricted optimizer. Each critical point is certified by a stationarity multiplier in the full parameter space, and the resulting surface is verified to be a free boundary minimal surface (Theorem~\ref{thm:certificate}, Section~\ref{ssec:checks}). The certificates are numerical. Their residuals are reported in Table~\ref{tab:best}, and the discretization error is estimated by refinement rather than bounded rigorously (Section~\ref{ssec:validation}). Table~\ref{tab:map} arranges the ring families by their numbers of rings and polar caps and names the surface each converges to. Figures~\ref{fig:type-polyhedral}--\ref{fig:type-stacked} show the surfaces of each type, Figures~\ref{fig:best-a}--\ref{fig:best-last} show the best surface for each $b$, and Figure~\ref{fig:families} plots the values against $b$. The configurations, certificates and surfaces are in an online repository \cite{Chen_Osting_genus0_code}, linked from their names in the tables and figures. The computations reproduce the values of \cite{Girouard_Lagace_2021} and \cite{Kao_2020} and contain all $27$ surfaces listed in \cite{Schulz_gallery} with $b\ge2$ (see Section~\ref{ssec:zoo} and Table~\ref{tab:zoo}). Beyond these, the symmetry classes reach much larger $b$. With a symmetry-reduced solver, we compute the best configuration found for each $b$ up to $b=2252$, where $\bar\sigma_1^\star=24.97$ (Table~\ref{tab:best}). These surfaces are spheres with $b$ small necks to $\partial\ball$ (Figure~\ref{fig:headline}).

We also compute the Morse index $\operatorname{ind}(\Sigma)$ of each critical point $\Sigma$, the number of directions in which $\bar\sigma_1$ increases. Suppose that $\bar\sigma_1$, maximized over the weights, is a Morse function on the space of configurations of $b$ caps up to M\"obius maps, and that it decreases as a configuration degenerates. Then the Poincar\'e--Hopf theorem gives
\begin{equation}\label{eq:morse-count}
\sum_{\Sigma}\frac{(-1)^{\operatorname{ind}(\Sigma)}}{|G_\Sigma|}=\frac{(-1)^{b-1}}{2\,b(b-1)(b-2)} ,
\end{equation}
where the sum is over the surfaces $\Sigma$ with $b$ boundary components and $G_\Sigma$ is the symmetry group of $\Sigma$ (Section~\ref{ssec:catalog}). The surfaces of Table~\ref{tab:catalog} satisfy \eqref{eq:morse-count} for $3\le b\le10$, and twenty-two of the surfaces with $b=9$ and $10$ were found because it failed without them. This leads us to the following conjecture.
\begin{conjecture}\label{conj:morse}
For every $b\ge3$, the critical points of $\bar\sigma_1$ with $b$ boundary components satisfy \eqref{eq:morse-count}. Each has Morse index at most $b-3$, with equality only for the equatorial $b$-gon.
\end{conjecture}
The count is not proved. The function $\bar\sigma_1$ is not known to decrease as a configuration degenerates, and it is not smooth at its critical points, where the first eigenvalue has multiplicity at least $3$. Even then, the index bound would not follow from the count.

Several of our statements are observations or formal calculations rather than theorems. The convergence of the ring families to the surfaces of Section~\ref{sec:fcms} is observed, not proved, and so is their rate $j^{-1}$ (Section~\ref{ssec:ringlimits}). The deficit $\tfrac{4\pi}{b}e^{-b}$ of the equatorial family comes from matching the bridges of \cite{Folha_Pacard_Zolotareva_2017} to half-catenoids (Section~\ref{ssec:equatorial}). The accuracy of the monopole and dipole approximations is measured, not proved, and the expansion of $\bar\sigma_1^\star(b)$ is a formal evaluation of them at the Fekete configurations (Section~\ref{ssec:large}).

\emph{Method.} The computations follow the approach developed for the genus-one problem in \cite{Chen_Osting_2026}, adapted to the sphere. By the Koebe uniformization theorem \cite{HeSchramm1993}, $\Sigma$ is conformally the round sphere with $b$ disjoint round caps removed. Since $\bar\sigma_1$ depends on the metric only through the boundary weight $\rho$, the problem becomes a maximization over the centers and radii of the caps and the Fourier coefficients of $\log\rho$ on each circle. A normal form fixes the M\"obius freedom (problem \eqref{eq:problem}, Section~\ref{ssec:uniformization}). The eigenvalues are those of a pencil of single and double layer potentials for the Green's function $-\tfrac1{2\pi}\log|x-y|$ of the sphere \cite{Gemmrich_Nigam_Steinbach_2008,Kropinski_Nigam_2014} (Sections~\ref{ssec:layer} and \ref{ssec:implementation}). Theorem~\ref{thm:derivatives} gives the derivatives of $\bar\sigma_1$ with respect to the centers, radii and weights (Section~\ref{ssec:derivatives}). Every optimum is certified by a stationarity multiplier $\Lambda$ computed in the \emph{full} parameter space, symmetry-breaking directions included. The certificate gives the spectral map and so the surface (Section~\ref{ssec:optimization}), which we then verify as an FBMS (Section~\ref{ssec:checks}). Unlike genus one, where the maximizers of \cite{Chen_Osting_2026} lie in $\mathbb B^5$ and $\mathbb B^6$, in genus zero the first eigenspace has dimension at most three, so every certified critical point of rank three is an FBMS in $\ball$.

The closest earlier work is \cite{Kao_2020}, which computes maximizers on the unit disk with $b-1$ holes. There, a symmetry that exchanges the outer circle with a hole is a M\"obius transformation and not an isometry, so symmetric configurations can be recognized in the output but not imposed. On the sphere, every symmetry is an isometry, and a symmetry class is described by a handful of parameters. In place of the method of particular solutions, the explicit kernel and the Kress quadrature \cite{Kress_2014} give all eigenvalues to rounding. With exact derivatives on the parameters of a class, most of our optima are stationary to rounding rather than to $10^{-6}$, and the certificate gives the surface.

\begin{table}[tbp]\centering\footnotesize\setlength{\tabcolsep}{1.2pt}\renewcommand{\arraystretch}{1.25}
\begin{tabular}{@{}l|ccc@{}}\toprule
& no polar cap & one polar cap & two polar caps\\\midrule
one ring & \begin{tabular}[t]{@{}c@{}}equatorial $b$-gon (Fig.~\ref{fig:type-equatorial})\\ doubled disk, $4\pi$\end{tabular} & \begin{tabular}[t]{@{}c@{}}pyramid (Fig.~\ref{fig:type-pyramid})\\ cake stand, $15.782$\end{tabular} & \begin{tabular}[t]{@{}c@{}}bipyramid (Fig.~\ref{fig:type-bipyramid})\\ wheel, $17.543$\end{tabular}\\\addlinespace[6pt]
two rings & \begin{tabular}[t]{@{}c@{}}prism, antiprism (Figs.~\ref{fig:type-prism}, \ref{fig:type-antiprism})\\ critical drum, $19.659$\end{tabular} & \begin{tabular}[t]{@{}c@{}}prism, antiprism + 1PC (Figs.~\ref{fig:type-prism}, \ref{fig:type-antiprism})\\ drum with one horn, $20.486$\end{tabular} & \begin{tabular}[t]{@{}c@{}}prism, antiprism + 2PC (Figs.~\ref{fig:type-prism}, \ref{fig:type-antiprism})\\ drum with two horns, $21.070$\end{tabular}\\\addlinespace[6pt]
three rings & \begin{tabular}[t]{@{}c@{}}stacked rings (Fig.~\ref{fig:type-stacked})\\ drum stack, $22.229$\end{tabular} & \begin{tabular}[t]{@{}c@{}}stacked rings + 1PC (Fig.~\ref{fig:type-stacked})\\ drum stack with one horn, $22.506$\end{tabular} & \begin{tabular}[t]{@{}c@{}}stacked rings + 2PC (Fig.~\ref{fig:type-stacked})\\ drum stack with two horns, $22.731$\end{tabular}\\\bottomrule
\end{tabular}
\caption{Map of the zoo (Sections~\ref{ssec:families} and \ref{ssec:catalog}): the ring families by their numbers of rings and polar caps (PC), with the limit network as $j\to\infty$ and its value $\bar\sigma_\infty$ (Section~\ref{sec:fcms}). Figure~\ref{fig:fcms} shows the meridian sections of the limit networks. The critical catenoid ($10.475$) is the bipyramid with $j=0$ and the antiprism with $j=1$. The classes are defined in Table~\ref{tab:classes}, and the values are in Table~\ref{tab:catalog} and Figure~\ref{fig:families}.} \label{tab:map}
\end{table}

\emph{Outline.} Section~\ref{sec:fcms} constructs the free curve minimal surfaces to which the ring families converge. Section~\ref{sec:method} uniformizes the problem to the punctured sphere, introduces the symmetry classes and the layer potential formulation, and proves the multipole form, the derivative formulas and the stationarity certificate. Section~\ref{sec:numerics} describes the discretization, the symmetry-reduced solver, the optimizer, the certification rule and the Morse index, and validates them. Section~\ref{sec:experiments} reports the computations family by family, the best configuration for each $b$, the rate at large $b$ and a map of the zoo. Section~\ref{sec:discussion} concludes with a discussion.

\emph{Notation.} We write $[n]:=\{1,\dots,n\}$. A point of $\sphere$ has polar angle $\beta\in[0,\pi]$ from the north pole $e_3$ and azimuth $\phi$ from $e_1$. Point groups are written in Schoenflies notation ($C_{jv}$, $D_{jh}$, $D_{jd}$, $T_d$, $O_h$, $I_h$, $I$). A superscript ${}^\star$ marks quantities at a computed critical configuration, such as $\bar\sigma_1^\star$.

\section{Free curve minimal surfaces}\label{sec:fcms}
\begin{figure}[tbp]
\input{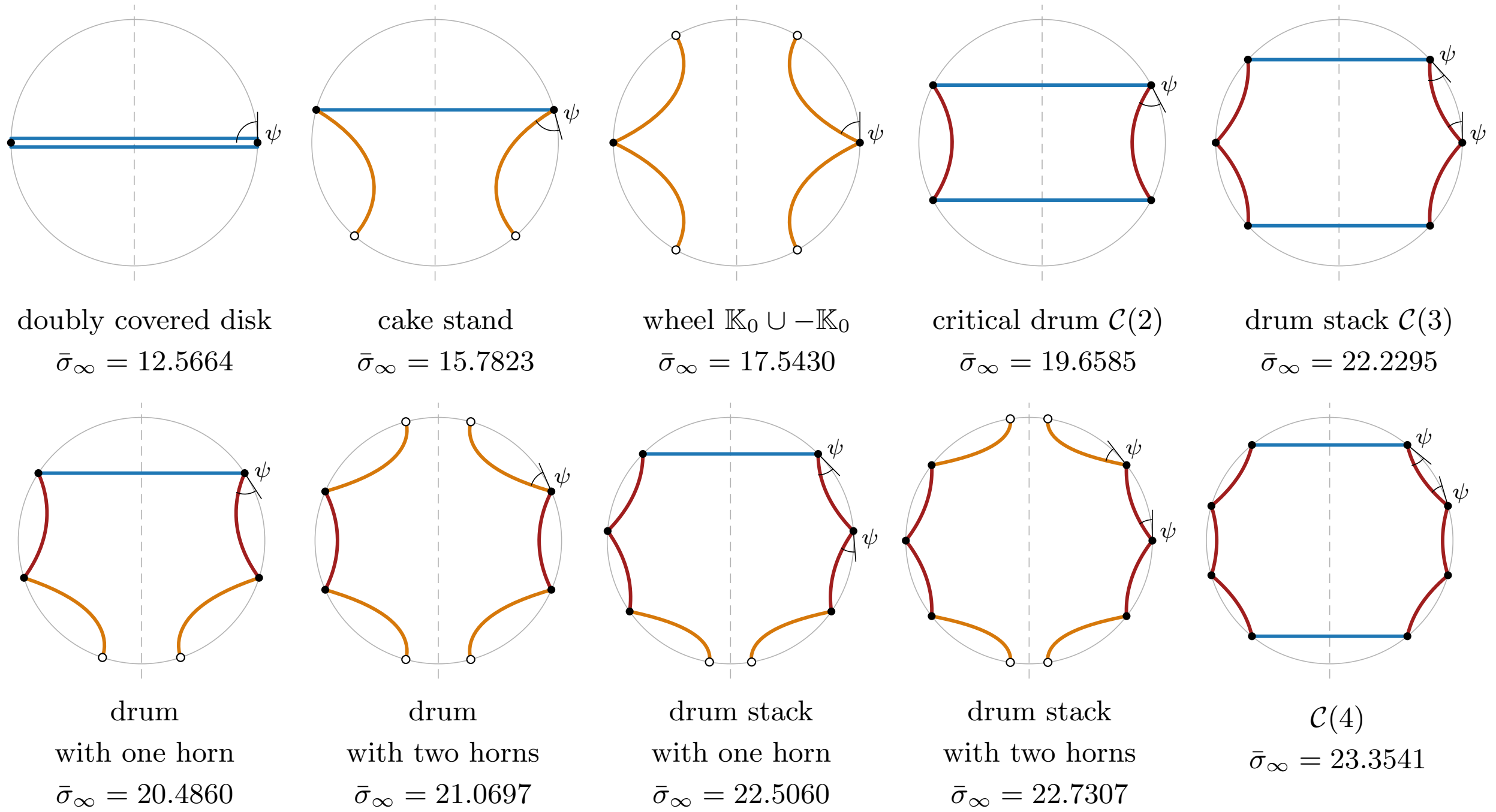}
\caption{Meridian sections of the free curve minimal surfaces of Section~\ref{sec:fcms}, with their values $\bar\sigma_\infty$. They are the limit networks of the ring families in Table~\ref{tab:map}, and Table~\ref{tab:limitnetworks} lists their data. Flat disks are blue, catenoidal bands red and horns orange, inside $\partial\ball$ (gray), with the axis of revolution dashed. Filled dots are the circles of $\Gamma$, where two sheets meet and make the same angle $\psi$ with $\partial\ball$. Open dots are boundary circles, where a horn meets $\partial\ball$ orthogonally. The two sheets of the doubly covered disk are drawn apart. The networks with one horn are drawn with the horn at the bottom. Top row: the surfaces of Propositions~\ref{prop:disk}--\ref{prop:c3}. Bottom row: networks computed from Lemma~\ref{lem:profile}, described after Proposition~\ref{prop:c3}.}\label{fig:fcms}
\end{figure}

In the ring families of Table~\ref{tab:map}, each ring contains $j$ caps. As $j$ grows, the surfaces of these families appear to converge to unions of flat disks and catenoidal annuli that meet along circles on $\partial\ball$ (Section~\ref{ssec:families}). If the convergence is smooth away from the circles, each piece of a limit must be a flat disk or a catenoidal annulus, since the limit is rotationally symmetric and stationary. These unions are free curve minimal surfaces in the sense of \cite{Karpukhin_Noutchegueme_2025}, extended in Section~\ref{ssec:fcms-def} to allow boundary circles, along which the surface meets $\partial\ball$ orthogonally. Boundary circles occur when the family has polar caps. In this section, we construct five of these surfaces explicitly: the doubly covered disk and the critical drum in closed form (Propositions~\ref{prop:disk} and~\ref{prop:drum}), the wheel and the drum stack $\mathcal C(3)$ up to the root of an equation in one variable (Propositions~\ref{prop:wheel} and~\ref{prop:c3}), and the cake stand up to a solution of the system \eqref{eq:cake} of two equations (Proposition~\ref{prop:cake}). Figure~\ref{fig:fcms} shows these five and the limits of the other ring families we computed, and Table~\ref{tab:limitnetworks} lists their data.

\subsection{Definition and the transmission problem}\label{ssec:fcms-def}
Let $M$ be a compact surface, possibly with boundary, and let $\Gamma$ be a finite union of disjoint closed curves in the interior of $M$. A \emph{free curve minimal surface} in $\ball$ is a continuous map $X\colon M\to\overline{\ball}$ with $X(\Gamma\cup\partial M)\subset\partial\ball$ that is conformal and harmonic on $M\setminus\Gamma$ and satisfies
\begin{subequations}\label{eq:fc}
\begin{align}
\partial_{\nu^+}X+\partial_{\nu^-}X\ &\parallel\ X &&\text{on }\Gamma,\label{eq:fc-curve}\\
\partial_\nu X\ &\parallel\ X &&\text{on }\partial M .\label{eq:fc-bdry}
\end{align}
\end{subequations}
Here $\nu^\pm$ are the outward unit normals of the two sides of $\Gamma$ and $\nu$ is the outward unit normal of $\partial M$, for any metric in the conformal class of $M$. For closed $M$, this is the definition of \cite[Definitions 1.1 and 1.5]{Karpukhin_Noutchegueme_2025}. The image is a union of minimal surfaces. Along $X(\Gamma)$ two sheets meet on $\partial\ball$, and the sum of their conormals is normal to $\partial\ball$. We call this the balancing condition. Along $X(\partial M)$, the surface meets $\partial\ball$ orthogonally, just as a free boundary minimal surface does.

These surfaces come with an eigenvalue problem. Let $g$ be a metric on $M$ and $\lambda>0$ a density on $\Gamma\cup\partial M$. The transmission--Steklov problem for a continuous function $u$ is
\begin{subequations}\label{eq:transmission}
\begin{align}
\Delta_gu&=0 &&\text{in }M\setminus\Gamma,\label{eq:tr-int}\\
\partial_{\nu^+}u+\partial_{\nu^-}u&=\sigma\lambda u &&\text{on }\Gamma,\label{eq:tr-curve}\\
\partial_\nu u&=\sigma\lambda u &&\text{on }\partial M .\label{eq:tr-bdry}
\end{align}
\end{subequations}
Its eigenvalues are $0=\sigma_0<\sigma_1\le\sigma_2\le\cdots$, and we normalize them by the total weight, $\bar\sigma_k=\sigma_k\int_{\Gamma\cup\partial M}\lambda\,ds$, where $ds$ is the arclength of $g$. For closed $M$ and $\lambda\equiv1$, these are the Steklov transmission eigenvalues of \cite{Karpukhin_Noutchegueme_2025}, and for $\Gamma=\emptyset$ the problem is the weighted Steklov problem of Section~\ref{ssec:uniformization}. Since only the Dirichlet energy and the measure $\lambda\,ds$ enter, the problem is conformally invariant.

\begin{lemma}\label{lem:fcms-eig}
Let $X$ be a free curve minimal surface and $g$ the metric it induces on $M\setminus\Gamma$. Put $\lambda=X\cdot(\partial_{\nu^+}X+\partial_{\nu^-}X)$ on $\Gamma$ and $\lambda=1$ on $\partial M$. Then the coordinates of $X$ solve \eqref{eq:transmission} with $\sigma=1$, and $\int_{\Gamma\cup\partial M}\lambda\,ds=2\operatorname{Area}(X)$, with the area counted with multiplicity.
\end{lemma}
\begin{proof}
The coordinates of a conformal harmonic map are harmonic for the induced metric, which is \eqref{eq:tr-int}. On $\Gamma$ we have $|X|=1$, so \eqref{eq:fc-curve} gives $\partial_{\nu^+}X+\partial_{\nu^-}X=\lambda X$, which is \eqref{eq:tr-curve} with $\sigma=1$. On $\partial M$, the conormal is a unit vector parallel to $X$ that points out of $\ball$, so it equals $X$, which is \eqref{eq:tr-bdry}. On each piece $\Sigma$ of the image, $\operatorname{div}_\Sigma X=2$, and the divergence theorem gives $2\operatorname{Area}(\Sigma)=\int_{\partial\Sigma}X\cdot\partial_\nu X\,ds$. Summing over the pieces gives the claim.
\end{proof}
So the normalized eigenvalue of the coordinates is twice the area. In the examples below $\lambda>0$.

\subsection{Surfaces of revolution}\label{ssec:fcms-rev}
With $t=\log\cot\tfrac\beta2$, the round metric of $\sphere$ is $\operatorname{sech}^2t\,(dt^2+d\phi^2)$. So $(t,\phi)$ is a conformal coordinate on $\sphere$ minus the poles, the circles of latitude are the lines of constant $t$, and the north pole is $t=+\infty$. Let $\Gamma$ consist of the circles $t=t_i$, $i=1,\dots,n$, with $t_1<\dots<t_n$, and let $M$ be the part $t_-\le t\le t_+$ of $\sphere$. A finite end $t_\pm$ is a boundary circle of $M$, the rim of a polar cap, and an infinite end means that $M$ contains the pole. Write the measure $\lambda\,ds$ as $a_i\,d\phi$ on the circle $t=t_i$ and as $a_\pm\,d\phi$ on a boundary circle. Then \eqref{eq:transmission} reads
\begin{subequations}\label{eq:cyl}
\begin{align}
u_{tt}+u_{\phi\phi}&=0 &&\text{for }t_-<t<t_+,\ t\ne t_i,\label{eq:cyl-int}\\
u_t(t_i^+)-u_t(t_i^-)&=-\sigma a_i\,u(t_i) &&\text{on the circles }t=t_i,\label{eq:cyl-curve}\\
\pm u_t(t_\pm)&=\sigma a_\pm\,u(t_\pm) &&\text{at a finite end }t_\pm,\label{eq:cyl-bdry}
\end{align}
\end{subequations}
with $u$ bounded at an infinite end and $\bar\sigma=2\pi\sigma\bigl(\sum_ia_i+\sum_\pm a_\pm\bigr)$. We look for rotationally symmetric maps
\begin{equation}\label{eq:rev}
X(t,\phi)=\bigl(R(t)\cos\phi,\ R(t)\sin\phi,\ Z(t)\bigr),\qquad R>0 ,
\end{equation}
and write $P=(R,Z)$ for the profile and $P'=(R',Z')$ for its tangent. 

\begin{lemma}\label{lem:profile}
The map \eqref{eq:rev} is a free curve minimal surface if and only if the following hold.
\begin{enumerate}[label=(\roman*),nosep]
\item (Conformality and harmonicity) On each interval between consecutive circles and ends, $R''=R$, $Z''=0$ and $|P'|=R$.
\item (Balancing) At each circle, $P$ is continuous, $|P(t_i)|=1$ and $P'(t_i^+)-P'(t_i^-)=-a_iP(t_i)$ for a constant $a_i>0$.
\item (Free boundary) At a finite end, $|P(t_\pm)|=1$ and $\pm P'(t_\pm)=a_\pm P(t_\pm)$ for a constant $a_\pm>0$.
\item (Pole) At an infinite end, $R\to0$ and $Z$ is constant.
\end{enumerate}
If (i)--(iv) hold, the surface lies in $\overline{\ball}$, an unbounded piece is a flat disk, and every other piece is part of a catenoid $R=\varrho\cosh\tfrac{Z-Z_*}{\varrho}$ about the axis, with waist radius $\varrho$ at height $Z_*$. At the circle $t_i$, of radius $R_i=R(t_i)$, the two sheets make the same angle $\psi_i$ with $\partial\ball$, and $a_i=2R_i\sin\psi_i$. A boundary circle has radius $a_\pm$. The coordinates of $X$ solve \eqref{eq:cyl} with $\sigma=1$, and
\begin{equation}\label{eq:fcms-value}
\bar\sigma=2\pi\Bigl(\sum_ia_i+\sum_\pm a_\pm\Bigr)=2\operatorname{Area}(X).
\end{equation}
\end{lemma}
\begin{proof}
Away from the poles, $(t,\phi)$ are conformal coordinates on $M$, so we check conformality, harmonicity and \eqref{eq:fc} with the flat metric $dt^2+d\phi^2$. The vectors $X_t=(R'\cos\phi,R'\sin\phi,Z')$ and $X_\phi=(-R\sin\phi,R\cos\phi,0)$ are orthogonal with $|X_t|=|P'|$ and $|X_\phi|=R$, and $X_{tt}+X_{\phi\phi}=\bigl((R''-R)\cos\phi,(R''-R)\sin\phi,Z''\bigr)$, so conformality and harmonicity are (i). At the circle $t_i$ the outward unit normals of the two sides are $\pm\partial_t$, so \eqref{eq:fc-curve} reads $X_t(t_i^-)-X_t(t_i^+)=a_iX$ with $a_i$ a function of $\phi$. Its first two components give $R'(t_i^-)-R'(t_i^+)=a_iR(t_i)$, so $a_i$ is constant because $R>0$, and with the third component this is the jump condition of (ii). At a finite end the outward unit normal is $\pm\partial_t$, and \eqref{eq:fc-bdry} is $\pm P'(t_\pm)=a_\pm P(t_\pm)$. Since $X$ maps $M$ into $\overline{\ball}$ and $|P|=1$ at the circles and the finite ends, $|P|^2$ has a maximum there, one-sided at an end. At $t_i$ this gives $a_i=P\cdot P'(t_i^-)-P\cdot P'(t_i^+)\ge0$. If $a_i=0$, then $P\cdot P'(t_i)=0$ while $(|P|^2)''=2|P'|^2+2RR''+2ZZ''=4R^2>0$ by (i), so $|P|>1$ near $t_i$, which is impossible. At $t_\pm$ it gives $a_\pm\ge0$, and $|a_\pm|=|P'(t_\pm)|=R(t_\pm)>0$. Continuity of $P$ and $|P|=1$ at the circles and the finite ends say that $X$ is continuous and maps $\Gamma\cup\partial M$ into $\partial\ball$, and by (iv) $X$ tends to a point on the axis at a pole in $M$. This gives the equivalence, once we show that $X$ maps $M$ into $\overline{\ball}$ and that an unbounded piece is a flat disk, smooth at the pole.

By (i), on each piece $Z$ is affine and $R=pe^t+qe^{-t}$ for constants $p$ and $q$, and $|P'|=R$ reads $Z'^2=4pq$. On an unbounded piece, say $t\ge t_n$, (iv) gives $p=0$ and $Z=Z(t_n)$, so $X=\bigl(R_ne^{-(t-t_n)}\cos\phi,\,R_ne^{-(t-t_n)}\sin\phi,\,Z(t_n)\bigr)$ is the flat disk of radius $R_n$ at height $Z(t_n)$. It is affine in the stereographic coordinate $e^{-t}e^{i\phi}=\tan\tfrac\beta2\,e^{i\phi}$ about the pole, so it is smooth there, and the case $t\le t_1$ is the same. On a bounded piece $Z'\ne0$, since otherwise $pq=0$ and $R$ would be strictly monotone between two points of $\partial\ball$ at the same height. So $p,q>0$ and $R=\varrho\cosh(t-t_*)$ with $\varrho=|Z'|$, and with $Z_*=Z(t_*)$ this is the catenoid $R=\varrho\cosh\tfrac{Z-Z_*}{\varrho}$. On a bounded piece $|P|^2$ is strictly convex and equals $1$ at both ends, and on an unbounded piece $|P|^2=R^2+Z(t_n)^2$ decreases from $1$, so $X$ maps $M$ into $\overline{\ball}$. At $t_i$ both tangents have length $R_i$ and differ by $-a_iP$ with $|P|=1$, so $P'(t_i^-)\cdot P=-P'(t_i^+)\cdot P=\tfrac{a_i}2$, which gives the common angle with $\sin\psi_i=\tfrac{a_i}{2R_i}$. For the induced metric $R^2(dt^2+d\phi^2)$, the measure of Lemma~\ref{lem:fcms-eig} is $a_i\,d\phi$ on the circle $t_i$ and $a_\pm\,d\phi$ on a boundary circle, which gives the last two claims.
\end{proof}
We call a catenoidal piece between two circles of $\Gamma$ a \emph{band}, and one between a circle of $\Gamma$ and a boundary circle a \emph{horn}. Each is conformally a flat cylinder, whose \emph{length} is the distance in $t$ between its two circles. The translations in $t$ are M\"obius transformations of $\sphere$, so we fix one circle at $t=0$. Conditions (ii)--(iv) are also the first-order conditions for $\bar\sigma_1$ to be critical in \eqref{eq:cyl} when $\sigma=1$ is its first nonzero eigenvalue, as in \cite[Theorems 1.4 and 1.7]{Karpukhin_Noutchegueme_2025} for closed $M$. Stationarity in the weights puts the circles on $\partial\ball$, and stationarity in the positions gives the equal angles.

The equations match the unknowns. For a surface with $n$ circles and $n_c\in\{0,1,2\}$ boundary circles, the unknowns are the $n-1+n_c$ lengths of the bounded pieces, the $n+n_c$ weights $a_i$ and $a_\pm$, and the scale factors of $R$ and $Z$, $2n+2n_c+1$ in all. Apart from $|P'|=R$ and $|P|=1$, conditions (i)--(iv) are linear and do not mix $R$ and $Z$. They are \eqref{eq:cyl} with $\sigma=1$, for $R$ in the modes $e^{\pm i\phi}$ and for $Z$ in the rotationally symmetric mode. Shooting from one end leaves one condition at the other end for each, $2$ equations in all. The condition $|P|=1$ gives $n+n_c$ equations. Since $|P'|^2-R^2=Z'^2-4pq$ is constant on each piece and vanishes on an unbounded one, $|P'|=R$ gives one equation for each of the $n-1+n_c$ bounded pieces. This makes $2n+2n_c+1$ equations.

Each proposition below gives a profile, and its proof checks conditions (i)--(iv). The top row of Figure~\ref{fig:fcms} shows their meridian sections, and Table~\ref{tab:limitnetworks} lists the heights of their circles and their angles $\psi_i$, which follow from $\sin\psi_i=\tfrac{a_i}{2R_i}$. In the propositions, we index the circles from the equatorial plane outward, as in Table~\ref{tab:limitnetworks}, with $t_0=0$ when $t=0$ is a circle. Most pieces are of the form $R=A(\cosh t-s\sinh t)$ with $A>0$ and $|s|<1$, and $Z$ affine. Then $R''=R$, $Z''=0$ and $R>0$, and $|P'|^2-R^2=A^2(s^2-1)+Z'^2$, so (i) holds if and only if $Z'^2=A^2(1-s^2)$. If one side of a circle $t_i$ is a flat disk, then $P'(t_i^-)=(R_i,0)$ or $P'(t_i^+)=(-R_i,0)$, so the proof of Lemma~\ref{lem:profile} gives $\tfrac{a_i}2=R_i^2$ and $\sin\psi_i=R_i$. Then $\psi_i$ is the angular distance of the circle from the nearer pole.

\begin{proposition}[The doubly covered disk]\label{prop:disk}
Let $M=\sphere$ with the circle $t_0=0$. Then $P=(e^{-|t|},0)$ is a free curve minimal surface, the equatorial disk covered twice. The weight of the equator is $a_0=2$, and $\bar\sigma=4\pi$.
\end{proposition}
\begin{proof}
Conditions (i) and (iv) of Lemma~\ref{lem:profile} are clear, and (iii) is vacuous since $M$ has no boundary. For (ii), $P(0)=(1,0)$ and $P'(0^+)-P'(0^-)=(-2,0)=-2P(0)$.
\end{proof}
This is $\mathcal C(1)$ of \cite[Theorem 1.10]{Karpukhin_Noutchegueme_2025}. Its two sheets coincide and lie in a plane.

\begin{proposition}[The cake stand]\label{prop:cake}
Let $(a_0,L)$ with $0<a_0<2$ and $a_0L>1$ solve
\begin{subequations}\label{eq:cake}
\begin{align}
(a_0L-1)\bigl(\sinh L+(1-a_0)\cosh L\bigr)&=a_0\bigl(\cosh L+(1-a_0)\sinh L\bigr),\label{eq:cake-a}\\
\tfrac{a_0}2\bigl(\cosh L+(1-a_0)\sinh L\bigr)^2+\bigl(1-\tfrac{a_0}2\bigr)(1-a_0L)^2&=1 .\label{eq:cake-b}
\end{align}
\end{subequations}
Put $R_0=\sqrt{\tfrac{a_0}2}$ and $Z_0=\sqrt{1-\tfrac{a_0}2}$. Let $M$ be the part $t\le L$ of $\sphere$ with the circle $t_0=0$. Then
\begin{equation*}
P=\bigl(R_0e^{t},\ Z_0\bigr)\ \ (t\le0),\qquad P=\bigl(R_0(\cosh t+(1-a_0)\sinh t),\ Z_0(1-a_0t)\bigr)\ \ (0\le t\le L)
\end{equation*}
is a free curve minimal surface. It consists of a flat disk of radius $R_0$ and a horn from the rim of the disk to a boundary circle of radius $a_+=\tfrac{a_0}{a_0L-1}$, and $\bar\sigma=2\pi(a_0+a_+)$. The system \eqref{eq:cake} has the solution $a_0=1.858284$, $L=2.068251$, for which $\bar\sigma=15.782277$.
\end{proposition}
\begin{proof}
We check the conditions of Lemma~\ref{lem:profile}. The disk satisfies (i) and (iv). The horn has $A=R_0$, $s=a_0-1$ and $Z'^2=a_0^2Z_0^2=A^2(1-s^2)$, so it satisfies (i). At the circle, $|P(0)|^2=R_0^2+Z_0^2=1$ and $P'(0^+)-P'(0^-)=\bigl(R_0(1-a_0)-R_0,\,-a_0Z_0\bigr)=-a_0P(0)$. At the end, $Z'=a_+Z$ holds for $a_+=\tfrac{a_0}{a_0L-1}>0$, then $R'=a_+R$ is \eqref{eq:cake-a}, and $|P(L)|=1$ is \eqref{eq:cake-b}. For the existence, \eqref{eq:cake-a} is a quadratic equation for $a_0$. Along its larger root $a_0(L)$, which is real and continuous for $2.06\le L\le2.07$ with $1.85<a_0(L)<1.86$, the difference of the two sides of \eqref{eq:cake-b} is positive at $L=2.06$ and negative at $L=2.07$. So \eqref{eq:cake} has a solution in this interval, and it satisfies $0<a_0<2$ and $a_0L>1$.
\end{proof}
Among Propositions~\ref{prop:disk}--\ref{prop:c3}, only the cake stand has no reflection symmetry in a horizontal plane. The balancing lifts the plate above the equatorial plane, to the side opposite the boundary circle. The surface combines a disk end, as in \cite{Kapouleas_Zou_2021}, with a horn that meets $\partial\ball$ orthogonally, as in \cite{Carlotto_Schulz_Wiygul_2025}. We found no other solution of \eqref{eq:cake}.

\begin{proposition}[The wheel]\label{prop:wheel}
The equation
\begin{equation}\label{eq:wheel}
L^2=L\cosh L-\sinh L
\end{equation}
has exactly one positive root, $L=2.0301901$. Put $s=\tfrac{L\sinh L-\cosh L}{L^2}$. Let $M$ be the part $|t|\le L$ of $\sphere$ with the circle $t_0=0$. Then
\begin{equation*}
P=\Bigl(\cosh t-s\sinh|t|,\ \tfrac{\sqrt{L^2-1}}{L^2}\,t\Bigr),\qquad|t|\le L,
\end{equation*}
is a free curve minimal surface. It consists of two horns glued along the equator, each ending on a boundary circle of radius $a_\pm=\tfrac1L$. The weight of the equator is $a_0=2s$, and
\begin{equation*}
\bar\sigma=2\pi\Bigl(2s+\frac2L\Bigr)=\frac{4\pi}{L^2}\bigl(L+L\sinh L-\cosh L\bigr)\approx17.542995 .
\end{equation*}
\end{proposition}
\begin{proof}
The function $h(L)=L\cosh L-\sinh L-L^2$ has $h(0)=0$ and $h'(L)=L(\sinh L-2)$. So $h$ decreases on $[0,\operatorname{arsinh}2]$ and then increases to infinity. Since $h(0)=0$, it has exactly one positive root, and this root is larger than $\operatorname{arsinh}2>1$. We check the conditions of Lemma~\ref{lem:profile}. Condition (iv) is vacuous since both ends are finite. Since $(L\cosh L-\sinh L)^2-(L\sinh L-\cosh L)^2=L^2-1$ and $L\cosh L-\sinh L=L^2$, we get $1-s^2=\tfrac{L^2-1}{L^4}$, and $|s|<1$ since $L>1$. So each horn, with $A=1$ and $Z'^2=\tfrac{L^2-1}{L^4}=1-s^2$, satisfies (i). At the equator $P(0)=(1,0)$ and $P'(0^+)-P'(0^-)=(-2s,0)=-2sP(0)$. At $t=L$, using $L^2\cosh L=L\cosh^2L-\sinh L\cosh L$ and $L^2\sinh L=L\sinh L\cosh L-\sinh^2L$,
\begin{equation*}
R(L)=\cosh L-s\sinh L=\frac1L,\qquad R'(L)=\sinh L-s\cosh L=\frac1{L^2}.
\end{equation*}
So $|P(L)|^2=\tfrac1{L^2}+\tfrac{L^2-1}{L^2}=1$ and $P'(L)=\tfrac1LP(L)$. The end $t=-L$ follows by symmetry.
\end{proof}
Each horn is the catenoidal annulus $\mathbb K_0$ of \cite[Lemma 3.3]{Carlotto_Schulz_Wiygul_2025}, which meets $\partial\ball$ along the equator and orthogonally along a circle in the upper hemisphere. Its scaling and height parameters in \cite{Carlotto_Schulz_Wiygul_2025} have product $L$, and $\operatorname{Area}(\mathbb K_0)=\tfrac14\bar\sigma$ agrees with \cite[Remark 3.9]{Carlotto_Schulz_Wiygul_2025}. Carlotto, Schulz and Wiygul desingularize the wheel into genus-zero free boundary minimal surfaces \cite[Theorem 6.1]{Carlotto_Schulz_Wiygul_2025}, and the wheel together with the equatorial disk, where three sheets meet along the equator, into surfaces of large genus with three boundary components \cite[Theorem 5.1]{Carlotto_Schulz_Wiygul_2025}.

\begin{proposition}[The critical drum]\label{prop:drum}
Let $t_1=\tfrac12\bigl(1+W(\tfrac1e)\bigr)=0.6392323$, with $W$ the Lambert function, and $\alpha=(\cosh^2t_1+t_1^2)^{-1/2}$. Let $M=\sphere$ with the circles $t=\pm t_1$. Then
\begin{equation*}
P=\alpha\bigl(\cosh t,\ t\bigr)\ \ (|t|\le t_1),\qquad P=\alpha\bigl(e^{t_1-|t|}\cosh t_1,\ \pm t_1\bigr)\ \ (\pm t\ge t_1)
\end{equation*}
is a free curve minimal surface. It consists of a catenoidal band closed by two flat disks. Both circles have the weight $a_1=\tfrac1{t_1}$, and
\begin{equation*}
\bar\sigma=\frac{4\pi}{t_1}=\frac{8\pi}{1+W(\frac1e)}\approx19.658536 .
\end{equation*}
\end{proposition}
\begin{proof}
The band has $A=\alpha$, $s=0$ and $Z'=\alpha$, so it satisfies (i) of Lemma~\ref{lem:profile}. The disks satisfy (i) and (iv), and (iii) is vacuous since $M$ has no boundary. Both formulas give $P(\pm t_1)=\alpha(\cosh t_1,\pm t_1)$, so $P$ is continuous, and the choice of $\alpha$ gives $|P(\pm t_1)|=1$. Since $R$ is even and $Z$ is odd, it suffices to check the jump at $t_1$. There the jump of $Z'$ is $-\alpha$ and the jump of $R'$ is $-\alpha(\cosh t_1+\sinh t_1)$. They equal $-a_1Z(t_1)$ and $-a_1R(t_1)$ if and only if $a_1=\tfrac1{t_1}$ and $a_1=1+\tanh t_1$. Since $1+\tanh t=\tfrac{2e^{2t}}{e^{2t}+1}$, the equation $\tfrac1{t_1}=1+\tanh t_1$ is $2t_1-1=e^{-2t_1}$, that is, $(2t_1-1)e^{2t_1-1}=e^{-1}$. As $ue^u<0$ for $u<0$ and $ue^u$ increases for $u\ge0$, its only solution is $2t_1-1=W(\tfrac1e)$.
\end{proof}
This is the critical drum $\mathcal C(2)$ of \cite{Karpukhin_Noutchegueme_2025}. Its value $\tfrac{4\pi}{t_1}$ is \cite[Theorem A.5]{Karpukhin_Noutchegueme_2025}, where $t_1$ is the positive root of $t=e^{-t}\cosh t$. Since $e^{-t}\cosh t=\tfrac1{1+\tanh t}$, this is the equation $\tfrac1{t_1}=1+\tanh t_1$ of the proof, and the Lambert form makes $t_1$ explicit. The two expressions for $a_1$ in the proof come from the coordinate $z$ and from the coordinates $x$, $y$. Their equality is the balancing of the drum, which makes $\sigma_1$ triple.

\begin{proposition}[The drum stack $\mathcal C(3)$]\label{prop:c3}
Let $t_1=0.9160481$ be the root in $(\tfrac12,\infty)$ of
\begin{equation}\label{eq:c3}
t_1e^{t_1}-\sinh t_1=\sqrt2\,t_1^{3/2} ,
\end{equation}
and put $D=t_1e^{t_1}-\sinh t_1$ and $s_0=\tfrac{t_1e^{t_1}-\cosh t_1}{D}$. Let $M=\sphere$ with the circles $t_0=0$ and $t=\pm t_1$. Then
\begin{equation*}
P=\Bigl(\cosh t-s_0\sinh|t|,\ \tfrac{\sqrt{2t_1-1}}{D}\,t\Bigr)\ \ (|t|\le t_1),
\end{equation*}
continued by flat disks for $|t|\ge t_1$, is a free curve minimal surface. It consists of two catenoidal bands glued along the equator and closed by two flat disks. The weights are $a_0=2s_0=2-\sqrt2\,e^{-t_1}t_1^{-3/2}$ at the equator and $a_1=\tfrac1{t_1}$ at the outer circles, and
\begin{equation*}
\bar\sigma=2\pi(a_0+2a_1)=2\pi\Bigl(2+\frac2{t_1}-\sqrt2\,e^{-t_1}t_1^{-3/2}\Bigr)\approx22.229471 .
\end{equation*}
\end{proposition}
\begin{proof}
Put $h(T)=Te^T-\sinh T-\sqrt2\,T^{3/2}$. Then $h(\tfrac12)=\tfrac12e^{1/2}-\sinh\tfrac12-\tfrac12\approx-0.197<0$, $h(T)\to\infty$, and $h$ is convex on $[\tfrac12,\infty)$, because $h'''(T)=(3+T)e^T-\cosh T+\tfrac{3\sqrt2}8T^{-3/2}>0$ and $h''(\tfrac12)=\tfrac52e^{1/2}-\sinh\tfrac12-\tfrac32\approx2.10>0$. So $h$ has exactly one root $t_1$ in $(\tfrac12,\infty)$. We check the conditions of Lemma~\ref{lem:profile}. The disks satisfy (i) and (iv), and (iii) is vacuous since $M$ has no boundary. The bands have $A=1$, $s=\pm s_0$ and $Z'^2=\tfrac{2t_1-1}{D^2}$, and $D^2(1-s_0^2)=(t_1e^{t_1}-\sinh t_1)^2-(t_1e^{t_1}-\cosh t_1)^2=2t_1e^{t_1}(\cosh t_1-\sinh t_1)-1=2t_1-1>0$. So $|s_0|<1$ and $Z'^2=1-s_0^2$, and the bands satisfy (i). At the equator $P(0)=(1,0)$ and the jump of $P'$ is $-2s_0P(0)$. At $t=t_1$, since $\cosh t_1-\sinh t_1=e^{-t_1}$, the definitions of $D$ and $s_0$ give $DR(t_1)=D\cosh t_1-(t_1e^{t_1}-\cosh t_1)\sinh t_1=t_1e^{t_1}(\cosh t_1-\sinh t_1)=t_1$ and $DR'(t_1^-)=D\sinh t_1-(t_1e^{t_1}-\cosh t_1)\cosh t_1=1-t_1$. The disk has $R'(t_1^+)=-R(t_1)$ and $Z'=0$, so the jump of $P'$ is $-\tfrac1D\bigl(1,\sqrt{2t_1-1}\bigr)=-\tfrac1{t_1}P(t_1)$. Also $|P(t_1)|^2=\tfrac{t_1^2+(2t_1-1)t_1^2}{D^2}=\tfrac{2t_1^3}{D^2}=1$ by the equation for $t_1$. The circle $t=-t_1$ follows by symmetry, since $R$ is even and $Z$ is odd. Finally $2-a_0=2-2s_0=\tfrac2D(\cosh t_1-\sinh t_1)=\tfrac{2e^{-t_1}}D=\sqrt2\,e^{-t_1}t_1^{-3/2}$.
\end{proof}
This is the drum stack $\mathcal C(3)$ of \cite{Karpukhin_Noutchegueme_2025}. Its outer circles lie at polar angle $47.63^\circ$ on $\partial\ball$, inside the interval $(\tfrac\pi4,\tfrac{2\pi}7)$ in which Kapouleas and Zou locate them (proof of \cite[Corollary 2.24]{Kapouleas_Zou_2021}).

The remaining networks of Table~\ref{tab:limitnetworks} are the drum and the drum stack with one or both disks replaced by horns, and the stack $\mathcal C(4)$. For these we solved the conditions of Lemma~\ref{lem:profile} by Newton's method, with residuals below $10^{-14}$, and the bottom row of Figure~\ref{fig:fcms} shows their meridian sections. For each type of network with horns, we found a single solution with positive weights. The stacks closed by disks are unique for each number of circles \cite[Theorem 1.10]{Karpukhin_Noutchegueme_2025}, \cite[Corollary 2.24]{Kapouleas_Zou_2021}. For them $\sigma=1$ is the first nonzero eigenvalue of \eqref{eq:transmission}, and $\mathcal C(N)$ maximizes $\bar\sigma_1$ among rotationally symmetric metrics on $\sphere$ with $N$ parallels \cite[Theorem 1.10]{Karpukhin_Noutchegueme_2025}. The networks with horns are not closed, so \cite{Karpukhin_Noutchegueme_2025} does not apply to them, and the value $\bar\sigma$ of a network is its $\bar\sigma_1$ only if $\sigma=1$ is the first nonzero eigenvalue. We checked this for the ten networks of Table~\ref{tab:limitnetworks} by solving \eqref{eq:cyl} mode by mode. In the mode $e^{im\phi}$ the Rayleigh quotient is $\int(f'^2+m^2f^2)\,dt$ divided by $\sum_ia_if(t_i)^2+\sum_\pm a_\pm f(t_\pm)^2$, so the first eigenvalue of this mode increases with $m\ge1$. It is at least $2$ for $m=2$ in every case, so only the modes with $m\le1$ matter. For every network except the doubly covered disk, $\sigma=1$ is the first nonzero eigenvalue, with multiplicity three: $z$ gives it in the rotationally symmetric mode, and $x$ and $y$ in the modes $e^{\pm i\phi}$. The next eigenvalue is at least $1.545$, the value for the wheel. The doubly covered disk lies in the plane $z=0$, and its nonzero eigenvalues are $\sigma=m$ in the modes $e^{\pm im\phi}$, so there $\sigma=1$ has multiplicity two.

\section{Layer potential formulation}\label{sec:method}
We pose the extremal problem on the sphere with $b$ round caps removed and describe its symmetry classes (Sections~\ref{ssec:uniformization}--\ref{ssec:classes}). We then write the Steklov problem as a boundary integral equation and in multipole form, with its monopole and dipole approximations (Sections~\ref{ssec:layer}--\ref{ssec:approx}), and give the derivatives and the certificate that turn a computed critical point into a free boundary minimal surface (Sections~\ref{ssec:derivatives}--\ref{ssec:checks}). Most of these results adapt \cite{Chen_Osting_2026,kim2025layer,Gemmrich_Nigam_Steinbach_2008} to the sphere, and we refer there for proofs. The multipole form of Section~\ref{ssec:multipole} is new and is proved here. We follow the notation of \cite{Chen_Osting_2026}, with spherical centers $c_\ell\in\sphere$ and angular radii $r_\ell$ in place of planar ones.

\subsection{Uniformization and the weighted problem}\label{ssec:uniformization}
By Koebe uniformization \cite[Theorem 0.2]{HeSchramm1993}, $\Sigma$ is conformally $\Omega:=\sphere\setminus\bigcup_{\ell\in[b]}\overline{C_\ell}$ with $C_\ell=\{x\in\sphere\colon x\cdot c_\ell>\cos r_\ell\}$ disjoint round caps, unique up to M\"obius transformations. With a frame $(e_{1,\ell},e_{2,\ell},c_\ell)$, the boundary circle $\Gamma_\ell$ is
\begin{equation}\label{eq:circle}
\xi_\ell(\theta)=\cos r_\ell\,c_\ell+\sin r_\ell\bigl(\cos\theta\,e_{1,\ell}+\sin\theta\,e_{2,\ell}\bigr),\qquad ds=\sin r_\ell\,d\theta,\qquad \nu=-\partial_{r_\ell}\xi_\ell,
\end{equation}
with $\nu$ the outward conormal of $\Omega$ and $\cot r_\ell$ the geodesic curvature of $\Gamma_\ell$. Since the Steklov spectrum depends on the conformal factor only through its boundary trace $\rho=e^\omega$ \cite[\S2.1]{Chen_Osting_2026}, the extremal problem is
\begin{subequations}\label{eq:problem}
\begin{gather}
\max_{\mathbf p}\ \bar\sigma_1(\mathbf p)\quad\text{subject to}\quad
\begin{cases}
\Delta_{\sphere}u=0 & \text{in }\Omega,\\
\partial_\nu u=\sigma\rho\,u & \text{on }\partial\Omega,
\end{cases}\label{eq:weighted}\\
c_1=e_3,\qquad c_2=-e_3,\qquad r_1=r_2,\qquad c_3\cdot e_2=0,\qquad c_3\cdot e_1>0\label{eq:normalization}
\end{gather}
\end{subequations}
over $\mathbf p=(\mathbf c,\mathbf r,\boldsymbol\alpha,\boldsymbol\beta)$. We normalize the eigenfunctions by $\int_{\partial\Omega}\rho\,u_ku_{k'}\,ds=\delta_{kk'}$. On each circle $\Gamma_\ell$, the log-weight $\omega_\ell=\omega\circ\xi_\ell$ is the truncated Fourier series
\[
\omega_\ell(\theta)=\alpha_{\ell,0}+\sum_{n=1}^{n_{\max}}\bigl(\alpha_{\ell,n}\cos n\theta+\beta_{\ell,n}\sin n\theta\bigr)
\]
in the angle $\theta$ of \eqref{eq:circle}, with coefficients $\boldsymbol\alpha=(\alpha_{\ell,n})$ and $\boldsymbol\beta=(\beta_{\ell,n})$. Here $\bar\sigma_k=L_\rho\sigma_k$ with $L_\rho=\int_{\partial\Omega}\rho\,ds$, and the caps are disjoint, $\arccos(c_\ell\cdot c_{\ell'})>r_\ell+r_{\ell'}$. The moduli space has dimension $3b-6$ ($b\ge3$). $\bar\sigma_k$ is invariant under the M\"obius group acting on $(\Omega,\rho)$ by $(\mathcal T\Omega,(\tfrac{\rho}{|d\mathcal T|})\circ\mathcal T^{-1})$. The normal form \eqref{eq:normalization} removes this invariance by fixing the first two centers and the direction of the third. This is always possible: two disjoint circles are M\"obius equivalent to concentric ones, a dilation ($z\mapsto\lambda z$ in a stereographic chart centered at its fixed points) equalizes the radii, and a rotation about $e_3$ fixes $c_3$.

\emph{Example.} With \eqref{eq:normalization}, the case $b=2$ is the symmetric annulus $\{|x\cdot e_3|<\cos r\}$, with two caps of equal radius $r$ at the poles. The critical catenoid is the member $r^\star\approx0.59$ of this one-parameter family, with constant weight $\rho\equiv1$. Its radius solves $\cos r^\star\,\log\cot\tfrac{r^\star}{2}=1$, which is the critical catenoid's equation $t\tanh t=1$ for $t=\log\cot\tfrac{r^\star}{2}$. At this radius the separated branches of Section~\ref{ssec:validation} coalesce into the triple $\bar\sigma_1^\star=4\pi\cos r^\star\approx10.47$.

\subsection{Symmetry classes}\label{ssec:classes}
Most of our computations impose a symmetry, which cuts down the number of parameters and lets us reach surfaces that are saddle points of the full problem. No symmetry is lost by working on the sphere. A finite group of conformal and anticonformal maps of $\sphere$ is conjugate by a M\"obius transformation into $\mathrm O(3)$, and its orientation-preserving part into $\mathrm{SO}(3)$. So every finite symmetry group $\mathcal G$ of a surface acts on $\Omega$ by rotations and reflections.

Under such a group the caps fall into $\mathcal G$-orbits. Caps in one orbit have the same radius, and $\omega$ has one constant level per orbit. Adding a constant to $\omega$ multiplies $\rho=e^\omega$ by a constant and leaves $\bar\sigma_1$ unchanged, so we measure the levels from that of the first orbit. The weight on a circle must also respect the stabilizer of its center, the elements of $\mathcal G$ that fix $c_\ell$. These are rotations about the axis through $c_\ell$, possibly together with reflections in planes containing that axis. We call such a plane, whose reflection lies in $\mathcal G$, a \emph{mirror}. If the stabilizer is $C_{qv}$, generated by the rotation by $\tfrac{2\pi}{q}$ and $q$ mirror planes, we choose the frame of \eqref{eq:circle} with $e_{1,\ell}$ in a mirror plane through $c_\ell$, and then $\omega_\ell(\theta)=\sum_n\alpha_{\ell,n}\cos(qn\theta)$. If it is the cyclic group $C_q$, without mirrors, the terms $\sin(qn\theta)$ appear as well, and $q=1$ gives the full series of Section~\ref{ssec:uniformization}. A \emph{class} is fixed by the group, the orbits and the freedom of each center within its orbit. Its geometric parameters are one radius per orbit, the free positions of the centers, and one level per orbit after the first.

Table~\ref{tab:classes} lists the classes we compute. From the most rigid to the most flexible, they are the polyhedral configurations, the equatorial family of caps on one great circle, and the ring families: pyramids, bipyramids, prisms, antiprisms and stacked rings. A \emph{ring} of $j$ caps puts the centers at a common polar angle $\beta$, with azimuths $\phi_0+\tfrac{2\pi k}{j}$, $k\in[j]$. Its latitude is free when $\beta$ is a parameter of the class. Two rings at $\beta$ and $\pi-\beta$ form a prism pair when aligned and an antiprism pair when rotated by $\tfrac{\pi}{j}$. Stacks are several rings with alternating offsets. Any of these may also have a cap at one or both poles. The polyhedral classes start from the vertices of the Platonic and Archimedean solids and of the geodesic icosahedra of frequency $\nu$, with $b=10\nu^2+2$. The geodesic icosahedron, the icosphere of mesh generation, is made by dividing each edge of the icosahedron into $\nu$ equal parts and each face into $\nu^2$ triangles, and projecting to $\sphere$. Its orbits either lie on the arcs where the mirrors of $I_h$ meet $\sphere$, free to slide along them, or are generic. At $b=72$ a generic orbit of the rotation group $I$ has two free position parameters. Section~\ref{ssec:best} describes the Fekete classes and the unrestricted problem, which do not presuppose a family.

\subsection{Layer potentials on the sphere}\label{ssec:layer}
Since $\sphere$ is closed, the Laplace--Beltrami operator has no fundamental solution, and one works instead with the modified Green's function
\begin{equation}\label{eq:green}
G(x,y):=-\frac1{2\pi}\log|x-y|,\qquad -\Delta_{\sphere,x}G(x,y)=\delta_y(x)-\frac1{4\pi},
\end{equation}
with $|x-y|$ the chordal distance in $\R^3$, as in the boundary integral methods on the sphere of Gemmrich, Nigam and Steinbach \cite{Gemmrich_Nigam_Steinbach_2008} and Kropinski and Nigam \cite{Kropinski_Nigam_2014}. The constant $\tfrac1{4\pi}$ is subtracted so that the right-hand side has mean zero. With $ds$ the round arclength, the single layer, the adjoint double layer and the modified single layer are
\begin{equation}\label{eq:layer}
\mathcal S[\varphi](x)=\int_{\partial\Omega}G(x,y)\varphi(y)\,ds(y),\quad
\mathcal D^*[\varphi](x)=\frac1{2\pi}\int_{\partial\Omega}\frac{\nu(x)\cdot y}{|x-y|^2}\varphi(y)\,ds(y),\quad
\mathcal S_0:=\mathcal S(I-M)+M,
\end{equation}
with $M\varphi:=\frac1{|\partial\Omega|}\int_{\partial\Omega}\varphi\,ds$ the mean of $\varphi$ with respect to $ds$, a constant function, and $|\partial\Omega|=\sum_\ell2\pi\sin r_\ell$. By \eqref{eq:green}, $\mathcal S_0[\varphi]$ is harmonic in $\Omega$ for every $\varphi$. On a circle, by \eqref{eq:circle},
\begin{equation}\label{eq:self-blocks}
G(x,y)=-\frac1{2\pi}\log\sin r_\ell-\frac1{4\pi}\log\Bigl(4\sin^2\frac{\theta-\theta'}2\Bigr),\qquad \partial_{\nu(x)}G(x,y)=\frac{\cot r_\ell}{4\pi},\qquad x,y\in\Gamma_\ell .
\end{equation}
The first is the self-interaction of a planar circle of radius $\sin r_\ell$, and the adjoint double-layer kernel is constant on $\Gamma_\ell$.

\begin{theorem}[Boundary integral formulation on the punctured sphere]\label{thm:bie}
$(\sigma_k,u_k)$ is an eigenpair of \eqref{eq:weighted} if and only if a nonzero density $\varphi_k$ satisfies $\mathsf A\varphi_k=\sigma_k\mathsf B\varphi_k$ with $\mathsf A:=(\mathcal D^*+\tfrac12I)(I-M)$, $\mathsf B:=\rho\,\mathcal S_0$ built from the operators \eqref{eq:layer}, and then $u_k=\mathcal S_0[\varphi_k]$.
\end{theorem}
This is the sphere version of \cite[Theorem 3]{kim2025layer}, and the proof there carries over. The jump relations are local. The representation of harmonic functions reduces to the planar case by stereographic projection from a point inside one cap, because the terms of \eqref{eq:green} other than the planar logarithm are constant on densities of zero mean.

\subsection{Multipole analysis}\label{ssec:multipole}
For round caps, the layer potentials of Theorem~\ref{thm:bie} take an explicit form in the Fourier modes on each circle, which we write as a multipole expansion. Seen from outside, a cap acts like a point charge at its center, whose strength is the flux through its circle, together with higher multipoles at the same point: a dipole, a quadrupole, and so on. Seen from one cap, the field of all the other caps is harmonic near it and has a Taylor expansion at its center. Proposition~\ref{prop:multipole} shows how the Steklov condition on the circle couples the two expansions mode by mode.

For each $\ell$, let $\zeta_\ell$ be the stereographic coordinate on $\sphere$ that is centered at $c_\ell$, isometric there, and aligned with the frame of \eqref{eq:circle}. In this coordinate $C_\ell=\{|\zeta_\ell|<R_\ell\}$ with $R_\ell=2\tan\tfrac{r_\ell}2$, the circle is $\zeta_\ell=R_\ell e^{i\theta}$ with the angle $\theta$ of \eqref{eq:circle}, and the antipode of $c_\ell$ is $\zeta_\ell=\infty$. Let $w_\ell=\tfrac1{L_\rho}\int_{\Gamma_\ell}\rho\,ds$ be the fraction of the total weight on $\Gamma_\ell$, so that $\sum_\ell w_\ell=1$, and expand the weight on $\Gamma_\ell$ in the angle,
\begin{equation}\label{eq:weight-modes}
\rho\,ds=\frac{L_\rho}{2\pi}\sum_{k\in\mathbb Z}w_{\ell,k}\,e^{ik\theta}\,d\theta,\qquad w_{\ell,0}=w_\ell,\qquad w_{\ell,-k}=\overline{w_{\ell,k}} .
\end{equation}
By Theorem~\ref{thm:bie}, every eigenfunction of \eqref{eq:weighted} is a single-layer potential $u=\mathcal S_0[\varphi]$ of a density $\varphi$ on $\partial\Omega$, so it suffices to expand these potentials.

\begin{proposition}[Multipole form]\label{prop:multipole}
Let $\varphi$ be a smooth density on $\partial\Omega$ and $u=\mathcal S_0[\varphi]$. Expand $\varphi|_{\Gamma_\ell}=\sum_{n\in\mathbb Z}\varphi_{\ell,n}e^{in\theta}$.
\begin{enumerate}[label=(\alph*),leftmargin=*]
\item (Charges and multipoles) The function $u$ has an expansion of the form
\begin{equation}\label{eq:mp-rep}
u=u_\infty+\sum_{\ell\in[b]}\Bigl(q_\ell\,G(\cdot,c_\ell)+\operatorname{Re}\sum_{n\ge1}a_{\ell,n}\,\zeta_\ell^{-n}\Bigr),
\end{equation}
where $u_\infty\in\mathbb R$, $q_\ell=\int_{\Gamma_\ell}(\varphi-M\varphi)\,ds$ and $a_{\ell,n}=\tfrac{\sin r_\ell}{n}R_\ell^{\,n}\,\overline{\varphi_{\ell,n}}$. The charges satisfy $\sum_\ell q_\ell=0$, and $q_\ell=\int_{\Gamma_\ell}\partial_\nu u\,ds$ is the flux of $u$ into $C_\ell$.
\item (The field at one cap) Let $d_m>R_m$ be the distance from $0$ to the images $\zeta_m(C_\ell)$, $\ell\ne m$. Then
\begin{equation}\label{eq:mp-local}
u=-\frac{q_m}{2\pi}\log|\zeta_m|+\operatorname{Re}\sum_{n\ge1}a_{m,n}\,\zeta_m^{-n}+\operatorname{Re}\sum_{n\ge0}h_{m,n}\,\zeta_m^{\,n},\qquad R_m\le|\zeta_m|<d_m,
\end{equation}
where the last sum converges on $|\zeta_m|<d_m$. Its value at $c_m$ is
\begin{equation}\label{eq:mp-e0}
h_{m,0}=u_\infty+\sum_{\ell\ne m}\Bigl(q_\ell\,G(c_m,c_\ell)+\operatorname{Re}\sum_{n\ge1}a_{\ell,n}\,\zeta_\ell(c_m)^{-n}\Bigr).
\end{equation}
\item (The Steklov condition) Expand $u|_{\Gamma_m}=\sum_{n\in\mathbb Z}U_{m,n}e^{in\theta}$, where $\zeta_m=R_me^{i\theta}$ on $\Gamma_m$. Then $U_{m,-n}=\overline{U_{m,n}}$ and
\begin{subequations}\label{eq:mp-modes}
\begin{align}
U_{m,0}&=\frac{q_m}{2\pi}\log\frac1{R_m}+h_{m,0},\label{eq:mp-mode0}\\
U_{m,n}&=\tfrac12\bigl(h_{m,n}R_m^n+\overline{a_{m,n}}R_m^{-n}\bigr)\qquad(n\ge1).\label{eq:mp-moden}
\end{align}
\end{subequations}
The pair $(\sigma,u)$ satisfies \eqref{eq:weighted} if and only if, with $\bar\sigma=L_\rho\sigma$, for all $m\in[b]$,
\begin{equation}\label{eq:mp-steklov}
q_m=\bar\sigma\sum_{k\in\mathbb Z}w_{m,-k}\,U_{m,k},\qquad \tfrac n2\bigl(\overline{a_{m,n}}R_m^{-n}-h_{m,n}R_m^n\bigr)=\frac{\bar\sigma}{2\pi}\sum_{k\in\mathbb Z}w_{m,n-k}\,U_{m,k}\quad(n\ge1).
\end{equation}
\end{enumerate}
\end{proposition}
Part (a) says that the single-layer potential of a density on the circles is the field of point sources at the centers of the caps. The charge $q_\ell$ and the multipole moment $a_{\ell,n}$ of order $n$ of the cap $\ell$ are the Fourier coefficients of the density on $\Gamma_\ell$, rescaled. In part (b), \eqref{eq:mp-local} is the Laurent expansion of $u$ in the annulus around one cap, the form given by the logarithmic conjugation theorem \cite{Axler_1986}. The logarithm and the negative powers are the charge and the multipoles of the cap $m$ itself, from \eqref{eq:mp-rep}. The nonnegative powers are the incoming field of the other caps, which is regular at $c_m$. The coefficient $h_{m,0}$ is the value of this field at $c_m$, given by \eqref{eq:mp-e0}, and $h_{m,1}$ encodes its gradient there. In part (c), the mode $n\ge1$ of $u$ on $\Gamma_m$ has an incoming part $h_{m,n}R_m^n$ and an outgoing part $\overline{a_{m,n}}R_m^{-n}$. The trace \eqref{eq:mp-moden} is half their sum, and the flux is $\tfrac n2$ times their difference. In the mode $0$, the trace \eqref{eq:mp-mode0} is the mean of $u$ on $\Gamma_m$, and the flux is $\tfrac{q_m}{2\pi}$. So \eqref{eq:mp-steklov} is $\partial_\nu u=\sigma\rho\,u$ written mode by mode. By \eqref{eq:mp-rep}--\eqref{eq:mp-modes}, $q$, $a$, $h$ and $U$ depend linearly on $\varphi$, so \eqref{eq:mp-steklov} is the generalized eigenvalue problem $\mathsf A\varphi=\sigma\mathsf B\varphi$ of Theorem~\ref{thm:bie}, written in the modes $\varphi_{\ell,n}$, with the flux on the left and the weighted trace on the right.
\begin{proof}[Proof of Proposition~\ref{prop:multipole}]
(a) In the coordinate $\zeta_\ell$ the chordal distance is \cite[Chapter~1, Section~2.4]{Ahlfors_1979}
\begin{equation}\label{eq:chordal}
|x-y|=\frac{|\zeta_\ell(x)-\zeta_\ell(y)|}{\bigl(1+\tfrac14|\zeta_\ell(x)|^2\bigr)^{1/2}\bigl(1+\tfrac14|\zeta_\ell(y)|^2\bigr)^{1/2}} .
\end{equation}
On $\Gamma_\ell$, $\zeta_\ell=R_\ell e^{i\theta}$ and $ds=\sin r_\ell\,d\theta$. For $|\zeta_\ell(x)|>R_\ell$, $\log|\zeta_\ell(x)-R_\ell e^{i\theta}|=\log|\zeta_\ell(x)|-\operatorname{Re}\sum_{n\ge1}\tfrac1nR_\ell^ne^{in\theta}\zeta_\ell(x)^{-n}$, and in \eqref{eq:chordal} the factor of $x$ does not depend on $\theta$ and the factor of $y$ is constant on $\Gamma_\ell$. So a density $e^{in\theta}$ on $\Gamma_\ell$ with $n\ge1$ has the single-layer potential $\tfrac{\sin r_\ell}{2n}R_\ell^n\,\overline{\zeta_\ell(x)}^{\,-n}$ outside $C_\ell$, and the density $1$ on $\Gamma_\ell$ has the potential $2\pi\sin r_\ell\,G(x,c_\ell)$ plus a constant, since the mean of $\log|x-y|$ over $y\in\Gamma_\ell$ is $\log|x-c_\ell|$ plus a constant. Since $\mathcal S_0[\varphi]=\mathcal S[\varphi-M\varphi]+M\varphi$, summing over the modes, with $\varphi_{\ell,-n}=\overline{\varphi_{\ell,n}}$, gives \eqref{eq:mp-rep}. The series of the cap $\ell$ converges for $|\zeta_\ell|\ge R_\ell$ because $\varphi$ is smooth, and $\sum_\ell q_\ell=\int_{\partial\Omega}(\varphi-M\varphi)\,ds=0$. By \eqref{eq:green}, the flux of $G(\cdot,c_\ell)$ into $C_m$ is $\delta_{\ell m}-\tfrac{|C_m|}{4\pi}$, and the multipoles have no flux, so the flux of $u$ into $C_m$ is $q_m-\tfrac{|C_m|}{4\pi}\sum_\ell q_\ell=q_m$.

(b) Write $u=-\tfrac{q_m}{2\pi}\log|\zeta_m|+\operatorname{Re}\sum_{n\ge1}a_{m,n}\zeta_m^{-n}+v_m$, so that $v_m$ is the field of the other caps. By \eqref{eq:mp-rep}, $v_m$ is the sum of $u_\infty$, of the multipole series of the caps $\ell\ne m$, and of $\sum_\ell q_\ell\,G(\cdot,c_\ell)+\tfrac{q_m}{2\pi}\log|\zeta_m|$. The multipole series of the cap $\ell\ne m$ is harmonic outside $C_\ell$. By \eqref{eq:chordal} in the coordinate $\zeta_m$, each $G(\cdot,c_\ell)$ is $-\tfrac1{2\pi}\log|\zeta_m-\zeta_m(c_\ell)|$, plus a constant, plus $\tfrac1{4\pi}\log\bigl(1+\tfrac14|\zeta_m|^2\bigr)$. For $c_\ell=-c_m$ the first term is absent. The last term is the same for every $\ell$, and its total coefficient $\sum_\ell q_\ell$ vanishes. So the logarithmic part of $v_m$ is harmonic away from the centers $c_\ell$, $\ell\ne m$, and it tends to $\sum_{\ell\ne m}q_\ell\,G(c_m,c_\ell)$ at $c_m$. Hence $v_m$ is harmonic on the disk $|\zeta_m|<d_m$, and its value at $c_m$ is \eqref{eq:mp-e0}. Harmonicity is conformally invariant in two dimensions, so $v_m$ is harmonic in $\zeta_m$ in the plane sense. On the disk it is therefore the real part of a holomorphic function, whose Taylor series at $0$ gives $v_m=\operatorname{Re}\sum_{n\ge0}h_{m,n}\zeta_m^{\,n}$. This is \eqref{eq:mp-local}.

(c) On $\Gamma_m$ put $\zeta_m=R_me^{i\theta}$ in \eqref{eq:mp-local}. Since $\operatorname{Re}(h\zeta^n)=\tfrac12(h\zeta^n+\bar h\bar\zeta^n)$, the coefficient of $e^{in\theta}$, $n\ge1$, is $\tfrac12h_{m,n}R_m^n$ from the field of the other caps and $\tfrac12\overline{a_{m,n}}R_m^{-n}$ from the multipole. This gives \eqref{eq:mp-modes}. The conormal $\nu$ points into $C_m$, toward smaller $|\zeta_m|$, and the flux form $\partial_\nu u\,ds$ is conformally invariant. So $\partial_\nu u\,\tfrac{ds}{d\theta}=-|\zeta_m|\,\partial_{|\zeta_m|}u$ at $|\zeta_m|=R_m$. Its mode $0$ is $\tfrac{q_m}{2\pi}$, and its mode $n\ge1$ is $\tfrac n2\bigl(\overline{a_{m,n}}R_m^{-n}-h_{m,n}R_m^n\bigr)$. By \eqref{eq:weight-modes}, the condition $\partial_\nu u\,ds=\sigma\rho\,u\,ds$ on $\Gamma_m$ is therefore \eqref{eq:mp-steklov} in the modes $n\ge0$, and its modes $n<0$ are the complex conjugates.
\end{proof}

\begin{remark}[Particular solutions]\label{rem:mps}
The method of particular solutions of \cite{Kao_2020} works on a planar disk with circular holes. It approximates the eigenfunction by the real and imaginary parts of $1$, $z^n$, $\log|z-z_i|$ and $(z-z_i)^{-n}$, with $z_i$ a point in the $i$-th hole, and imposes the Steklov condition by least squares at collocation points. The stereographic coordinate $z$ centered at $-c_1$ maps $\Omega$ onto such a domain, with $\Gamma_1$ as the outer circle. In it the multipoles of the cap $1$ are multiples of $z^n$, and those of a cap $\ell\ne1$ are polynomials in $(z-z_\ell)^{-1}$, with $z_\ell$ the image of $c_\ell$. By \eqref{eq:chordal} and $\sum_\ell q_\ell=0$, the charges give $\sum_{\ell\ne1}-\tfrac{q_\ell}{2\pi}\log|z-z_\ell|$ plus a constant. So \eqref{eq:mp-rep} spans the same functions. Here, however, the coefficients are the Fourier modes of the density, and Proposition~\ref{prop:multipole}(c) imposes the Steklov condition exactly, mode by mode.
\end{remark}

\begin{corollary}[Constant weights]\label{cor:const}
Let $\rho$ be constant on each circle, i.e., $\rho|_{\Gamma_\ell}=\tfrac{L_\rho\,w_\ell}{2\pi\sin r_\ell}$, and let $u$ be as in Proposition~\ref{prop:multipole}. Then $(\sigma,u)$ satisfies \eqref{eq:weighted} if and only if, with $\bar\sigma=L_\rho\sigma$ and $\varepsilon_m=\tfrac{\bar\sigma\,w_m}{2\pi}$,
\begin{equation}\label{eq:mp-const}
q_m=\bar\sigma\,w_m\,U_{m,0},\qquad(n-\varepsilon_m)\,\overline{a_{m,n}}=(n+\varepsilon_m)\,R_m^{2n}\,h_{m,n}\qquad\text{for all }m\in[b]\text{ and }n\ge1 .
\end{equation}
In particular, if $\varepsilon_m\ne n$, then $\overline{a_{m,n}}=\kappa_{m,n}\,R_m^{2n}\,h_{m,n}$ with $\kappa_{m,n}=\tfrac{n+\varepsilon_m}{n-\varepsilon_m}$.
\end{corollary}
\begin{proof}
Now $w_{m,k}=0$ for $k\ne0$, so \eqref{eq:mp-steklov} reads $q_m=\bar\sigma w_mU_{m,0}$ and $\tfrac n2\bigl(\overline{a_{m,n}}R_m^{-n}-h_{m,n}R_m^n\bigr)=\varepsilon_mU_{m,n}$. Insert \eqref{eq:mp-moden}.
\end{proof}
Thus the flux of a cap is proportional to the mean of $u$ on its circle, and each multipole of a cap is proportional to the Taylor coefficient of the same order of the field of the other caps. Only the centers $c_\ell$, the radii $R_\ell$ and the fractions $w_\ell$ enter. When $w_m$ is small, as for many caps, $\varepsilon_m$ is small and $\kappa_{m,n}$ is close to $1$, the response of an insulating hole.

\subsection{Monopole and dipole approximations}\label{ssec:approx}
Let the weight be constant on each circle, as in Corollary~\ref{cor:const}, which costs little when it is nearly constant. Let $\tilde\rho$ be the mean of $\rho$ on each circle. If $(1-\delta)\tilde\rho\le\rho\le(1+\delta)\tilde\rho$, then $\tfrac{1}{1+\delta}\bar\sigma_k(\tilde\rho)\le\bar\sigma_k(\rho)\le\tfrac{1}{1-\delta}\bar\sigma_k(\tilde\rho)$ for every $k$, by the min--max principle, since the two weights have the same total. Since the multipole of order $n$ carries the factor $R_m^{2n}$, the high orders matter little for small caps, and we keep only the first of them.

In the \emph{monopole approximation}, every cap has a charge and no multipoles: $a_{\ell,n}=0$ for all $n\ge1$. Let $\bar u=(U_{m,0})_m$ be the vector of the means of $u$ on the circles and $q=(q_\ell)_\ell$ the vector of the charges. By \eqref{eq:mp-mode0} and \eqref{eq:mp-e0}, $\bar u=\mathsf Gq+u_\infty\mathbf 1$, where $\mathbf 1$ is the vector of ones and
\begin{equation}\label{eq:G-matrix}
\mathsf G_{m\ell}=G(c_m,c_\ell)\quad(m\ne\ell),\qquad\mathsf G_{\ell\ell}=\frac1{2\pi}\log\frac1{R_\ell} .
\end{equation}
The entry $\mathsf G_{m\ell}$ is the potential at $c_m$ of a unit charge at $c_\ell$, and $\mathsf G_{\ell\ell}$ is the mean of $-\tfrac1{2\pi}\log|\zeta_\ell|$ on $\Gamma_\ell$. The term $\tfrac1{4\pi}\log(1+\tfrac14R_m^2)$ of the potentials drops out of these means, as in the proof of Proposition~\ref{prop:multipole}(b). By \eqref{eq:mp-const}, the charges are proportional to the means, $q=\bar\sigma\operatorname{diag}(w)\,\bar u$, and the constant $u_\infty$ is fixed by the condition that the charges add up to zero. Together these give the eigenvalue problem
\begin{equation}\label{eq:mp-model}
\mathsf P\,\mathsf G\operatorname{diag}(w)\,\bar u=\frac1{\bar\sigma}\,\bar u,\qquad\mathsf P=I-\mathbf 1w^{\mathsf T},
\end{equation}
for $b\times b$ matrices, where $\mathsf P$ removes the mean weighted by $w$.

In the \emph{dipole approximation}, every cap has a charge and a dipole: $a_{\ell,n}=0$ for all $n\ge2$. Let $a=(a_{\ell,1})_\ell$ be the vector of the dipoles and $h=(h_{\ell,1})_\ell$ the vector of the gradients of the incoming fields at the centers. Then $\bar u=\mathsf Gq+u_\infty\mathbf 1+\mathsf Va$ and $h=\mathsf D_1q+\mathsf D_2a$ for real linear maps that depend only on the centers and the frames of \eqref{eq:circle}. Here $\mathsf V$ gives the values at the centers of the fields of the dipoles, and $\mathsf D_1$ and $\mathsf D_2$ give the gradients at the centers of the fields of the charges and of the dipoles. By Corollary~\ref{cor:const} with $n=1$, $a=\mathsf Kh$, where $\mathsf K$ records the response of each cap, $a_{\ell,1}=\kappa_{\ell,1}R_\ell^2\,\overline{h_{\ell,1}}$. Eliminating $a$ and $h$ gives \eqref{eq:mp-model} with $\mathsf G$ replaced by
\begin{equation}\label{eq:geff}
\mathsf G_{\rm eff}=\mathsf G+\mathsf V\,(I-\mathsf K\mathsf D_2)^{-1}\mathsf K\,\mathsf D_1 .
\end{equation}
The matrix $\mathsf G_{\rm eff}$ maps the charges to the means of $u$ on the circles, taking into account the dipoles that the charges induce. It depends on $\bar\sigma$ through the factors $\kappa_{\ell,1}$, and we solve the eigenvalue problem by iteration.

The reciprocals of the largest eigenvalues of the two problems are
\begin{equation}\label{eq:models}
\bar\sigma^{(0)}_1=\frac1{\lambda_{\max}\bigl(\mathsf P\,\mathsf G\operatorname{diag}(w)\bigr)},\qquad
\bar\sigma^{(1)}_1=\frac1{\lambda_{\max}\bigl(\mathsf P\,\mathsf G_{\rm eff}\operatorname{diag}(w)\bigr)}.
\end{equation}
Here $\mathsf G_{\rm eff}$ is evaluated at $\bar\sigma=\bar\sigma^{(1)}_1$. Section~\ref{ssec:large} uses them to approximate $\bar\sigma_1$ at large $b$. The same elimination, carried to multipoles of any order, gives further approximations. With all the multipoles it recovers the problem \eqref{eq:weighted} itself, whenever the inverse in \eqref{eq:geff} exists.

\emph{Relation to Theorem~\ref{thm:bie}.} By \eqref{eq:mp-rep}, the densities that are constant on each circle correspond to the functions \eqref{eq:mp-rep} without multipoles, and those with Fourier modes of order at most $1$ to the functions with charges and dipoles. The monopole and dipole approximations are thus the Galerkin projections of the pencil $(\mathsf A,\mathsf B)$ onto these two spaces of densities. For a weight that is not constant on the circles, the projected pencil has the same form, with the modes on one circle coupled by the Fourier coefficients of the weight as in \eqref{eq:mp-steklov}.

\begin{remark}[Related expansions]\label{rem:hv}
Huang and Venkatraman expand the Steklov eigenvalues of a closed manifold with many small holes about their limits \cite[Theorem 1.3]{Huang_Venkatraman_2026}. Their holes are unweighted geodesic balls at maximally separated points. Each radius balances the boundary area of the hole against the weighted volume of its Voronoi cell. On the sphere with constant weight this gives $r=2b^{-1}$. The leading correction, of order $b^{-1}\log b$, is a Coulomb-type energy, as in our monopole approximation. The next term contains the integral of the squared gradient of the limiting eigenfunction over the holes, with coefficient $-\tfrac{d}{d-1}$ in dimension $d$, so $-2$ on a surface \cite[(1.12) and (1.14)]{Huang_Venkatraman_2026}. It corresponds to our dipoles. Section~\ref{ssec:large} computes the coefficient $16\pi$ and the constant of the term of order $b^{-1}$, and optimizes over the radii. Huang and Venkatraman pose optimization problems of this kind as an open direction. The dipoles change the value at order $b^{-1}$ by the factor $1+2f$, in agreement with their coefficient $-2$. The optimal radius is $r^\star=2b^{-1}$, their balanced radius.
\end{remark}

\subsection{Derivatives}\label{ssec:derivatives}
The parameters of cap $\ell$ are its center $c_\ell$, moved by $\dot c_\ell\in T_{c_\ell}\sphere$ through the rigid rotation with generator $g=c_\ell\times\dot c_\ell$ of the circle and its frame, its radius $r_\ell$, and the coefficients $(\alpha_{\ell,n},\beta_{\ell,n})$ of $\omega_\ell$. In the notation of \cite[\S2.3]{Chen_Osting_2026}, $\operatorname{op}(k)$ is the operator with kernel $k(x,\nu(x),y)$ and measure $\sin r_{\ell'}d\theta'$, $\mathcal S=\operatorname{op}(k_{\mathcal S})$, $\mathcal D^*=\operatorname{op}(k_{\mathcal D^*})$ with
\begin{equation}\label{eq:kernels}
k_{\mathcal S}(x,y)=-\frac1{2\pi}\log|x-y|,\qquad k_{\mathcal D^*}(x,\nu,y)=\frac1{2\pi}\frac{\nu\cdot y}{|x-y|^2},
\end{equation}
$\mathbf 1_\ell$ is the indicator of $\Gamma_\ell$ and $[\mathbf 1_\ell,\mathcal K]=\mathbf 1_\ell\mathcal K-\mathcal K\mathbf 1_\ell$. For a rotation generator $g$, the rotated kernel is $k_g(x,\nu,y):=\nabla_x k\cdot(g\times x)+\nabla_\nu k\cdot(g\times\nu)$. For $\mathcal K\in\{\mathcal S,\mathcal D^*\}$ with kernel $k$ from \eqref{eq:kernels}, the derivatives with respect to the radius are
\begin{align}
\partial_{r_\ell}\mathcal K&=\mathbf 1_\ell\operatorname{op}\bigl(-\nabla_x k\cdot\nu(x)+\nabla_\nu k\cdot x\bigr)-\operatorname{op}\bigl(\nabla_y k\cdot\nu(y)\bigr)\mathbf 1_\ell+\mathcal K\,\mathbf 1_\ell\cot r_\ell,\label{eq:radius-derivative}\\
\partial_{r_\ell}M&=M\Bigl(\mathbf 1_\ell\cot r_\ell-\frac{2\pi\cos r_\ell}{|\partial\Omega|}\Bigr).\notag
\end{align}

\begin{theorem}[Derivatives through the boundary operators]\label{thm:derivatives}
Let $(\sigma_k,u_k)$ be a simple eigenpair of \eqref{eq:weighted}, with $\int_{\partial\Omega}\rho\,u_k^2\,ds=1$, and $\varphi_k$ its density in Theorem~\ref{thm:bie}. For every parameter $p$,
\begin{equation}\label{eq:prop1}
\partial_p\sigma_k=\int_{\partial\Omega}u_k\bigl[(\partial_p\mathsf A-\sigma_k\partial_p\mathsf B)\varphi_k\bigr]ds,\qquad\partial_p\bar\sigma_k=L_\rho\,\partial_p\sigma_k+\sigma_k\,\partial_pL_\rho,
\end{equation}
where, with $g=c_\ell\times\dot c_\ell$,
\begin{center}\footnotesize\renewcommand{\arraystretch}{1.6}
\begin{tabular}{@{}llll@{}}\toprule
$p$ & $\partial_p\mathsf A$ & $\partial_p\mathsf B$ & $\partial_pL_\rho$\\\midrule
$\alpha_{\ell,n}$ & $0$ & $\rho_\ell(\theta)\cos(n\theta)\,\mathbf 1_\ell\,\mathcal S_0$ & $\sin r_\ell\int_0^{2\pi}\rho_\ell\cos(n\theta)\,d\theta$\\
$\beta_{\ell,n}$ & $0$ & $\rho_\ell(\theta)\sin(n\theta)\,\mathbf 1_\ell\,\mathcal S_0$ & $\sin r_\ell\int_0^{2\pi}\rho_\ell\sin(n\theta)\,d\theta$\\
$\dot c_\ell\in T_{c_\ell}\sphere$ & $\bigl[\mathbf 1_\ell,\operatorname{op}(k_{\mathcal D^*,g})\bigr](I-M)$ & $\rho\,\bigl[\mathbf 1_\ell,\operatorname{op}(k_{\mathcal S,g})\bigr](I-M)$ & $0$\\
$r_\ell$ & $(\partial_{r_\ell}\mathcal D^*)(I-M)-(\mathcal D^*+\tfrac12I)\,\partial_{r_\ell}M$ & $\rho\,\bigl[(\partial_{r_\ell}\mathcal S)(I-M)+(I-\mathcal S)\,\partial_{r_\ell}M\bigr]$ & $\cos r_\ell\int_0^{2\pi}\rho_\ell\,d\theta$\\\bottomrule
\end{tabular}
\end{center}
At an $m$-fold eigenvalue $\sigma$, let the columns of $\Phi$ be the densities $\varphi_k$ of the cluster and $\Psi=\mathcal S_0\Phi$ their boundary values, with $\Psi^{\mathsf T}\mathsf B\Phi=I$. Then the cluster matrix $\bar{\mathsf S}=L_\rho\Psi^{\mathsf T}\mathsf A\Phi$ of \cite[(21)--(22)]{Chen_Osting_2026} has the derivative $\partial_p\bar{\mathsf S}=L_\rho\Psi^{\mathsf T}(\partial_p\mathsf A-\sigma\,\partial_p\mathsf B)\Phi+\sigma\,(\partial_pL_\rho)\,I$.
\end{theorem}

Formula \eqref{eq:prop1} and the weight rows are \cite[Proposition 1]{Chen_Osting_2026}. The proof there uses only that $\mathcal S_0$ is symmetric for the bilinear pairing in $ds$ and that $u_k=\mathcal S_0\varphi_k$ is the left eigenvector of the pencil, and both facts hold on the sphere with the round measure. In a symmetry class, the derivative with respect to a class parameter is the sum of \eqref{eq:prop1} over the orbit. In \eqref{eq:radius-derivative}, the self-interaction block of $\Gamma_\ell$, with the measure included, is $\partial_{r_\ell}(\mathcal S)_{\ell\ell}=\cot r_\ell\,(\mathcal S)_{\ell\ell}-\tfrac{\cos r_\ell}{2\pi}\,d\theta'$ and $\partial_{r_\ell}(\mathcal D^*)_{\ell\ell}=-\tfrac{\sin r_\ell}{4\pi}\,d\theta'$. In the rotation rows the self-interaction block cancels.

\subsection{Certificates}\label{ssec:optimization}
By the multiplicity bound of Section~\ref{sec:intro}, the first eigenspace of a maximizer has dimension at most three, and exactly three since a free boundary minimal surface with $b\ge2$ is not planar. We therefore maximize the smallest of three eigenvalues that coalesce at the optimum, and the nonsmoothness is confined to a cluster of size three, compared with clusters of size up to six in the genus-one problem of \cite{Chen_Osting_2026}. Within a symmetry class (Section~\ref{ssec:classes}) the cluster is simpler. In a class whose group has a three-dimensional irreducible representation, as in the polyhedral classes, $\sigma_1$ is triple for every member of the class, and $\bar\sigma_1$ is a smooth function to maximize. In the dihedral and pyramidal classes, the triple splits into a two-dimensional and a one-dimensional representation. Then we maximize the smaller of two smooth branches, which we tell apart by the action of the rotation by $\tfrac{2\pi}{j}$ about the $e_3$-axis on the eigenspace. Section~\ref{ssec:optimizer} describes how the critical points are computed.

Since $\bar\sigma_1$ is not differentiable at a triple eigenvalue, a maximizer is not a point where gradients vanish. The right notion is Clarke stationarity, and its certificate turns out to be the surface itself. Let $\bar{\mathsf S}(\mathbf p)$ be the $3\times3$ cluster matrix of the first three eigenvalues, the compression of the pencil to their eigenspace. Its eigenvalues are $\bar\sigma_1,\bar\sigma_2,\bar\sigma_3$, and its derivatives $\partial_p\bar{\mathsf S}$ come from Section~\ref{ssec:derivatives}.
\begin{theorem}[Certificate and spectral map]\label{thm:certificate}
Let $\mathbf p^\star$ be a configuration at which $\bar\sigma_1=\bar\sigma_2=\bar\sigma_3$, with eigenfunctions $u_1,u_2,u_3$ satisfying $\int_{\partial\Omega}\rho\,u_ku_{k'}\,ds=\delta_{kk'}$, and let $\Lambda$ be a symmetric $3\times3$ matrix with $\Lambda\succeq0$, $\operatorname{tr}\Lambda=1$ and $\langle\Lambda,\partial_p\bar{\mathsf S}(\mathbf p^\star)\rangle_F=0$ for every parameter $p$ (centers, radii and weight modes of every cap). Then $\mathbf p^\star$ is a Clarke stationary point of $\bar\sigma_1$, and if $\operatorname{rank}\Lambda=3$ the spectral map
\begin{equation}\label{eq:spectral-map}
X:=(L_\rho\Lambda)^{1/2}(u_1,u_2,u_3)\colon\overline\Omega\to\R^3
\end{equation}
is a conformal, possibly branched, free boundary minimal immersion into $\ball$. It satisfies $|X|=1$ on $\partial\Omega$, meets $\sphere$ orthogonally, and $\int_{\partial\Omega}\rho XX^{\mathsf T}ds=L_\rho\Lambda$. Its boundary length and area are
\[
L\bigl(X(\partial\Omega)\bigr)=\bar\sigma_1^\star,\qquad \operatorname{Area}\bigl(X(\Omega)\bigr)=\tfrac12\bar\sigma_1^\star.
\]
\end{theorem}
This is the adaptation to the punctured sphere of \cite[Theorem 2]{Chen_Osting_2026}, which in turn is the Fraser--Schoen correspondence \cite{Fraser_2011,Fraser_2015} between maximizers of $\bar\sigma_1$ and free boundary minimal surfaces, stated at a Clarke stationary point with the multiplier $\Lambda$ in place of maximality. The proof there applies verbatim. Its only input specific to the surface is that the directions $(\mathbf c,\mathbf r)$ span the tangent space of the moduli space of circle domains \cite{HeSchramm1993}. We compute $\Lambda$ by least squares over the \emph{full} parameter space, that is, all centers, radii and cosine and sine modes of every cap, not only the parameters of the class. We report the residual
\begin{equation}\label{eq:K}
\mathcal K:=\max_p\bigl|\langle\Lambda,\partial_p\bar{\mathsf S}\rangle_F\bigr|,
\end{equation}
the maximum over all parameters $p$. A class optimum is a critical point of the full problem exactly when $\mathcal K$ vanishes, and this is how we check the symmetry-breaking directions. By Theorem~\ref{thm:certificate}, $\Lambda$ is the inertia tensor of the boundary curves of the surface \eqref{eq:spectral-map}. A small eigenvalue of $\Lambda$ indicates a surface that is nearly flat in that direction (Section~\ref{ssec:equatorial}), and $\operatorname{rank}\Lambda=3$ makes $X$ a surface in $\ball$ rather than in a plane. Since the coordinates of $X$ are first eigenfunctions, $X$ is then an embedding, star-shaped and free of branch points \cite[Proposition 8.1]{Fraser_2015}. The area density of \eqref{eq:checks} tests the absence of branch points.

\subsection{Verification of the surface}\label{ssec:checks}
In exact arithmetic, Theorem~\ref{thm:certificate} makes $X$ a free boundary minimal surface. To measure how far the computed surface is from one, we test its defining properties on the traces, as in \cite[\S2.6]{Chen_Osting_2026}. The boundary must lie on the unit sphere. The surface must meet it orthogonally and conformally, meaning that along each boundary circle the tangential and normal derivatives of $X$ have equal length. The first eigenspace must also be a well-separated triple. We measure these properties by
\begin{subequations}\label{eq:checks}
\begin{align}
\mathcal R_{\rm sph}&:=\max_{\partial\Omega}\bigl||X|-1\bigr|,\qquad
\mathcal C_X:=\max_{\partial\Omega}\frac{\bigl|\tfrac{|\partial_\theta X|^2}{\sin^2r_\ell}-|\partial_\nu X|^2\bigr|}{\tfrac{|\partial_\theta X|^2}{\sin^2r_\ell}+|\partial_\nu X|^2},\label{eq:checks-geometry}\\
\delta_3&:=\frac{\bar\sigma_3-\bar\sigma_1}{\bar\sigma_1},\qquad \lambda_{\min}(\Lambda),\qquad \frac{\bar\sigma_4}{\bar\sigma_1},\qquad \min_\Omega J,\label{eq:checks-spectral}
\end{align}
\end{subequations}
where $\mathcal C_X$ is the relative conformality defect, $\delta_3$ the cluster spread, $\lambda_{\min}(\Lambda)$ the smallest eigenvalue of the certificate (the surface is an immersion into $\ball$ only when $\operatorname{rank}\Lambda=3$), $\tfrac{\bar\sigma_4}{\bar\sigma_1}$ the gap to the next eigenvalue, and $J:=\tfrac12|\nabla X|^2$ the area density of $X$ with respect to the round metric, which vanishes at a branch point. The boundary length of a free boundary minimal surface built from first eigenfunctions equals $\bar\sigma_1$, so the length of $X(\partial\Omega)$ minus $\bar\sigma_1^\star$ tests the scaling. Section~\ref{ssec:verification} and Table~\ref{tab:best} report the geometric quantities \eqref{eq:checks-geometry}, the spectral quantities \eqref{eq:checks-spectral} and the residual $\mathcal K$ of \eqref{eq:K}.

\section{Numerical methods}\label{sec:numerics}
This section describes the discretization of the boundary integral formulation of Section~\ref{sec:method} (Section~\ref{ssec:implementation}), the computation and certification of critical points (Section~\ref{ssec:optimizer}) and of their Morse indices (Section~\ref{ssec:morse}), and the tests of the solver (Section~\ref{ssec:validation}). The method follows \cite[\S\S2--3]{Chen_Osting_2026}, adapted to the sphere.

\subsection{Discretization}\label{ssec:implementation}
Each circle carries $N_\ell$ equispaced nodes, $N=\sum_\ell N_\ell$ in all. We discretize the pencil by the Nystr\"om method of \cite[\S3]{Chen_Osting_2026}, with Kress' product quadrature \cite[\S12.3]{Kress_2014} for the logarithm in \eqref{eq:self-blocks} and the trapezoidal rule for everything else. The eigenvalues converge exponentially in $N_\ell$ (Section~\ref{ssec:validation}). A symmetric eigensolver gives the full spectrum of the dense pencil, and every optimum is re-evaluated with $2N_\ell$ nodes. The area density $J$ of \eqref{eq:checks} is evaluated on a triangulation of $\Omega$, with $X$ computed at its vertices from the layer potentials.

\emph{Symmetry reduction.} For the icosahedral classes at large $b$, we replace the dense pencil by a smaller one. The nodes of each orbit are placed by the group, so the group permutes the nodes and the pencil commutes with this permutation. The triple $\bar\sigma_1$ belongs to the vector representation $T_{1u}$ of $I_h$, and on the $T_{1u}$-isotypic component the pencil factors as $I_3\otimes(A_{\rm red}-\sigma B_{\rm red})$. We solve one of the three identical blocks. Its dimension is about $\tfrac{3N}{|I_h|}$, since each free orbit of $|I_h|$ nodes carries the regular representation, which contains $T_{1u}$ three times. The triple thus becomes a simple eigenvalue of the reduced pencil. To confirm that $\bar\sigma_1$ lies in $T_{1u}$, we also compute the lowest eigenvalue of every other isotypic component. The gradient is computed in adjoint form on the reduced pencil. For $b=12$, $42$, $92$ and $162$, the reduced values agree with the dense solver to $10^{-12}$ and the gradients to $10^{-11}$. One eigensolve at $b=162$ takes $0.9$ seconds instead of $8.9$, and the gain grows with $b$. For these classes $\Lambda=\tfrac13I$, and the residual \eqref{eq:K} in the symmetry-breaking directions vanishes, because it is the derivative of an invariant function at a symmetric configuration. The class gradient is then the full certificate.

\subsection{Optimization}\label{ssec:optimizer}
We use the two-stage procedure of \cite[\S2.4]{Chen_Osting_2026}: gradient ascent, followed by maximin and sequential quadratic programming steps on the cluster matrix $\bar{\mathsf S}$ of Section~\ref{ssec:optimization}, with the gradients of Theorem~\ref{thm:derivatives} and the disjointness of the caps as a constraint. The procedure is the same within a symmetry class and for the unrestricted problem, in which every cap is its own orbit and the M\"obius freedom is fixed by \eqref{eq:normalization}. A Newton refinement of the stationarity system $\langle\Lambda,\partial_p\bar{\mathsf S}\rangle_F=0$, $\bar{\mathsf S}=\bar\sigma_1I$ in $(\mathbf p,\Lambda)$, as in \cite[\S3.2]{Chen_Osting_2026}, then brings the residual \eqref{eq:K} to rounding in one or two steps. In the unrestricted problem, the normal form transports the conformal factor of the dilation into the weight, so twelve weight modes per cap are needed to reach rounding. Every optimum is re-evaluated with twice the nodes per circle, and restarted there if the first three $\bar\sigma_k$ move by more than $10^{-9}$. It is also re-optimized with twice the weight modes per cap, and the modes are doubled again until $\bar\sigma_1^\star$ moves by less than $10^{-9}$. A configuration is \emph{certified} when three conditions hold: its cluster spread $\tfrac{\bar\sigma_3-\bar\sigma_1}{\bar\sigma_1}$ is below $10^{-11}$, the first three $\bar\sigma_k$ move by less than $10^{-9}$ when the nodes are doubled, and the residual $\mathcal K$ of \eqref{eq:K} over the full parameter space, with the cosine and sine modes $n\le12$ of every cap, is below $10^{-7}$. The tables report $\bar\sigma_1^\star$ to eight decimals, the accuracy of the doubling test, and the text rounds it. The tables also report $\mathcal K$ when it falls short of the threshold. Class optima that are not critical points of the full problem at all, with $\mathcal K$ of order one, are left out. The symmetry group of each result is determined a posteriori, since a class optimum may be more symmetric than its class.

\subsection{Morse index}\label{ssec:morse}
At a certified critical point the multiplier $\Lambda^\star$ is positive definite, so every direction that splits the triple $\bar\sigma_1=\bar\sigma_2=\bar\sigma_3$ decreases $\bar\sigma_1$ to first order. The second-order behavior is that of the Lagrangian $\ell(\mathbf p)=\langle\Lambda^\star,\bar{\mathsf S}(\mathbf p)\rangle_F$ on the tangent space of the manifold where the triple persists \cite{Overton_Womersley_1993,Lewis_Overton_1996}. We form the Hessian of $\ell$ on an orthonormal basis of this tangent space by differences of the gradient of Theorem~\ref{thm:derivatives}, with $\Lambda^\star$ fixed and step $10^{-4}$, and symmetrize it. The cluster basis of each displaced configuration is aligned to the reference one by an orthogonal Procrustes rotation. The three rotations and the three dilations leave $\bar\sigma_1$ unchanged and are projected out, as in \cite[\S4.5]{Chen_Osting_2026}. The \emph{Morse index} is the number of eigenvalues of the remaining matrix above $10^{-8}$ times its largest eigenvalue in absolute value. It counts the directions in which $\bar\sigma_1$ increases, so a local maximizer has index $0$. The Hessian uses six weight modes per cap. For the geodesic icosahedra with $b\ge362$ we compute it block by block in the irreducible representations of $I_h$, with three modes.

\subsection{Validation}\label{ssec:validation}
For caps at $\pm e_3$ of common radius $r$ and $\rho\equiv1$, the separated solution of \cite[\S3.1]{Kao_2020} gives the branches $\bar\sigma=4\pi\cos r$ (double) and $\tfrac{4\pi}{\log\cot(r/2)}$ (radial) of the Example in Section~\ref{ssec:uniformization}. The solver reproduces both to rounding at $N_\ell=16$, and maximizing over $r$ and a free weight returns the critical catenoid with constant weight.

For three caps with random weights, the derivatives \eqref{eq:prop1} with respect to the radii, the three rotation generators of each cap and the weight coefficients agree with central differences of the discrete $\bar\sigma_1$ to $10^{-9}$--$10^{-10}$. The sum over the caps for a common rotation generator vanishes to $10^{-15}$. Transporting a configuration by the dilation $z\mapsto1.7z$, with caps and weight transported as in Section~\ref{ssec:uniformization}, leaves the first four $\bar\sigma_k$ unchanged to all printed digits, as M\"obius invariance requires.

Started from the class maximizers for $b=4,6,12$ and the square antiprism, with centers and radii perturbed by $10^{-2}$, the unrestricted optimizer passes through both stages and the refinement and returns the class values to $10^{-12}$ with $\mathcal K\le3\times10^{-13}$. From a perturbed cube it leaves the cube, which is a critical point of its class but a saddle of the unrestricted problem, and reaches the square antiprism.

At the configurations of Table~\ref{tab:catalog}, $|\bar\sigma_1(N_\ell)-\bar\sigma_1(128)|$ decays exponentially in the number of nodes per circle, as in \cite[\S5]{kim2025layer}. It reaches rounding from $N_\ell=32$ for $b=3$ and from $N_\ell=24$ for $b\ge4$.

\section{Computational experiments and surface zoo}\label{sec:experiments}
We now apply the method of Sections~\ref{sec:method} and \ref{sec:numerics} to many symmetry classes. The classes are those of Section~\ref{ssec:classes} and Table~\ref{tab:classes}. Section~\ref{ssec:zoo} compares the results with previous work. Section~\ref{ssec:families} follows the families as $b$ grows and relates them to the limit surfaces of Section~\ref{sec:fcms}. Section~\ref{ssec:best} describes the searches beyond the families and reports the best configuration found for each $b$, and Section~\ref{ssec:large} the behavior for large $b$. Section~\ref{ssec:catalog} collects every critical point in the catalog of the zoo.

\subsection{Comparison with previous work}\label{ssec:zoo}
Table~\ref{tab:zoo} lists every configuration reported by \cite{Girouard_Lagace_2021}, \cite{Kao_2020} or \cite{Schulz_gallery}, together with the critical point of this paper that reproduces it. By \eqref{eq:area}, the areas of \cite{Schulz_gallery} are directly comparable with our values. The earlier values agree with ours to within the accuracy of the methods used there. At $b=6$ the value of \cite{Kao_2020} lies $4.3\times10^{-3}$ above our octahedral value. Their $b=6$ domain, one small hole surrounded by four equal holes on a square, is the octahedral picture in the disk chart, so that entry appears to be inaccurate, and the value of \cite{Schulz_gallery} agrees with ours. Their ``irregular'' configuration at $b=20$ is presumably our Fekete optimum.

\begin{table}[tbp]\centering\scriptsize\setlength{\tabcolsep}{4pt}
\begin{tabular}{@{}rlC{5.2cm}rrrr@{}}\toprule
$b$ & label & configuration found here & $\bar\sigma_1^\star$ & $2\operatorname{Area}$ \cite{Schulz_gallery} & \cite{Kao_2020} & \cite{Girouard_Lagace_2021}\\\midrule
2 & -- & \zoo{b0002/antipodal_pair}{antipodal pair}, $D_{\infty h}$ & 10.47478066 & 10.4748 & 10.4748 & \\
\addlinespace[1.5pt]\cmidrule(lr){1-7}\addlinespace[1.5pt]
3 & P3 & \zoo{b0003/equatorial_3_gon}{equatorial 3-gon}, $D_{3h}$ & 12.01198922 & 12.0119 & 12.0120 & \\
\addlinespace[1.5pt]\cmidrule(lr){1-7}\addlinespace[1.5pt]
4 & Y3 & \zoo{b0004/tetrahedron}{tetrahedron}, $T_d$ & 13.66761484 & 13.6672 & 13.6676 & 13.6675\\
4 & P4 & \zoo{b0004/equatorial_4_gon}{equatorial 4-gon}, $D_{4h}$ & 12.42990413 & 12.4300 &  & \\
\addlinespace[1.5pt]\cmidrule(lr){1-7}\addlinespace[1.5pt]
5 & P3 & \zoo{b0005/triangular_bipyramid}{triangular bipyramid}, $D_{3h}$ & 14.46866597 & 14.4683 & 14.4687 & \\
5 & Y4 & \zoo{b0005/square_pyramid}{square pyramid}, $C_{4v}$ & 14.43163401 & 14.4312 &  & \\
5 & P5 & \zoo{b0005/equatorial_5_gon}{equatorial 5-gon}, $D_{5h}$ & 12.53139727 & 12.5312 &  & \\
\addlinespace[1.5pt]\cmidrule(lr){1-7}\addlinespace[1.5pt]
6 & P4 & \zoo{b0006/octahedron}{octahedron}, $O_h$ & 15.42492045 & 15.4249 & 15.4292 & 15.4249\\
6 & P3 & \zoo{b0006/triangular_prism}{triangular prism}, $D_{3h}$ & 15.12598742 & 15.1261 &  & \\
6 & Y5 & \zoo{b0006/pentagonal_pyramid}{pentagonal pyramid}, $C_{5v}$ & 14.81894243 & 14.8194 &  & \\
6 & P6 & \zoo{b0006/equatorial_6_gon}{equatorial 6-gon}, $D_{6h}$ & 12.55681968 & 12.5576 &  & \\
\addlinespace[1.5pt]\cmidrule(lr){1-7}\addlinespace[1.5pt]
7 & C5 & \zoo{b0007/pentagonal_bipyramid}{pentagonal bipyramid}, $D_{5h}$ & 15.95213199 & 15.9524 & 15.9520 & \\
7 & C2 & \zoo{b0007/quadruple_plus_pair_plus_1pc}{quadruple + pair + 1PC}, $C_{2v}$ & 15.91819689 & 15.9182 &  & \\
7 & C3a & \zoo{b0007/triangular_antiprism_plus_1pc}{triangular antiprism + 1PC}, $C_{3v}$ & 15.91577614 & 15.9153 &  & \\
7 & Y2 & \zoo{b0007/quadruple_plus_pair_plus_1pc_pair_near_the_pole}{quadruple + pair + 1PC, pair near the pole}, $C_{2v}$ & 15.46032046 & 15.4598 &  & \\
7 & C3 & \zoo{b0007/triangular_prism_plus_1pc}{triangular prism + 1PC}, $C_{3v}$ & 15.45115152 & 15.4510 &  & \\
7 & C6 & \zoo{b0007/hexagonal_pyramid}{hexagonal pyramid}, $C_{6v}$ & 15.04160665 & 15.0419 &  & \\
7 & C7 & \zoo{b0007/equatorial_7_gon}{equatorial 7-gon}, $D_{7h}$ & 12.56361680 & 12.5657 &  & \\
\addlinespace[1.5pt]\cmidrule(lr){1-7}\addlinespace[1.5pt]
8 & C4a & \zoo{b0008/square_antiprism}{square antiprism}, $D_{4d}$ & 16.49540952 & 16.4956 & 16.4954 & \\
8 & C2 & \zoo{b0008/four_stacked_pairs}{four stacked pairs}, $D_{2d}$ & 16.49051351 & 16.4905 &  & \\
8 & C4 & \zoo{b0008/cube}{cube}, $O_h$ & 16.42485423 & 16.4249 &  & 16.4249\\
8 & C2 & \zoo{b0008/three_pairs_plus_2pc}{three pairs + 2PC}, $D_{2d}$ & 16.29009761 & 16.2901 &  & \\
8 & Y2 & \zoo{b0008/quadruple_plus_pair_plus_2pc}{quadruple + pair + 2PC}, $C_{2v}$ & 16.28575637 & 16.2857 &  & \\
8 & C6 & \zoo{b0008/hexagonal_bipyramid}{hexagonal bipyramid}, $D_{6h}$ & 16.27617583 & 16.2765 &  & \\
8 & C3 & \zoo{b0008/triangular_prism_plus_2pc}{triangular prism + 2PC}, $D_{3h}$ & 15.68321560 & 15.6831 &  & \\
8 & C7 & \zoo{b0008/heptagonal_pyramid}{heptagonal pyramid}, $C_{7v}$ & 15.18338502 & 15.1833 &  & \\
8 & C8 & \zoo{b0008/equatorial_8_gon}{equatorial 8-gon}, $D_{8h}$ & 12.56554312 & 12.5657 &  & \\
\addlinespace[1.5pt]\cmidrule(lr){1-7}\addlinespace[1.5pt]
9 &  & \zoo{b0009/three_stacked_triangles}{three stacked triangles}, $D_{3h}$ & 16.97069066 &  & 16.9707 & \\
\addlinespace[1.5pt]\cmidrule(lr){1-7}\addlinespace[1.5pt]
12 &  & \zoo{b0012/icosahedron}{icosahedron}, $I_h$ & 18.06871939 &  & 18.0687 & 18.0686\\
\addlinespace[1.5pt]\cmidrule(lr){1-7}\addlinespace[1.5pt]
15 &  & \zoo{b0015/five_stacked_triangles}{five stacked triangles}, $D_{3h}$ & 18.79339513 &  & 18.7934 & \\
\addlinespace[1.5pt]\cmidrule(lr){1-7}\addlinespace[1.5pt]
20 &  & \zoo{b0020/fekete_class}{Fekete class}, $D_{3h}$ & 19.70770682 &  & 19.7076 & \\
20 &  & \zoo{b0020/dodecahedron}{dodecahedron}, $I_h$ & 19.57267499 &  &  & 19.5718\\
\addlinespace[1.5pt]\cmidrule(lr){1-7}\addlinespace[1.5pt]
32 &  & \zoo{b0032/icosahedron_plus_dodecahedron}{icosahedron + dodecahedron}, $I_h$ & 21.01400764 &  &  & 21.0135\\
\bottomrule
\end{tabular}
\caption{Comparison with previous work (Sections~\ref{ssec:platonic} and \ref{ssec:zoo}). Each configuration reported in \cite{Schulz_gallery}, by \cite{Kao_2020} or by \cite{Girouard_Lagace_2021} is listed with the critical point of Table~\ref{tab:catalog} that reproduces it and its symmetry group. The values from \cite{Schulz_gallery} are twice the listed areas and are matched within $3\times10^{-3}$. The digit of a label of \cite{Schulz_gallery} is the order of the main rotation axis. Within each $b$ the rows are in decreasing order of $\bar\sigma_1^\star$.}\label{tab:zoo}
\end{table}

\cite{Schulz_gallery} lists the $27$ embedded genus-zero FBMS with $b\le8$ obtained by the equivariant min-max and gluing constructions of \cite{karpukhin2024,Carlotto_Schulz_Wiygul_2025,Franz_Schulz_2026}. Most of them come from the families (Table~\ref{tab:zoo}). The $b=8$ surfaces C2 and C2-gyrobifastigium come from four stacked pairs of caps and from three pairs of caps with caps at both poles. The boundary curves of the second have the $D_{2d}$ symmetry of the gyrobifastigium, although its class imposes only $C_2$. Three others, C2 and Y2 for $b=7$ and Y2 for $b=8$, come from classes with the group $C_{2v}$: four caps in a generic orbit and a pair of caps on a mirror plane, with one or two polar caps (Table~\ref{tab:classes}). The surface Y2 for $b=7$ is a saddle of index $2$, which maximization does not reach. Every surface of \cite{Schulz_gallery} with $b\le8$ is therefore among our critical points. The $C_{2v}$ classes also give two critical points for $b=8$ not in \cite{Schulz_gallery}, one of them distinct from but only $1.5\times10^{-6}$ below the triangular prism with two poles.

\subsection{The families}\label{ssec:families}
\begin{table}[tbp]\centering\scriptsize\setlength{\tabcolsep}{2.7pt}
\begin{tabular}{@{}rC{5.2cm}lrrr@{}}\toprule
$b$ & best configuration found & symmetry & $\bar\sigma_1^\star$ & $\mathcal K$ \eqref{eq:K} & $\lambda_{\min}(\Lambda^\star)$ \eqref{eq:checks}\\\midrule
2 & \zoo{b0002/antipodal_pair}{antipodal pair} & $D_{\infty h}$ & 10.47478066 & $1\times10^{-14}$ & 0.3052\\
3 & \zoo{b0003/equatorial_3_gon}{equatorial 3-gon} & $D_{3h}$ & 12.01198922 & $3\times10^{-10}$ & 0.1514\\
4 & \zoo{b0004/tetrahedron}{tetrahedron} & $T_d$ & 13.66761484 & $7\times10^{-14}$ & 0.3333\\
5 & \zoo{b0005/triangular_bipyramid}{triangular bipyramid} & $D_{3h}$ & 14.46866597 & $2\times10^{-13}$ & 0.3079\\
6 & \zoo{b0006/octahedron}{octahedron} & $O_h$ & 15.42492045 & $4\times10^{-10}$ & 0.3333\\
7 & \zoo{b0007/pentagonal_bipyramid}{pentagonal bipyramid} & $D_{5h}$ & 15.95213199 & $2\times10^{-13}$ & 0.3078\\
8 & \zoo{b0008/square_antiprism}{square antiprism} & $D_{4d}$ & 16.49540952 & $2\times10^{-13}$ & 0.3168\\
9 & \zoo{b0009/three_stacked_triangles}{three stacked triangles} & $D_{3h}$ & 16.97069066 & $8\times10^{-14}$ & 0.3296\\
10 & \zoo{b0010/square_antiprism_plus_2pc}{square antiprism + 2PC} & $D_{4d}$ & 17.37203030 & $2\times10^{-13}$ & 0.3298\\
11 & \zoo{b0011/c2v_17p6912}{$C_{2v}$: $4+3\times2+1$} & $C_{2v}$ & 17.69115333 & $4\times10^{-13}$ & 0.3257\\
12 & \zoo{b0012/icosahedron}{icosahedron} & $I_h$ & 18.06871939 & $3\times10^{-13}$ & 0.3333\\
13 & \zoo{b0013/c2v_18p3064}{$C_{2v}$: $4+4\times2+1$} & $C_{2v}$ & 18.30638945 & $9\times10^{-13}$ & 0.3295\\
14 & \zoo{b0014/hexagonal_antiprism_plus_2pc}{hexagonal antiprism + 2PC} & $D_{6d}$ & 18.56544315 & $2\times10^{-13}$ & 0.3252\\
15 & \zoo{b0015/five_stacked_triangles}{five stacked triangles} & $D_{3h}$ & 18.79339513 & $2\times10^{-13}$ & 0.3303\\
16 & \zoo{b0016/td_19p0089}{$T_d$: $12+4$} & $T_d$ & 19.00885253 & $1\times10^{-12}$ & 0.3333\\
17 & \zoo{b0017/c2v_19p2012}{$C_{2v}$: $2\times4+4\times2+1$} & $C_{2v}$ & 19.20121924 & $9\times10^{-13}$ & 0.3310\\
18 & \zoo{b0018/four_stacked_squares_plus_2pc}{four stacked squares + 2PC} & $D_{4d}$ & 19.38539453 & $1\times10^{-10}$ & 0.3326\\
19 & \zoo{b0019/fekete_class_19p5446}{Fekete class} & $C_{2v}$ & 19.54459459 & $7\times10^{-12}$ & 0.3300\\
20 & \zoo{b0020/fekete_class}{Fekete class} & $D_{3h}$ & 19.70770682 & $5\times10^{-13}$ & 0.3327\\
21 & \zoo{b0021/fekete_class}{Fekete class} & $C_{2v}$ & 19.85473178 & $5\times10^{-13}$ & 0.3331\\
22 & \zoo{b0022/fekete_class_19p9935}{Fekete class} & $T_d$ & 19.99354385 & $5\times10^{-13}$ & 0.3333\\
23 & \zoo{b0023/fekete_class}{Fekete class} & $D_{3}$ & 20.11635049 & $4\times10^{-13}$ & 0.3311\\
24 & \zoo{b0024/fekete_class}{Fekete class} & $O$ & 20.24322425 & $8\times10^{-13}$ & 0.3333\\
25 & \zoo{b0025/c1_20p3531}{$C_1$ (20.3531)} & $C_1$ & 20.35311176 & $1\times10^{-12}$ & 0.3329\\
26 & \zoo{b0026/fekete_class}{Fekete class} & $C_{2}$ & 20.46172826 & $9\times10^{-13}$ & 0.3329\\
27 & \zoo{b0027/fekete_class}{Fekete class} & $D_{5h}$ & 20.56911733 & $7\times10^{-13}$ & 0.3326\\
28 & \zoo{b0028/fekete_class}{Fekete class} & $T$ & 20.66487632 & $6\times10^{-13}$ & 0.3333\\
29 & \zoo{b0029/fekete_class}{Fekete class} & $D_{3}$ & 20.75364318 & $4\times10^{-12}$ & 0.3328\\
30 & \zoo{b0030/fekete_class}{Fekete class} & $D_{2}$ & 20.84340804 & $9\times10^{-13}$ & 0.3324\\
31 & \zoo{b0031/fekete_class_20p9278}{Fekete class} & $C_{3v}$ & 20.92777889 & $1\times10^{-12}$ & 0.3328\\
32 & \zoo{b0032/icosahedron_plus_dodecahedron}{icosahedron + dodecahedron} & $I_h$ & 21.01400764 & $5\times10^{-13}$ & 0.3333\\
33 & \zoo{b0033/c1_21p0835}{$C_1$ (21.0835)} & $C_1$ & 21.08353796 & $8\times10^{-12}$ & 0.3328\\
34 & \zoo{b0034/fekete_class}{Fekete class} & $D_{2}$ & 21.15725596 & $8\times10^{-13}$ & 0.3318\\
35 & \zoo{b0035/fekete_class}{Fekete class} & $C_{2}$ & 21.22645244 & $2\times10^{-12}$ & 0.3329\\
36 & \zoo{b0036/fekete_class}{Fekete class} & $D_{2}$ & 21.29451111 & $8\times10^{-13}$ & 0.3327\\
37 & \zoo{b0037/fekete_class_21p3600}{Fekete class} & $D_{5h}$ & 21.35998448 & $8\times10^{-13}$ & 0.3325\\
38 & \zoo{b0038/six_stacked_hexagons_plus_2pc_21p4221}{six stacked hexagons + 2PC} & $D_{6d}$ & 21.42213538 & $3\times10^{-12}$ & 0.3323\\
39 & \zoo{b0039/fekete_class_21p4831}{Fekete class} & $D_{3h}$ & 21.48308018 & $2\times10^{-12}$ & 0.3327\\
40 & \zoo{b0040/fekete_class_21p5411}{Fekete class} & $T_d$ & 21.54113373 & $1\times10^{-12}$ & 0.3333\\
42 & \zoo{b0042/geodesic_icosahedron_nu2}{geodesic icosahedron $\nu=2$} & $I_h$ & 21.64565725 & $2\times10^{-12}$ & 0.3333\\
60 & \zoo{b0060/60_orbit_on_the_vertex_to_vertex_arcs}{60-orbit on the vertex-to-vertex arcs} & $I_h$ & 22.28356339 & $8\times10^{-13}$ & 0.3333\\
72 & \zoo{b0072/icosahedron_plus_generic_60_orbit}{icosahedron + generic 60-orbit} & $I$ & 22.69609123 & $3\times10^{-12}$ & 0.3333\\
92 & \zoo{b0092/geodesic_icosahedron_nu3}{geodesic icosahedron $\nu=3$} & $I_h$ & 23.07686419 & $4\times10^{-11}$ & 0.3333\\
162 & \zoo{b0162/geodesic_icosahedron_nu4}{geodesic icosahedron $\nu=4$} & $I_h$ & 23.77221789 & $4\times10^{-10}$ & 0.3333\\
252 & \zoo{b0252/geodesic_icosahedron_nu5}{geodesic icosahedron $\nu=5$} & $I_h$ & 24.16173249 & $2\times10^{-9}$ & 0.3333\\
362 & \zoo{b0362/geodesic_icosahedron_nu6}{geodesic icosahedron $\nu=6$} & $I_h$ & 24.40228582 & $1\times10^{-10}$ & 0.3333\\
492 & \zoo{b0492/geodesic_icosahedron_nu7}{geodesic icosahedron $\nu=7$} & $I_h$ & 24.56163111 & $6\times10^{-9}$ & 0.3333\\
642 & \zoo{b0642/geodesic_icosahedron_nu8}{geodesic icosahedron $\nu=8$} & $I_h$ & 24.67287741 & $4\times10^{-9}$ & 0.3333\\
812 & \zoo{b0812/geodesic_icosahedron_nu9}{geodesic icosahedron $\nu=9$} & $I_h$ & 24.75377466 & $1\times10^{-10}$ & 0.3333\\
1002 & \zoo{b1002/geodesic_icosahedron_nu10}{geodesic icosahedron $\nu=10$} & $I_h$ & 24.81454590 & $5\times10^{-10}$ & 0.3333\\
1212 & \zoo{b1212/geodesic_icosahedron_nu11}{geodesic icosahedron $\nu=11$} & $I_h$ & 24.86142450 & $4\times10^{-10}$ & 0.3333\\
1442 & \zoo{b1442/geodesic_icosahedron_nu12}{geodesic icosahedron $\nu=12$} & $I_h$ & 24.89839157 & $5\times10^{-10}$ & 0.3333\\
1692 & \zoo{b1692/geodesic_icosahedron_nu13}{geodesic icosahedron $\nu=13$} & $I_h$ & 24.92808940 & $3\times10^{-10}$ & 0.3333\\
1962 & \zoo{b1962/geodesic_icosahedron_nu14}{geodesic icosahedron $\nu=14$} & $I_h$ & 24.95232944 & $7\times10^{-10}$ & 0.3333\\
2252 & \zoo{b2252/geodesic_icosahedron_nu15}{geodesic icosahedron $\nu=15$} & $I_h$ & 24.97238832 & $4\times10^{-9}$ & 0.3333\\
\bottomrule
\end{tabular}
\caption{The best configuration found for each $b$ (Sections~\ref{ssec:best} and \ref{ssec:thomson}), with the symmetry group of its surface, which may exceed that of its class (Table~\ref{tab:classes}). The value $\bar\sigma_1^\star$ uses $2N_\ell$ nodes per circle, and $\lambda_{\min}(\Lambda^\star)$ is the smallest eigenvalue of the certificate. Previous values are in Table~\ref{tab:zoo}, and the surfaces in Figures~\ref{fig:best-a}--\ref{fig:best-last}. They are the upper envelope of Figure~\ref{fig:families}.}\label{tab:best}
\end{table}

We follow the families in the order of Section~\ref{ssec:classes}. When a class contains a configuration of a more symmetric class, the optimizer goes there, so the symmetric configurations are stable within their larger classes.

\begin{remark}[Existence]\label{rem:existence}
Karpukhin, Kusner, McGrath and Stern prove that $\bar\sigma_1$ attains its supremum among metrics invariant under a finite group generated by reflections, provided the boundary components form enough orbits \cite[Theorem 9.15]{karpukhin2024}. Among our classes, this covers the polyhedral classes other than the chiral one at $b=72$, the equatorial polygons, the bipyramids, the prisms with two polar caps, the stacks with an odd number of rings, and the prisms, antiprisms and three stacked rings with one polar cap (their cases (3), (4), (5) and (7)). It does not cover the pyramids, or the prisms and antiprisms without polar caps, whose caps form too few orbits. The antiprisms with two polar caps and the stacks with an even number of rings have symmetry groups that are not generated by reflections, and the theorem gives a maximizer only for a subgroup.
\end{remark}

\emph{Polyhedral configurations.}\label{ssec:platonic} The classes with a polyhedral group are the Platonic solids, the icosidodecahedron, the icosahedron with the dodecahedron ($b=32$), a free $60$-orbit on the vertex-to-vertex arcs of the icosahedron ($b=60$), the chiral configuration with $b=72$, and the geodesic icosahedra with $\nu=2,\dots,15$ ($b=42$ to $2252$). Their surfaces are shown in Figure~\ref{fig:type-polyhedral} and, for large $b$, in Figures~\ref{fig:best-b42}--\ref{fig:best-b2252}. The certificate is $\Lambda=\tfrac13I$ to $10^{-10}$, as it must be for a three-dimensional irreducible representation. In particular, the Platonic configurations are stable within their larger classes. Even so, the cube is a saddle beaten by the square antiprism, and the dodecahedron and the icosidodecahedron are beaten by stacked rings. The truncated icosahedron, with the $60$-orbit one third of the way along each edge, is not a critical point: freed, the orbit moves to $19.6^\circ$ from the nearer vertex and the value rises by $6\times10^{-3}$. At $b=72$ the chiral class, started from the Thomson position, beats an orbit on a mirror arc and has symmetry $I$.

\emph{The equatorial family.}\label{ssec:equatorial} With all caps on the equator, the class optima increase toward $4\pi$, the value of the doubly covered disk, and agree with \cite{Schulz_gallery} to all its digits. The deficit $4\pi-\bar\sigma_1^\star$ falls exponentially, and matched asymptotics predict it without a fitted parameter. Write the two sheets as graphs $\pm u$ over the equatorial disk. To first order, the free boundary condition is $\partial_r u=u$ on $r=1$, which a constant height $h>0$ violates. The $b$ bridges at the roots of unity act as boundary sources of a strength $A$: $\partial_r u-u=-A\sum_{j}\delta_{\theta_j}$ with $\theta_j=\tfrac{2\pi j}{b}$. By Fourier series, the solution on $r=1$ is the Green's function of \cite[(4.7)]{Folha_Pacard_Zolotareva_2017},
\begin{equation*}
u(\theta)=\frac{A}{\pi}\Bigl(\frac b2-\sum_{n\ge1}\frac{\cos nb\theta}{n-1/b}\Bigr).
\end{equation*}
So the sheets lie at height $h=\tfrac{Ab}{2\pi}$, and near the bridge at $\theta_j$
\begin{equation*}
u=\frac A\pi\Bigl(\frac b2+\log\bigl(b|\theta-\theta_j|\bigr)+\eta_b\Bigr),\qquad\eta_b=-\sum_{k\ge2}\frac{\zeta(k)}{b^{k-1}}=\psi\bigl(1-\tfrac1b\bigr)+\gamma,
\end{equation*}
with $\zeta$ the Riemann zeta function, $\psi$ the digamma function and $\gamma$ Euler's constant. This must match a half-catenoid of waist $\tau$, cut in half by $\partial\ball$, with axis the vertical line through $e^{i\theta_j}$. Its radius at height $z$ is $\tau\cosh\tfrac z\tau$, so at distance $s$ from the axis its height is $\tau\operatorname{arcosh}\tfrac s\tau=\tau\log\tfrac{2s}{\tau}+O\bigl(\tfrac{\tau^3}{s^2}\bigr)$. In the overlap $\tau\ll|\theta-\theta_j|\ll\tfrac1b$, the logarithms match when $A=\pi\tau$ and the constants match when $\log\tfrac2\tau=\tfrac b2+\log b+\eta_b$. Hence $\tau=\tfrac2b\,e^{-b/2-\eta_b}$ and $h=e^{-b/2-\eta_b}$, which is the neck scale of the construction in \cite[Section~7]{Folha_Pacard_Zolotareva_2017}. By \eqref{eq:area}, the deficit is twice the area lost against the doubly covered disk. To compute it, remove half-disks of radius $\varepsilon$ around the bridges, $\tau\ll\varepsilon\ll\tfrac1b$. To second order, the area of each sheet over the rest of the disk differs from the area of the rest of the disk by $\tfrac12\int|\nabla u|^2-\tfrac12\int_{r=1}u^2\,d\theta$. After integration by parts only the half-circles contribute, in total $-\tfrac\pi2b\tau^2\bigl(\tfrac b2+\log(b\varepsilon)+\eta_b\bigr)$ per sheet. Within distance $\varepsilon$ of its axis, each half-catenoid exceeds the two half-disks it replaces by $\pi\tau^2\bigl(\log\tfrac{2\varepsilon}{\tau}-\tfrac12\bigr)$. In the sum $\varepsilon$ cancels, and by the matching so does $\tfrac b2$. The surface loses $\tfrac\pi2b\tau^2$, the area of $b$ half-disks of radius $\tau$:
\begin{equation}\label{eq:equatorial}
4\pi-\bar\sigma_1^\star(b)\approx\pi b\tau^2=\frac{4\pi}{b}\,e^{-b-2\eta_b}=\frac{4\pi}{b}\,e^{-b}\Bigl(1+\frac{\pi^2}{3b}+O(b^{-2})\Bigr).
\end{equation}
This is the dashed curve in Figure~\ref{fig:families}. It agrees with the computed deficits to a relative error of $6\times10^{-3}$ at $b=7$ and $2\times10^{-3}$ at $b=8$, of the order of $\tau$, which is $1.1\times10^{-2}$ and $5.7\times10^{-3}$ there. The relative error is at most $10^{-3}$ for $9\le b\le16$ and $5\times10^{-3}$ at $b=18$, where the deficit is only $1.3\times10^{-8}$.

The surfaces (Figure~\ref{fig:type-equatorial}) are two nearly flat disks joined by $b$ boundary bridges, the family of \cite{Folha_Pacard_Zolotareva_2017,Ketover_2016}. Franz and Schulz construct such a surface for every $b\ge2$ with the rotational symmetry of the $b$-gon, conjecturally the same family \cite[Theorem 5.1 and Remark 5.2]{Franz_Schulz_2026}, \cite[Conjecture 7.8]{Carlotto_Schulz_Wiygul_2025}. For $b=3$ the triangle is the maximizer, and for $b\ge4$ the equatorial surfaces are the lowest critical points found, as \cite[Conjecture 1.21]{karpukhin2024} predicts for the equatorial $b$-noid. They are basic reflection surfaces \cite[Definition 5.4]{karpukhin2024}, and Karpukhin, McGrath and Stern bound the supremum of $\bar\sigma_1$ over their symmetry class between $4\pi-e^{-\mu_2b}$ and $4\pi-e^{-\mu_1b}$ \cite[Theorem 1.12]{Karpukhin_McGrath_Stern_2025}. The asymptotics \eqref{eq:equatorial} suggest that the constant in the exponent is $1$, with prefactor $\tfrac{4\pi}{b}e^{-2\eta_b}$. This exponential decay contrasts with the $j^{-1}$ of the ring families in \eqref{eq:ringlaw}.

\emph{Rings and their limit problem.} The remaining families are built from rings of $j$ caps. As $j\to\infty$, a ring becomes a circle of $\Gamma$ in the transmission--Steklov problem \eqref{eq:transmission} on $\sphere$, with $\lambda$ the weight per unit length of the circle, while a polar cap stays a boundary circle. We call this the \emph{limit problem}, and its optimum over the positions and weights of the circles $\bar\sigma_\infty$. In the coordinates $(t,\phi)$ of Section~\ref{ssec:fcms-rev}, it is \eqref{eq:cyl}. The surface to which a family converges is its limit network, one of the free curve minimal surfaces of Section~\ref{sec:fcms}.

\emph{Pyramids.}\label{ssec:pyramids} A ring of $j$ caps with one cap at the pole gives the surfaces of Figure~\ref{fig:type-pyramid}. The pyramids match \cite{Schulz_gallery} to all digits and are the second-lowest critical points found at each $b$. As $j\to\infty$ they appear to converge to the cake stand of Proposition~\ref{prop:cake}, a horn with a round plate (Figure~\ref{fig:fcms}, Table~\ref{tab:limitnetworks}). The plate lies off the equatorial plane, on the side opposite the polar cap. For $j\ge4$ we know of no construction of the pyramids in $\ball$. For $j=3$ ($b=4$) the class contains the tetrahedron.

\emph{Bipyramids.}\label{ssec:bipyramids} A ring of $j$ caps on the equator with caps at both poles gives the surfaces of Figure~\ref{fig:type-bipyramid}, and the critical catenoid is the case $j=0$. Here the relative weight level of Section~\ref{ssec:classes} matters, since no M\"obius transformation other than a rotation preserves the class. Without it the values fall short by $10^{-2}$, except for $j=4$, whose optimum is the octahedron with all caps in one orbit. With it, the triangular and pentagonal bipyramids are the maximizers for $b=5$ and $7$. As $j\to\infty$ the bipyramids appear to converge to the wheel $\mathbb K_0\cup-\mathbb K_0$ of Proposition~\ref{prop:wheel}. Carlotto, Schulz and Wiygul construct, for every sufficiently large $j$, an embedded free boundary minimal surface of genus zero with $j+2$ boundary components and the symmetry of the bipyramid, which converges to this wheel \cite[Theorem 6.1]{Carlotto_Schulz_Wiygul_2025}. These surfaces appear to be the class maximizers, which exist by Remark~\ref{rem:existence}. Our bipyramids fit the conjectures of \cite[Conjectures 7.9 and 7.10]{Carlotto_Schulz_Wiygul_2025}: the family exists for every computed $j\ge3$, it is the octahedron at $j=4$, its values stay below $4\operatorname{Area}(\mathbb K_0)$, and its members are the best configurations found for $b=5$ and $7$ (and, as the octahedron, for $b=6$). The wheel gives $\bar\sigma_\infty=4\operatorname{Area}(\mathbb K_0)\approx17.5430$ (Figure~\ref{fig:fcms}), explicit in the positive root $L=2.0301901$ of \eqref{eq:wheel}. Each horn has conformal length $L$, corresponding to polar caps of angular radius $14.96^\circ$ on $\sphere$.

\emph{Prisms.}\label{ssec:prisms} Two aligned rings of $j$ caps at polar angles $\beta$ and $\pi-\beta$, with zero, one or two polar caps, give the surfaces of Figure~\ref{fig:type-prism}. With one polar cap the two rings become separate orbits, with their own polar angles, radii and weights. The optimum then moves one ring onto the equator, up to the M\"obius dilations along the axis, which commute with the symmetry. For $b=7$ the capped triangular prism and antiprism are the two surfaces of \cite{Schulz_gallery} with a threefold axis. Wherever the difference is resolved, the antiprism of the next paragraph beats the prism with the same caps, by a margin that shrinks with $j$ ($0.30$ for $j=3$, $2\times10^{-4}$ for $j=8$). With one polar cap the two agree to $2\times10^{-6}$ from $j=13$ on. Without polar caps the limit is the critical drum of Proposition~\ref{prop:drum} (Figure~\ref{fig:fcms}), with $\bar\sigma_\infty=\tfrac{8\pi}{1+W(\frac1e)}\approx19.6585$, where $W$ is the Lambert function.

\emph{Antiprisms.}\label{ssec:antiprisms} Rotating one ring of the prism by $\tfrac{\pi}{j}$ gives the antiprisms of Figure~\ref{fig:type-antiprism}. The square antiprism, with rings at $55.05^\circ$ against the cube's $54.74^\circ$, is the maximizer for $b=8$. With two polar caps, the square and hexagonal antiprisms are the maximizers for $b=10$ and $14$. For $j=3$ the class contains the octahedron. Without polar caps and for $j\ge4$, neither \cite{karpukhin2024} (Remark~\ref{rem:existence}) nor \cite{Kapouleas_Zou_2021} applies: the antiprisms lack the reflection in the equatorial plane, which \cite{Kapouleas_Zou_2021} imposes only to simplify the presentation.

\emph{Stacked rings.}\label{ssec:stacked} Stacks of three to six rings of $j$ caps (Table~\ref{tab:classes}) give the surfaces of Figure~\ref{fig:type-stacked}. Consecutive rings are rotated by $\tfrac{\pi}{j}$, every latitude is free subject to the symmetry, and there may be polar caps. The optimal latitudes spread the rings nearly evenly, and an optimum can be more symmetric than its class: six triangles with two poles merge pairwise into three hexagons. The stacks give the best values found for $b=9$, $15$, $18$ and $38$. At $b=9$ the maximizer is three stacked triangles, a triangular prism with a third triangle on the equator. The stacks beat the dodecahedron at $b=20$ and the icosidodecahedron at $b=30$ by about $0.1$, and the Fekete classes in turn beat the stacks there (Section~\ref{ssec:best}). At $b=32$ the six stacked pentagons with two poles return to the icosahedron with the dodecahedron. Rings of two caps are stacks too: four stacked pairs give the $b=8$ surface C2 of \cite{Schulz_gallery}. With one polar cap, the three rings become separate orbits. None of these is the best configuration for its $b$. Without polar caps the three rings converge to the drum stack $\mathcal C(3)$ (Table~\ref{tab:limitnetworks}). This is the free curve minimal surface of Proposition~\ref{prop:c3} (Figure~\ref{fig:fcms}), with $\bar\sigma_\infty\approx22.2295$ explicit in the root $t_1=0.9160481$ of \eqref{eq:c3} in $(\tfrac12,\infty)$.

\emph{Limits of the ring families.}\label{ssec:ringlimits} By \eqref{eq:fcms-value}, the limit values $\bar\sigma_\infty$ are twice the areas of the networks of Section~\ref{sec:fcms}. As an independent check, we maximized the first eigenvalue of the limit problem directly. In the coordinates $(t,\phi)$ the problem \eqref{eq:cyl} separates into Fourier modes, which we solve by shooting, and we optimize the positions and weights of the circles. The two computations agree to $10^{-6}$, and at each optimum the first eigenvalue is triple, except for the doubly covered disk, where it is double. So $\sigma=1$ is the first nonzero eigenvalue of every network of Section~\ref{sec:fcms}, including those with horns, and each network maximizes $\bar\sigma_1$ over the positions and weights of its circles. The class optima approach the limit values of Table~\ref{tab:limitnetworks} as
\begin{equation}\label{eq:ringlaw}
\bar\sigma_1^\star\approx\bar\sigma_\infty-A\,j^{-1}-B\,j^{-2},
\end{equation}
with the fitted coefficients $A$ and $B$ of Table~\ref{tab:limitnetworks}. The residuals are below $1.2\times10^{-2}$ for the families of one and two rings, below $2\times10^{-3}$ for three stacked rings, and below $8\times10^{-4}$ for three stacked rings with one polar cap. Fitting $\bar\sigma_\infty$ as well gives values within $9\times10^{-3}$ of the limits for the families of one and two rings (on $j\ge9$), and $22.514$ and $22.80$ for the stacks with one and two polar caps. The last fit uses only four values of $j$. The deficit is of order $j^{-1}$, without the logarithm of \eqref{eq:rate}, because the caps of a ring shrink like their spacing, with $j\,r\approx1.2$. The computed configurations also approach the limit networks in shape (Table~\ref{tab:limitnetworks}, Figure~\ref{fig:fcms}). The rate $j^{-1}$ agrees with a heuristic from the gluing constructions. Near each circle the necks are modeled on half of a singly periodic Scherk surface \cite[Definition 3.1]{Kapouleas_Zou_2021}, and each of the $j$ cores, scaled to the spacing of order $\tfrac1j$, changes the area by $O(j^{-2})$. The Scherk model also predicts that, as $j\to\infty$, each boundary curve of a ring spans the fraction $\tfrac{2\psi}{\pi}$ of its azimuthal period, with $\psi$ the angle between each sheet and $\partial\ball$ (Table~\ref{tab:limitnetworks}): $0.83$ for the pyramids, $0.72$ for the bipyramids, and $0.69$ for the prisms and antiprisms.

\begin{table}[tbp]\centering
\scriptsize\setlength{\tabcolsep}{2.5pt}\begin{tabular}{@{}C{3.5cm}C{2.85cm}r@{\hspace{9pt}}C{3.25cm}C{2.9cm}rrr@{}}\toprule
family & limit network & $\bar\sigma_\infty$ & heights of the circles & $\psi$ & $A$ & $B$ & fit $j$\\\midrule
equatorial & doubly covered disk & $12.566$ & $0$ & $90^\circ$ & -- & -- & --\\
pyramids & cake stand & $15.782$ & $0.266$; horn to $-0.757$ & $74.56^\circ$ & 2.8 & 10.1 & 5--20\\
bipyramids & wheel $\mathbb K_0\cup-\mathbb K_0$ & $17.543$ & $0$; horns to $\pm0.870$ & $64.62^\circ$ & 5.9 & 10.3 & 5--30\\
prisms, antiprisms & critical drum $\mathcal C(2)$ & $19.659$ & $\pm0.467$ & $62.18^\circ$ & 10.6 & 9.5 & 5--20\\
prisms, antiprisms + 1PC & drum with one horn & $20.486$ & $0.548$, $-0.302$; horn to $-0.948$ & $56.80^\circ$, $54.03^\circ$ & 13.0 & 5.6 & 5--19\\
prisms, antiprisms + 2PC & drum with two horns & $21.070$ & $\pm0.399$; horns to $\pm0.965$ & $50.40^\circ$ & 14.6 & 2.5 & 5--20\\
three stacked rings & drum stack $\mathcal C(3)$ & $22.229$ & $0$, $\pm0.674$ & $42.63^\circ$, $47.63^\circ$ & 18.3 & $-2.2$ & 7--13\\
three stacked rings + 1PC & drum stack with one horn & $22.506$ & $0.703$, $0.078$, $-0.575$; horn to $-0.985$ & $45.31^\circ$, $40.07^\circ$, $43.02^\circ$ & 18.8 & $-2.7$ & 9--20\\
three stacked rings + 2PC & drum stack with two horns & $22.731$ & $0$, $\pm0.613$; horns to $\pm0.988$ & $37.91^\circ$, $41.24^\circ$ & -- & -- & --\\
four stacked rings & $\mathcal C(4)$ & $23.354$ & $\pm0.281$, $\pm0.776$ & $33.65^\circ$, $39.08^\circ$ & -- & -- & --\\
polyhedral & $\partial\ball$ & $25.133$ & -- & $\to0^\circ$ & -- & -- & --\\\bottomrule
\end{tabular}
\caption[Limit networks of the families]{Limit networks as $j\to\infty$ (Section~\ref{sec:fcms}, Figure~\ref{fig:fcms}): the value $\bar\sigma_\infty$, twice the area counted with multiplicity, the heights of the circles on $\partial\ball$ (for a horn, also its boundary circle), and the angle $\psi$ of each sheet with $\partial\ball$ (the sheets meet at $180^\circ-2\psi$). Symmetric stacks are listed from the equator outward, and the networks with one horn from the disk to the horn, with the horn at the bottom as in Figure~\ref{fig:fcms}. $\mathcal C(k)$ is the stack of $k-1$ catenoidal bands closed by two flat disks \cite{Karpukhin_Noutchegueme_2025}, and the last row is its limit as $k\to\infty$. The values from shooting in the limit problem and from the balanced networks agree to $10^{-6}$. $A$ and $B$ are the coefficients of \eqref{eq:ringlaw} fitted over the range of $j$ in the last column (Section~\ref{ssec:ringlimits}), from the antiprisms in the rows with two rings.}\label{tab:limitnetworks}
\end{table}

\subsection{The best configuration found for each \texorpdfstring{$b$}{b}}\label{ssec:best}\label{ssec:verification}
\begin{figure}[tbp]\centering
\input{fig_tex/fig_families.tex}
\caption{Left: $\bar\sigma_1^\star$ against $b$ by family (Sections~\ref{ssec:families} and \ref{ssec:large}), with previous values \cite{Kao_2020,Girouard_Lagace_2021,Schulz_gallery} and the bounds $2\pi b$ \cite{Fraser_2015} and $8\pi$ \cite{Karpukhin_Stern_2021}. The curves are the asymptotic formula \eqref{eq:equatorial} (blue dashed), the rate formula \eqref{eq:rate-full} (dotted) and the fits \eqref{eq:ringlaw} (dash-dotted, Section~\ref{ssec:ringlimits}). The legend gives the limit network of each ring family as $j\to\infty$ (Figure~\ref{fig:fcms}, Table~\ref{tab:limitnetworks}). The prism and antiprism markers include zero, one and two polar caps. At large $b$ the antiprism rows have no polar caps (limit $19.659$) or two ($21.070$), and the prism row has one ($20.486$). Right: relative error of the multipole approximations \eqref{eq:models} of Section~\ref{ssec:approx}, evaluated on the centers, radii and circle weights of each configuration with $b\ge3$ (Section~\ref{ssec:thomson}). Gray dots are the monopole approximation $\bar\sigma^{(0)}_1$, and the family markers the dipole approximation $\bar\sigma^{(1)}_1$. The points are the rows of Table~\ref{tab:catalog}, with the best of each $b$ filled. The upper envelope on the left is Table~\ref{tab:best}.}\label{fig:families}
\end{figure}

The families need not contain the maximizer, so we also run searches that do not presuppose one. For $5\le b\le18$ we solve the unrestricted problem, every cap its own orbit, starting from the Thomson configuration, the minimizer of the Coulomb energy (Section~\ref{ssec:thomson}). We also start it from ten random configurations at each $b\le12$, and from saddle points pushed along their unstable directions. For $19\le b\le40$ we optimize instead in the class of the Fekete configuration (Section~\ref{ssec:large}). Its symmetry group splits the centers into orbits, each with one radius, one weight level, the free position of its stabilizer and the stabilizer-invariant weight modes. We start this \emph{Fekete class} from the Fekete configuration and from three random configurations of the class. At $b=32$ and $38$ the Fekete configuration already lies in a class of the families: the icosahedron with the dodecahedron, and six stacked hexagons with two polar caps. From $b=42$ on we run only icosahedral classes.

Table~\ref{tab:best} lists, for each $b$ computed, the configuration with the largest $\bar\sigma_1$ found, with its symmetry group and the verification quantities $\mathcal K$ and $\lambda_{\min}(\Lambda^\star)$ of \eqref{eq:K} and \eqref{eq:checks}. Figures~\ref{fig:best-a}--\ref{fig:best-last} show the surfaces. For $b=11$, $13$, $16$ and $17$ the best configuration comes from the unrestricted optimizer started at the Thomson configuration. To our knowledge, these values have not been reported before. For every $b$ from $19$ to $40$ except $25$, $32$, $33$ and $38$, it comes from the Fekete class, which beats the best of the families by between $3\times10^{-4}$ and $0.57$. At $b=32$ the best is the icosahedron with the dodecahedron, and at $b=38$ six stacked hexagons with two polar caps. At $b=25$ and $33$ the optimum of the Fekete class is a saddle of index $1$. Started from it along the unstable eigenvector, the unrestricted optimizer reaches a critical point without symmetry, with index $0$, which lies only $5.2\times10^{-7}$ and $1.2\times10^{-6}$ above the saddle. These two are the only configurations of Table~\ref{tab:best} with trivial symmetry group. Schulz conjectures \cite{Schulz_gallery} that a free boundary minimal surface of genus zero with trivial symmetry group exists with $61$ boundary components. Every surface of Table~\ref{tab:catalog} with $b\le18$ has a mirror plane. For $b=15$ and $16$ the best surfaces found have the groups $D_{3h}$ and $T_d$, larger than the chiral groups $D_3$ and $T$ of the Fekete configurations. The first chiral surface in Table~\ref{tab:catalog} has $b=19$ and the group $C_2$. The best configurations found at $b=23$, $24$ and $28$ have the rotation groups $D_3$, $O$ and $T$. At $b=24$ the cap centers form a snub cube, with $32$ triangular and $6$ square faces, as do the Fekete points.

The residual $\mathcal K$ is at most $4\times10^{-10}$ for $b\le92$ and at most $6\times10^{-9}$ for the geodesic icosahedra (Table~\ref{tab:best}). The other quantities of \eqref{eq:checks} are uniform, except for the boundary length. That length agrees with $\bar\sigma_1^\star$ to $4\times10^{-8}$ in every row except the equatorial $b$-gons with $14\le b\le20$, where it differs by up to $7\times10^{-5}$ ($b=18$). There the caps nearly touch, and the defect grows with $\mathcal K$. The cluster spread $\delta_3$ is below $3\times10^{-12}$ (for $b\ge362$ it vanishes by construction, and the triple is checked against the other isotypic components, Section~\ref{ssec:implementation}). The gap $\tfrac{\bar\sigma_4}{\bar\sigma_1}$ is smallest for the critical catenoid, where it equals $t^2=1.4392$ with $t\tanh t=1$, and generally increases with $b$. $\mathcal R_{\rm sph}$ and $\mathcal C_X$ are below $4\times10^{-10}$ and $2\times10^{-8}$ for $b\le42$ and reach $6\times10^{-6}$ and $4\times10^{-4}$ for $b=60$--$252$. For the geodesic icosahedra with $362\le b\le2252$ they are below $6\times10^{-9}$ and $3\times10^{-6}$. The area density satisfies $\min_\Omega J>0$ in every row with $b\le252$, in agreement with \cite[Proposition 8.1]{Fraser_2015}. Its value depends on the M\"obius gauge of $\Omega$, and it is not evaluated for $b\ge362$. All values are in the repository \cite{Chen_Osting_genus0_code}. Every row gives a lower bound $\sigma^\star(0,b)\ge\bar\sigma_1^\star$, up to the discretization error of Section~\ref{ssec:validation}, whether or not it is critical. Since $\sigma^\star(0,b)$ is twice the largest area of an embedded genus-zero free boundary minimal surface with $b$ boundary components \cite{Fraser_2015,Chodosh_Gianocca_2026} (Section~\ref{sec:intro}), Table~\ref{tab:best} also bounds these areas from below, and its surfaces are our candidates for the area maximizers conjectured in \cite[Conjecture 7.10]{Carlotto_Schulz_Wiygul_2025}.

\subsection{Large \texorpdfstring{$b$}{b}}\label{ssec:large}
\begin{figure}[tbp]\centering
\includegraphics[width=0.42\textwidth]{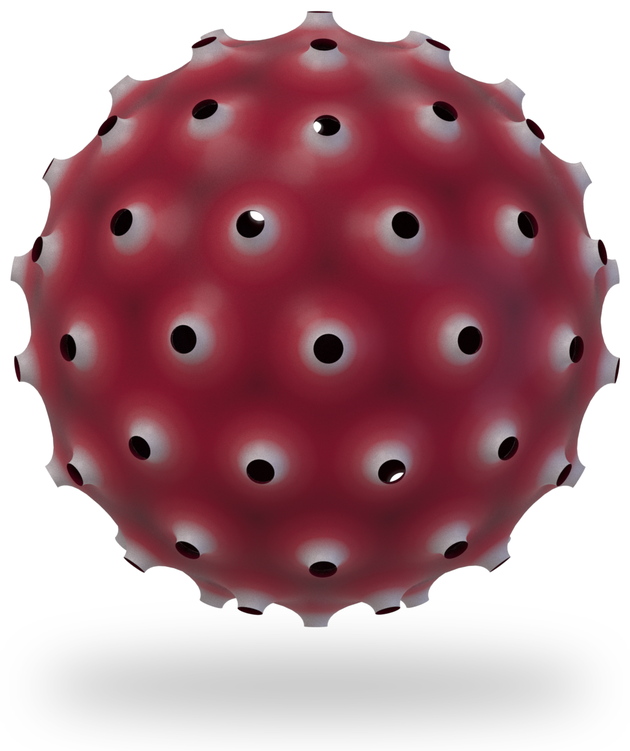}\hfill\includegraphics[width=0.56\textwidth]{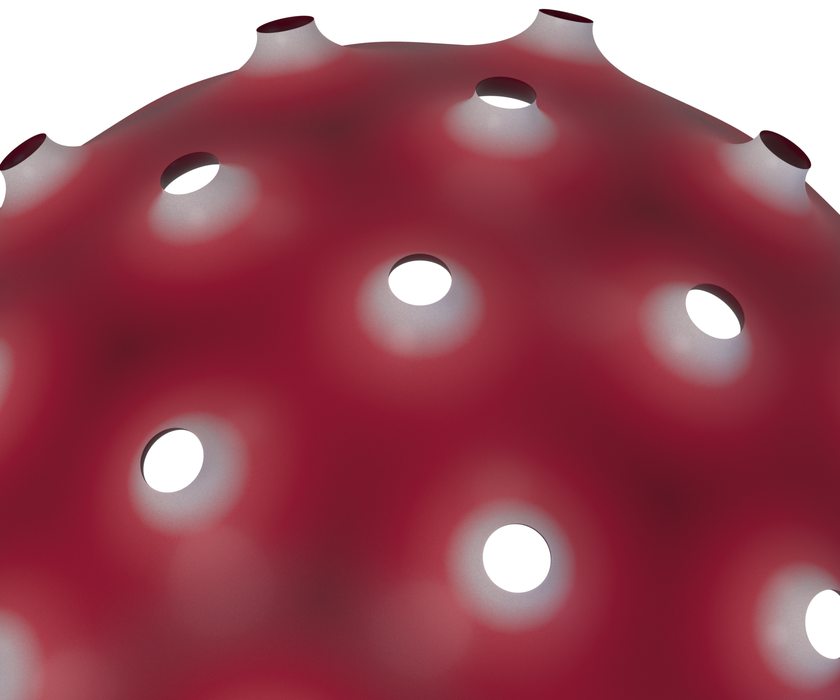}
\caption{The surface with $b=92$, the geodesic icosahedron with $\nu=3$ (Section~\ref{ssec:large}, Table~\ref{tab:best}). Left: the whole surface. Right: a close-up. The sheet lies just inside $\partial\ball$, and each boundary circle is the rim of a short neck to $\partial\ball$. The geometry is not exaggerated, and the color shows $|X|$ as in Figure~\ref{fig:headline}.}\label{fig:closeup}
\end{figure}

With icosahedral symmetry the cost grows with the number of orbits rather than with $b$. The reduced solver of Section~\ref{ssec:implementation} thus reaches the geodesic icosahedron with $\nu=15$, where $b=2252$ and $\bar\sigma_1^\star=24.97$ (Figures~\ref{fig:best-b42}--\ref{fig:best-b2252}). The surfaces are a sheet just inside $\partial\ball$ with $b$ short necks to it (Figure~\ref{fig:headline}), and Figure~\ref{fig:closeup} shows the surface with $b=92$ and a close-up of its necks. This is the picture behind the homogenization of \cite{Girouard_Lagace_2021} and the convergence theorem of \cite{Karpukhin_Stern_2021}. The cap radii scale like $b^{-1}$, with $b\,r$ between $1.7$ and $2.0$ from $b=12$ to $b=2252$.

Karpukhin and Stern proved that $8\pi-\sigma^\star(0,b)$ lies between two multiples of $b^{-1}\log b$ \cite[Theorem 1.8]{Karpukhin_Stern_2021}. Below we use the monopole and dipole approximations of Section~\ref{ssec:approx}, evaluate them formally for equal caps at the Fekete configurations, and compare with the computed values.

\emph{The approximations.} On the scale of their spacing, small boundary circles act as point sources. The monopole approximation of Section~\ref{ssec:approx} makes this precise. It keeps one number per cap, the flux $q_\ell=\bar\sigma_1w_\ell\,\bar u_\ell$, with $\bar u_\ell$ the mean of the eigenfunction on $\Gamma_\ell$, and gives the value $\bar\sigma^{(0)}_1$ of \eqref{eq:models}. The dipole approximation adds the response of each cap to the gradient of the field of the other caps, with the factor $\kappa_{\ell,1}=\tfrac{1+\varepsilon_\ell}{1-\varepsilon_\ell}$ of Corollary~\ref{cor:const}, and gives the value $\bar\sigma^{(1)}_1$ of \eqref{eq:models}. For $\varepsilon_\ell=0$ the response is that of an insulating hole. To leading order in the area fraction $f=\tfrac14\sum_\ell r_\ell^2$ of the caps, the dipoles divide the value of the monopole approximation by $1+2f$. Both approximations make two simplifications: they drop the higher multipoles of the exact system \eqref{eq:mp-steklov}, and they replace the weight on each circle by its mean. The first simplification has no error bound here, and the second only the bound of Section~\ref{ssec:approx}, which is far from sharp. Section~\ref{ssec:catalog} measures the resulting errors on every configuration of Table~\ref{tab:catalog}.

\emph{Equal caps at the Fekete configurations.} Let the caps have equal weights $w_\ell=\tfrac1b$ and equal radii $r=\tfrac{\tilde r}b$, so that $R_\ell\approx r$ and $\tfrac1{\bar\sigma^{(0)}_1}=\tfrac1b\lambda_{\max}(\mathsf P\mathsf G)$. Since the first eigenfunctions are close to the coordinate functions, we evaluate the Rayleigh quotient of $\mathsf G$ on the three vectors $(c_\ell\cdot e_i)_\ell$ and average. Using $c_m\cdot c_\ell=1-\tfrac12|c_m-c_\ell|^2$ and $\sum_i\sum_\ell(c_\ell\cdot e_i)^2=b$, we get
\begin{equation}\label{eq:log-energy}
\frac1{\bar\sigma^{(0)}_1}\approx\frac1{b^2}\Bigl(\frac b{2\pi}\log\frac1r+\frac1\pi L(\mathbf c)+\frac1{4\pi}\sum_{m\ne\ell}|c_m-c_\ell|^2\log|c_m-c_\ell|\Bigr),\qquad L(\mathbf c):=-\sum_{\ell<\ell'}\log|c_\ell-c_{\ell'}| .
\end{equation}
The centers $\mathbf c=(c_1,\dots,c_b)$ thus enter through their logarithmic energy $L$ and through a sum of a smooth kernel. For the Fekete configurations, the minimizers of $L$, one has $L(\mathbf c)=\bigl(\tfrac14-\tfrac12\log2\bigr)b^2-\tfrac14\,b\log b+\tfrac12C_{\log}\,b+o(b)$, where $C_{\log}$ is the constant in the expansion of the minimal logarithmic energy summed over ordered pairs \cite[Theorem 1.5]{Betermin_Sandier_2018}. B\'etermin and Sandier prove $C_{\log}\le C_{\rm BHS}=2\log2+\tfrac12\log\tfrac23+3\log\tfrac{\sqrt\pi}{\Gamma(1/3)}\approx-0.0556$, and Brauchart, Hardin and Saff conjecture equality \cite{Brauchart_Hardin_Saff_2012}. For equidistributed centers, the last sum in \eqref{eq:log-energy} is $b^2\bigl(2\log2-\tfrac12\bigr)$ to leading order, the mean of the kernel over $\sphere\times\sphere$. The terms of order $b^2$ then add up to $\tfrac{b^2}{8\pi}$, which gives the limit $8\pi$. Including the factor $1+2f=1+\tfrac{\tilde r^2}{2b}$ of the dipoles, the terms of order $b\log b$ and $b$ give
\begin{equation}\label{eq:rate-full}
\bar\sigma_1^\star\approx\frac{8\pi}{1+b^{-1}\bigl(2\log b-\frac{d(\tilde r)}{8\pi}\bigr)},\qquad d(\tilde r)=32\pi\bigl(\log\tilde r-\tfrac18\tilde r^2-C_{\log}\bigr),
\end{equation}
and to first order
\begin{equation}\label{eq:rate}
8\pi-\bar\sigma_1^\star(b)\approx\frac{16\pi\log b-d}{b} .
\end{equation}
The coefficient $16\pi\approx50.3$ does not depend on $\tilde r$. It combines the self term $\log\tfrac1r$ with the term $-\tfrac14\,b\log b$ of the energy. The constant does depend on $\tilde r$. Maximizing $d(\tilde r)$ gives $\tilde r=2$, that is, $r^\star=2b^{-1}$, and
\begin{equation}\label{eq:d}
d=32\pi\bigl(\log2-\tfrac12\bigr)-32\pi\,C_{\log}\ge25.0,
\end{equation}
with equality if the conjecture holds. The lower bound of \cite[Section 4]{Karpukhin_Stern_2021} is proved with caps of radius $b^{-3/2}$, with the remark that $b^{-1}$ might be optimal.

\emph{Comparison with the computed values.} The cap radii of the geodesic icosahedra approach the optimal radius $r^\star=2b^{-1}$ of the dipole approximation, with $b\,r=1.98$, $1.99$ and $2.00$ at $b=252$, $1002$ and $2252$. Along the geodesic icosahedra, the slope of $b\,(8\pi-\bar\sigma_1^\star)$ against $\log b$ between consecutive $b$ decreases steadily toward $16\pi$, from $55.3$ between $b=92$ and $162$ to $51.8$ between $b=1962$ and $2252$. With the coefficient fixed, a fit of $d$ over $162\le b\le2252$ gives $d\approx30$, with relative residuals in $8\pi-\bar\sigma_1^\star$ of at most $2.6\%$, and a fit at each $b$ drifts from $35$ at $b=162$ to $27$ at $b=2252$. This fitted value differs from the value \eqref{eq:d} of the evaluation because \eqref{eq:rate} is only the first-order expansion of \eqref{eq:rate-full}, whose next order is still large at these $b$. With $d$ from \eqref{eq:d}, the form \eqref{eq:rate-full} has no fitted parameter. It reproduces $\bar\sigma_1^\star$ to $1.3\times10^{-2}$ at $b=162$ and $2\times10^{-4}$ at $b=2252$, and the value of $d$ that it fits rises from $22.7$ at $b=162$ to $24.6$ at $b=2252$. It is the dotted curve in Figure~\ref{fig:families}. A power law does not fit: the local exponent $-\tfrac{\Delta\log(8\pi-\bar\sigma_1^\star)}{\Delta\log b}$ increases steadily along the geodesic icosahedra, from $0.73$ to $0.85$. By \eqref{eq:rate}, the deficit falls below $0.1$ only beyond $b\approx4\times10^3$. Since $\sigma^\star(0,b)\ge\bar\sigma_1^\star(b)$, the computed deficits bound the true ones from above. Together with the multipole approximations, they suggest that the limit asked for in \cite[Open Question 2]{Karpukhin_Stern_2021} is $16\pi$.

\subsection{Map of the zoo}\label{ssec:catalog}\label{ssec:thomson}
Table~\ref{tab:catalog} collects every critical point computed. The repository \cite{Chen_Osting_genus0_code} contains one folder for each configuration of Table~\ref{tab:catalog}, linked from its name there and in Tables~\ref{tab:zoo} and \ref{tab:best}, with its configuration parameters, spectrum, optimality certificate, and surface meshes and plots. Table~\ref{tab:map} orders the ring families and their limits, Figures~\ref{fig:type-polyhedral}--\ref{fig:type-stacked} show the surfaces of each type, and Figure~\ref{fig:families} plots the values against $b$ by family.

\emph{Logarithmic energy of the centers.} Oudet, Kao and Osting \cite{Kao_2020} observed that the maximizers for $b\le20$ resemble the Thomson configurations \cite{Wales_CCD}, the minimizers of the Coulomb energy. The monopole approximation \eqref{eq:models} of Section~\ref{ssec:approx} explains a relation of this kind and shows that the relevant energy is the logarithmic one. In it the caps interact through the matrix $\mathsf G$ of \eqref{eq:G-matrix}, whose off-diagonal entries are the Green's function of $\sphere$, logarithmic in the chordal distance between the centers. For equal caps, the Rayleigh quotient \eqref{eq:log-energy} reduces the dependence on the centers to their logarithmic energy, and a smaller energy gives a larger $\bar\sigma_1$. The centers of the boundary curves are the centroids of $X^\star(\Gamma_\ell)$. This explains why the best configurations lie close to the Fekete configurations without coinciding with them, and why the class of the Fekete configuration gives the best value found for most $b$ from $19$ to $40$ (Section~\ref{ssec:best}).

\emph{Accuracy of the multipole approximations.} The right panel of Figure~\ref{fig:families} tests the two approximations \eqref{eq:models} on every configuration of Table~\ref{tab:catalog}, with its own centers, radii and mean weight on each circle. The monopole approximation overestimates $\bar\sigma_1^\star$ by $5$ to $20\%$, nearly uniformly at fixed $b$, so it still orders the critical points of each $b$ except for near-ties. With the dipoles, $\bar\sigma^{(1)}_1$ of \eqref{eq:models}, the median relative error falls to $6\times10^{-4}$ for $13\le b\le20$, and to $10^{-12}$ for the geodesic icosahedron with $b=2252$. Because both approximations replace the weight on each circle by its mean, they are least accurate where the weight varies most, with errors up to $5\%$ for the equatorial polygons.

\emph{Morse index.} The Morse index of $\bar\sigma_1$ (Section~\ref{ssec:morse}) is the index column of Table~\ref{tab:catalog}. It is not the index of the surface as a critical point of area. Where computed, it tends to grow as the value decreases, from $0$ at the maximizer to $b-3$ at the equatorial $b$-gon. For the pyramids and the bipyramids with a ring of $j\ge6$ caps it is $j-3$. For $b\le12$ the maximizer is the only local maximum in Table~\ref{tab:catalog}. At $b=16$ and $22$ there is a second one, $0.006$ and $0.0015$ below the best. Every maximizer found for $16\le b\le40$ whose index is computed is a local maximum, except at $b=36$, where the best configuration found has index $1$. At $b=92$ and $162$ the geodesic classes also contain saddles, of index $50$ and $52$.

The indices and symmetry groups of Table~\ref{tab:catalog} satisfy the count \eqref{eq:morse-count} for small $b$. The configurations of $b$ labeled caps, up to M\"obius maps, form a space with Euler characteristic $(-1)^{b-3}(b-3)!$. On this space, $\bar\sigma_1$, maximized over the weights, is a function whose critical points correspond to free boundary minimal surfaces. A surface $\Sigma$ gives $\tfrac{2\,b!}{|G_\Sigma|}$ critical points, which differ by the labels of the caps. Suppose that $\bar\sigma_1$ is smooth with nondegenerate critical points, and that it decreases as a configuration degenerates, that is, as the radius of a cap or the distance between two caps tends to zero. Then, by the Poincar\'e--Hopf theorem, the critical points, each counted with the sign $(-1)^{\operatorname{ind}(\Sigma)}$, add up to the Euler characteristic, which is \eqref{eq:morse-count}. These hypotheses on $\bar\sigma_1$ are not proved, so \eqref{eq:morse-count} is a prediction. With the indices and groups of Table~\ref{tab:catalog}, it holds for $3\le b\le10$. For $b=4$, the tetrahedron and the equatorial $4$-gon give $\tfrac1{24}-\tfrac1{16}=-\tfrac1{48}$. At $b=9$ and $10$, seven of the fifteen and fifteen of the twenty-six surfaces were found because the count failed without them. Their starting points were the critical points of the logarithmic energy of $b$ points on the sphere. For every $b\ge11$ the count fails, so critical points are missing there. Table~\ref{tab:catalog} has $1$, $2$, $3$, $4$, $7$, $13$, $15$ and $26$ critical points for $b=3,\dots,10$.

The index bound of Conjecture~\ref{conj:morse} holds in every row of Table~\ref{tab:catalog} with a computed index. Here $b-3$ is the complex dimension of the space of $b$ points on the sphere up to M\"obius maps, and its cohomology vanishes above that degree.
At the equatorial $b$-gon with $4\le b\le12$, the $b-3$ increasing directions are odd under the reflection in the equatorial plane. For $b=3$ the conjecture would imply that the equatorial $3$-gon is the only critical point. Every term of \eqref{eq:morse-count} would be positive and at least $\tfrac1{12}$, the right-hand side.

\begin{figure}[tbp]\centering
\input{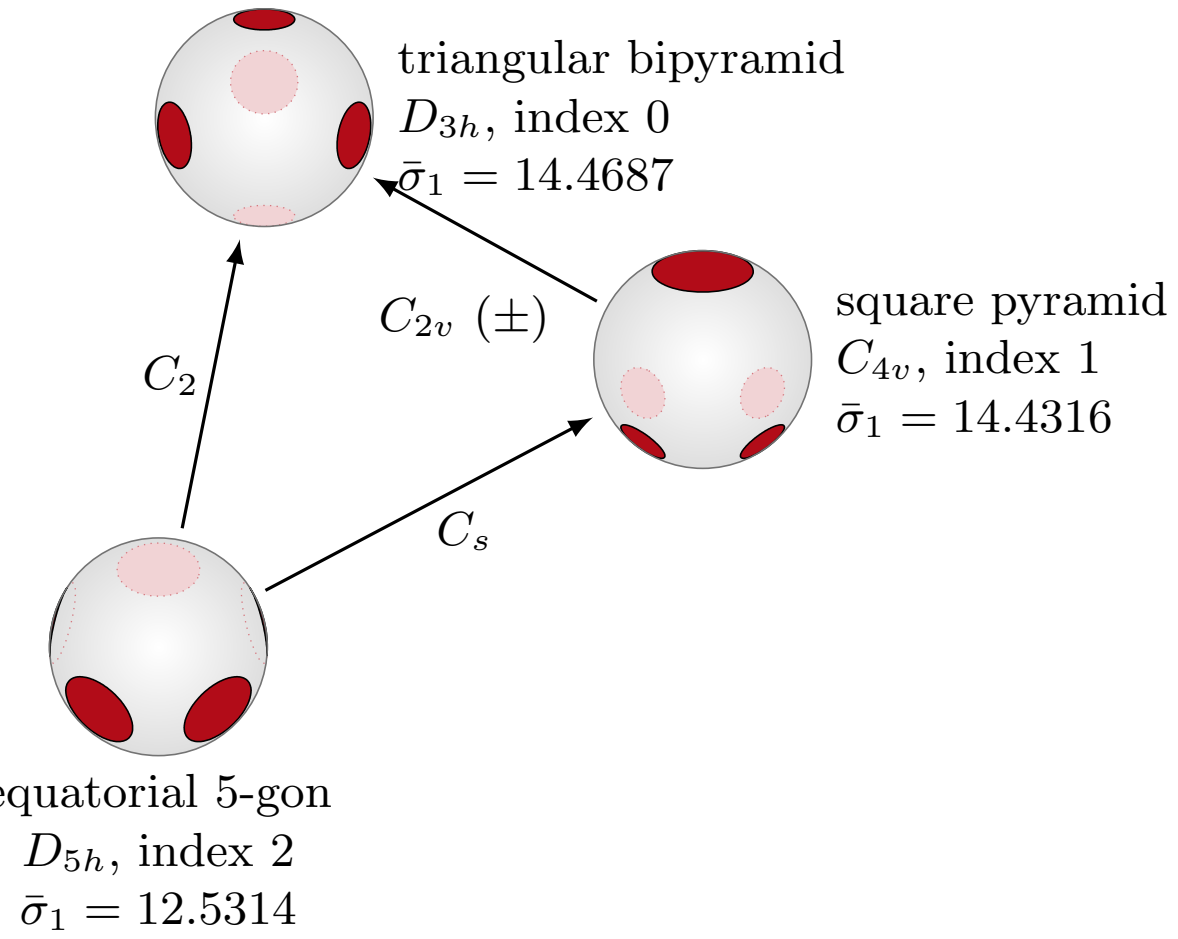}
\caption{Connections between the critical points with $b=5$ (Section~\ref{ssec:catalog}). Each sphere shows the caps of a critical configuration of Table~\ref{tab:catalog}, with the caps on the far side pale. An arrow is an ascent path of $\bar\sigma_1$, labeled by the symmetry group that it keeps. Both signs of the unstable direction of the square pyramid lead to the triangular bipyramid ($\pm$).}
\label{fig:connections}
\end{figure}

\emph{Connections.} For $b=4$, $5$ and $6$ we follow the ascent paths of $\bar\sigma_1$ that leave each saddle along its unstable directions and keep their symmetry. Figure~\ref{fig:connections} shows $b=5$. The equatorial $5$-gon connects to the square pyramid and to the triangular bipyramid. The square pyramid connects to the triangular bipyramid on both sides of its unstable direction. The two ends differ by which pair of opposite base caps becomes the poles.

\section{Discussion}\label{sec:discussion}
We computed critical points of $\bar\sigma_1$ in genus zero with a layer potential method on the sphere minus round caps, whose eigenvalues, derivatives and certificates (Theorems~\ref{thm:bie}, \ref{thm:derivatives} and~\ref{thm:certificate}) turn each computed critical point into a free boundary minimal surface. This gives a zoo of $268$ certified surfaces with up to $2252$ boundary components (Table~\ref{tab:catalog}). The results raise the following questions.

\emph{Limit surfaces.} The computed families suggest that the limits other than the wheel and the stacks of \cite{Kapouleas_Zou_2021} can be desingularized as well: the drum, the cake stand, and the drums and stacks with horns (Section~\ref{sec:fcms}). A proof would also give the convergence of the families, which we only observe. For the class maximizers, the convergence of the values might also be proved without gluing. Filling in the caps and the upper semicontinuity of eigenvalues of measures \cite{Girouard_Karpukhin_Lagace_2021} bound the limit from above by the limit problem of Section~\ref{ssec:families}. The lower bound would need a symmetric version of the homogenized test metrics of \cite{Girouard_Lagace_2021,Girouard_Karpukhin_Lagace_2021}. The networks with horns raise a spectral question too. For the stacks closed by disks, \cite[Theorem 1.10]{Karpukhin_Noutchegueme_2025} shows that $\sigma=1$ is the first eigenvalue of the transmission problem and that the stack maximizes $\bar\sigma_1$ among rotationally symmetric metrics. For the cake stand, the wheel and the other networks with horns we checked this only numerically (Section~\ref{ssec:ringlimits}), and we found one network of each type.

\emph{Large $b$.} The main analytic question is to make the analysis of Section~\ref{ssec:large} rigorous. This has two parts. The first is an expansion of $\bar\sigma_1$ for a sphere with $b$ small round caps, with a bound on the remainder. By Proposition~\ref{prop:multipole} the eigenvalue problem is an exact system for the charges and multipoles of the caps, and the multipole approximations keep only the charges and the dipoles. The remainder therefore comes from the multipoles of order two and higher and from the variation of the weight on each circle. For holes without weights and with balanced radii, such an expansion is proved in \cite{Huang_Venkatraman_2026} (Remark~\ref{rem:hv}). The accuracy of the dipole approximation (Section~\ref{ssec:thomson}) suggests that the charges and the dipoles carry the leading terms. Such an expansion would give $\sigma^\star(0,b)\ge8\pi-b^{-1}(16\pi\log b-d)$ up to lower order. The second part is the matching upper bound, which must hold for every metric. The upper bound of \cite[Theorem 1.8]{Karpukhin_Stern_2021} follows from a stability inequality that bounds the distance in $(W^{1,2})^*$ between the boundary measure and the area measure by the deficit $8\pi-\bar\sigma_1$. The constant $16\pi$ would need this inequality with its sharp constant. Together, the two parts would answer \cite[Open Question 2]{Karpukhin_Stern_2021}.

\emph{The zoo.} Unrestricted runs from the Fekete configurations at $b\ge19$, and from every critical point of the zoo, would test whether the best configurations found are the maximizers. Existence of the class maximizers is open for the pyramids, the prisms and antiprisms without polar caps, and the chiral class at $b=72$ (Remark~\ref{rem:existence}). Only two surfaces of the zoo, at $b=25$ and $33$, have trivial symmetry group (Section~\ref{ssec:best}). Whether such surfaces exist with fewer boundary components is open. Schulz conjectures one with $61$ \cite{Schulz_gallery}. The symmetry reduction of Section~\ref{ssec:implementation} applies to every class, not only the icosahedral ones, and to the higher eigenvalues $\bar\sigma_k$, as in \cite[\S5.2]{Kao_2020}. Finally, no closed form is known beyond the disk and the critical catenoid. The surfaces computed here are accurate to about twelve digits and could serve to test candidate Weierstrass data.

\subsection*{Data availability}
All configurations, certificates and surfaces are in the repository \cite{Chen_Osting_genus0_code}.

\subsection*{Acknowledgments}
The authors used Claude (Anthropic) to assist with editing the text, verifying references, checking derivations, developing parts of the numerical code, and exploring symmetry configurations. All AI-assisted text and code were checked by the authors, who take full responsibility for the content of the paper.

\FloatBarrier
\printbibliography
\clearpage
\section*{Appendix: figures and tables}
\begin{table}[H]\centering\scriptsize\setlength{\tabcolsep}{2.5pt}\renewcommand{\arraystretch}{0.72}
\begin{tabular}{@{}C{4.4cm}lC{2.6cm}C{1.95cm}C{3.17cm}@{}}\toprule
class & group & $b$ & geometric parameters & weight modes per orbit\\\midrule
\multicolumn{5}{@{}l}{\textbf{Polyhedral configurations}}\\
tetrahedron & $T_d$ & $4$ & $1$ & $\cos3n\theta$\\
octahedron, cube & $O_h$ & $6$, $8$ & $1$ & $\cos4n\theta$, $\cos3n\theta$\\
icosahedron, dodecahedron, icosi\-dodecahedron & $I_h$ & $12$, $20$, $30$ & $1$ & $\cos5n\theta$, $\cos3n\theta$, $\cos2n\theta$\\
icosahedron + dodecahedron & $I_h$ & $32$ & $3$ & $\cos5n\theta$, $\cos3n\theta$\\
$60$-orbit on the vertex-to-vertex arcs of the icosahedron, free along the arc & $I_h$ & $60$ & $2$ & $\cos n\theta$\\
icosahedron + generic $60$-orbit & $I$ & $72$ & $5$ & $\cos5n\theta$; $\cos n\theta$, $\sin n\theta$\\
geodesic icosahedron of frequency $\nu$: one orbit per vertex of a face & $I_h$ & $10\nu^2+2$, $\nu=2,\dots,15$ & $2o-1+p$ & $\cos qn\theta$, $q$  the number of rotations fixing a center of the orbit\\
\addlinespace[2pt]\multicolumn{5}{@{}l}{\textbf{The equatorial family}}\\
equatorial $b$-gon & $D_{bh}$ & $3,\dots,16$, $18$, $20$ & $1$ & $\cos2n\theta$\\
\addlinespace[2pt]\multicolumn{5}{@{}l}{\textbf{Pyramids and bipyramids}}\\
pyramid: ring of $j$ + one pole & $C_{jv}$ & $j+1$, $j=3,\dots,16$, $20$ & $4$ & ring $\cos n\theta$, pole $\cos jn\theta$\\
bipyramid: equatorial $j$-gon + two poles & $D_{jh}$ & $j+2$, $j=3,\dots,10$, $16$, $18$, $24$, $30$ & $3$ & ring $\cos2n\theta$, poles $\cos jn\theta$\\
\addlinespace[2pt]\multicolumn{5}{@{}l}{\textbf{Prisms and antiprisms}}\\
prism / antiprism: two rings of $j$, aligned / rotated by $\tfrac{\pi}{j}$ & $D_{jh}$ / $D_{jd}$ & $2j$, $j=3,\dots,9$ / $j=3,\dots,16$, $20$ & $2$ & $\cos n\theta$\\
prism / antiprism + two poles & $D_{jh}$ / $D_{jd}$ & $2j+2$, $j=3,\dots,9$ / $j=3,\dots,16$, $19$ & $4$ & rings $\cos n\theta$, poles $\cos jn\theta$\\
capped prism / antiprism: rings as separate orbits + one pole & $C_{jv}$ & $2j+1$, $j=3,\dots,19$ & $7$ & rings $\cos n\theta$, pole $\cos jn\theta$\\
\addlinespace[2pt]\multicolumn{5}{@{}l}{\textbf{Stacked rings}}\\
$m$ rings of two caps, alternating or independent, with or without polar caps & $D_{2d}$, $D_{2h}$, $C_{2v}$ & $2m$, $2m+1$, $2m+2$; $m=2,3,4$ & $3m-1$ ($+2$ with poles) & $\cos n\theta$\\
$m$ free pairs about one half-turn, + one or two poles ($m=3$) or none ($m=4$) & $C_2$ & $7$, $8$, $8$ & $11$, $13$, $13$ & $\cos n\theta$, $\sin n\theta$; poles $\cos2n\theta$, $\sin2n\theta$\\
$m$ rings of $j$, consecutive rings rotated by $\tfrac{\pi}{j}$, with or without two poles & $D_{jh}$ ($m$ odd), $D_{jd}$ ($m$ even) & $mj$, $mj+2$; $m=3,\dots,6$, $j=3,\dots,6$; $m=3$, $j=7,9,11$, $13$; $m=5$, $j=7$ & $2\lceil\tfrac m2\rceil+\lfloor\tfrac m2\rfloor-1$ ($+2$ with poles) & $\cos n\theta$, equatorial ring $\cos2n\theta$, poles $\cos jn\theta$\\
three rings of $j$ as separate orbits, the middle one rotated by $\tfrac{\pi}{j}$, + one pole & $C_{jv}$ & $3j+1$, $j=3,\dots,9$, $11$, $13$, $16$, $20$ & $9$ & rings $\cos n\theta$, pole $\cos jn\theta$\\
four caps in a generic orbit of $C_{2v}$ + a pair on a mirror plane, + one pole, two poles or a second pair & $C_{2v}$ & $7$, $8$, $8$ & $7$, $9$, $8$ & four caps $\cos n\theta$, $\sin n\theta$; pairs $\cos n\theta$; poles $\cos2n\theta$\\
\addlinespace[2pt]\multicolumn{5}{@{}l}{\textbf{Classes of the Fekete configurations}}\\
Fekete class: orbits of the Fekete configuration under its symmetry group $G$ & $G$ & $19,\dots,31$, $33,\dots,37$, $39$, $40$ & as in Section~\ref{ssec:classes} & stabilizer-invariant modes\\
class of the unrestricted surface, from its centers & $T_d$, $C_{2v}$ & $16$, $17$ & as above & as above\\
\addlinespace[2pt]\multicolumn{5}{@{}l}{\textbf{Unrestricted}}\\
every cap its own orbit, normal form \eqref{eq:normalization} & $C_1$ & any & $4b-7$ & $\cos n\theta$, $\sin n\theta$\\
\bottomrule
\end{tabular}
\caption{The symmetry classes by family, with their geometric parameters and weight modes (Section~\ref{ssec:classes}; the families are discussed in Section~\ref{ssec:families}). For the geodesic icosahedra, $o$ is the number of orbits and $p$ the number of free positions, one per $60$-orbit and two per generic $120$-orbit ($\nu\ge6$). The index $n$ runs to $n_{\max}$ (Section~\ref{ssec:optimizer}).}\label{tab:classes}
\end{table}

\begin{figure}[p]\centering
\setlength{\panelw}{\dimexpr(\textwidth-6\tabcolsep)*100/400\relax}
\setlength{\tabcolsep}{1pt}
\begin{tabular}{@{}cccc@{}}
\paneltitle{$b=4$, $\bar\sigma_1^\star=13.67$\\\zoo{b0004/tetrahedron}{tetrahedron}} & \paneltitle{$b=6$, $\bar\sigma_1^\star=15.42$\\\zoo{b0006/octahedron}{octahedron}} & \paneltitle{$b=8$, $\bar\sigma_1^\star=16.42$\\\zoo{b0008/cube}{cube}} & \paneltitle{$b=12$, $\bar\sigma_1^\star=18.07$\\\zoo{b0012/icosahedron}{icosahedron}}\\[1pt]
\includegraphics[width=\panelw]{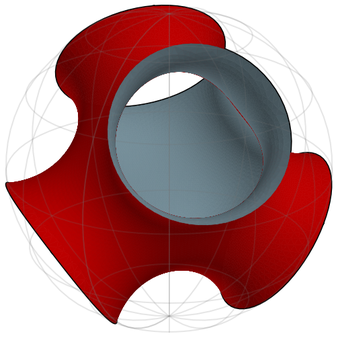} & \includegraphics[width=\panelw]{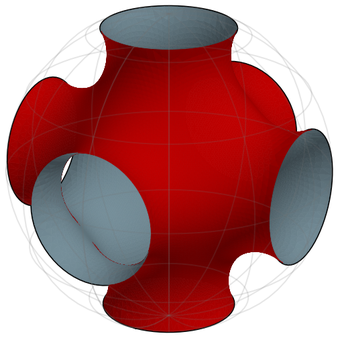} & \includegraphics[width=\panelw]{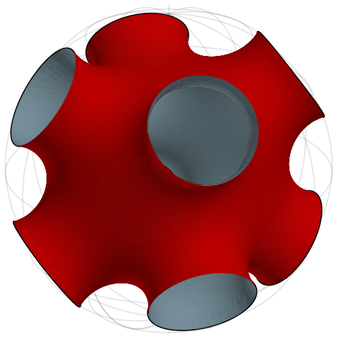} & \includegraphics[width=\panelw]{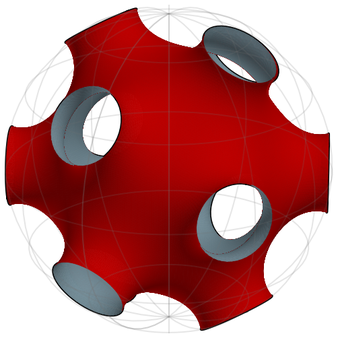}\\[-1pt]
\includegraphics[width=\panelw]{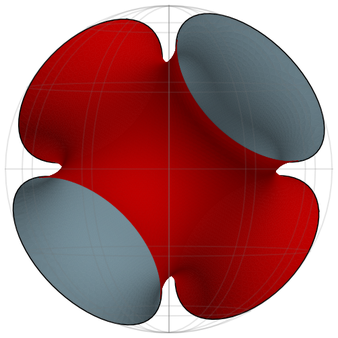} & \includegraphics[width=\panelw]{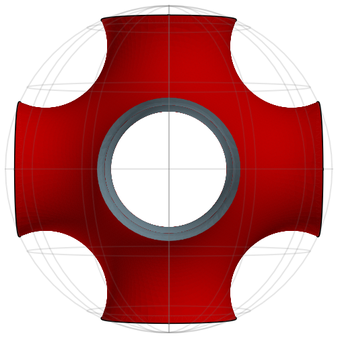} & \includegraphics[width=\panelw]{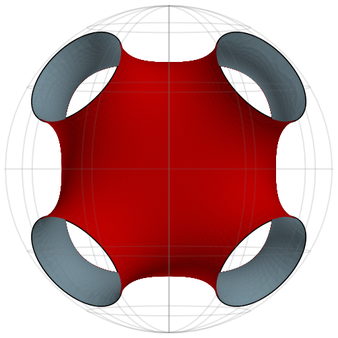} & \includegraphics[width=\panelw]{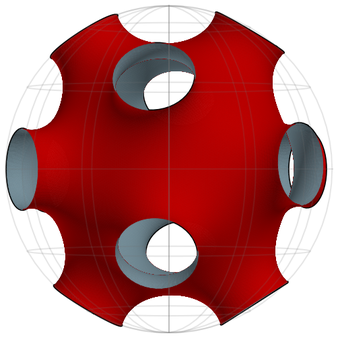}\\[10pt]
\paneltitle{$b=20$, $\bar\sigma_1^\star=19.57$\\\zoo{b0020/dodecahedron}{dodecahedron}} & \paneltitle{$b=12$, $\bar\sigma_1^\star=17.98$\\\zoo{b0012/three_stacked_squares}{cuboctahedron}} & \paneltitle{$b=30$, $\bar\sigma_1^\star=20.73$\\\zoo{b0030/icosidodecahedron}{icosidodecahedron}} & \paneltitle{$b=32$, $\bar\sigma_1^\star=21.01$\\\zoo{b0032/icosahedron_plus_dodecahedron}{icosahedron + dodecahedron}}\\[1pt]
\includegraphics[width=\panelw]{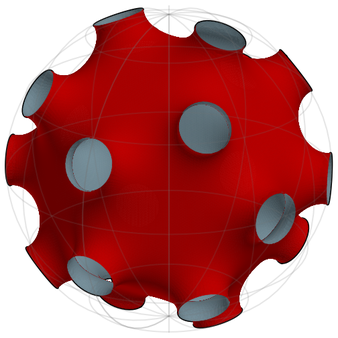} & \includegraphics[width=\panelw]{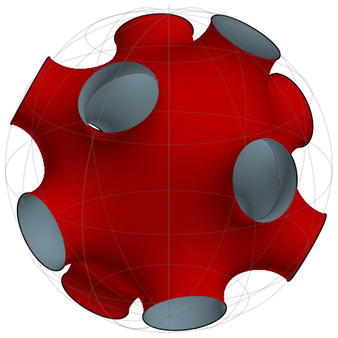} & \includegraphics[width=\panelw]{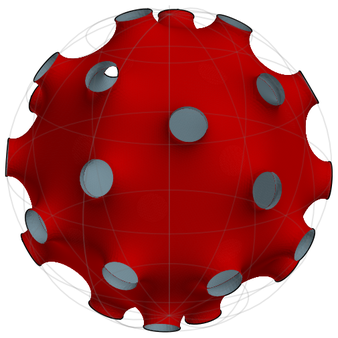} & \includegraphics[width=\panelw]{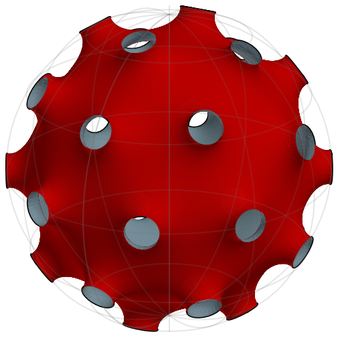}\\[-1pt]
\includegraphics[width=\panelw]{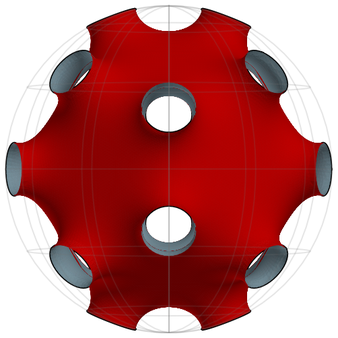} & \includegraphics[width=\panelw]{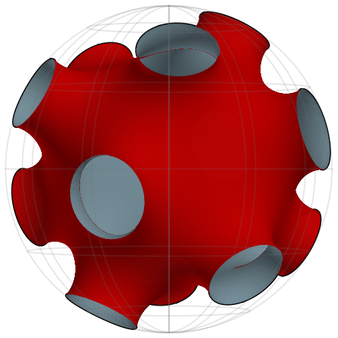} & \includegraphics[width=\panelw]{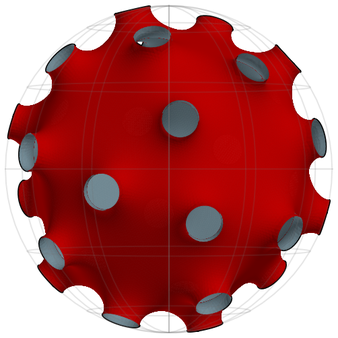} & \includegraphics[width=\panelw]{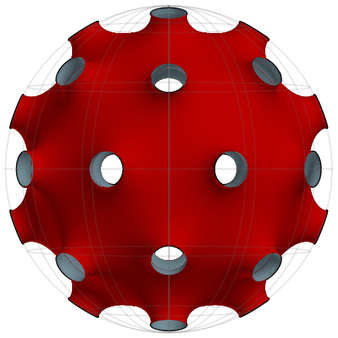}
\end{tabular}

\caption{Polyhedral configurations (Section~\ref{ssec:platonic}): the surfaces $X^\star(\overline\Omega)\subset\ball$ built from the certificates, oblique and side views, boundary circles in black. Values in Table~\ref{tab:catalog}.}\label{fig:type-polyhedral}
\end{figure}

\begin{figure}[t]\centering
\setlength{\panelw}{\dimexpr(\textwidth-6\tabcolsep)*100/400\relax}
\setlength{\tabcolsep}{1pt}
\begin{tabular}{@{}cccc@{}}
\paneltitle{$b=3$, $\bar\sigma_1^\star=12.01$\\\zoo{b0003/equatorial_3_gon}{triangle}} & \paneltitle{$b=4$, $\bar\sigma_1^\star=12.43$\\\zoo{b0004/equatorial_4_gon}{square}} & \paneltitle{$b=5$, $\bar\sigma_1^\star=12.53$\\\zoo{b0005/equatorial_5_gon}{pentagon}} & \paneltitle{$b=6$, $\bar\sigma_1^\star=12.56$\\\zoo{b0006/equatorial_6_gon}{hexagon}}\\[1pt]
\includegraphics[width=\panelw]{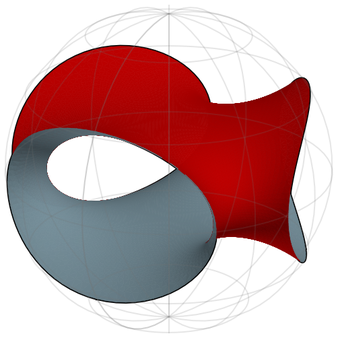} & \includegraphics[width=\panelw]{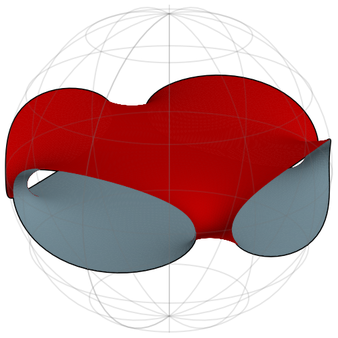} & \includegraphics[width=\panelw]{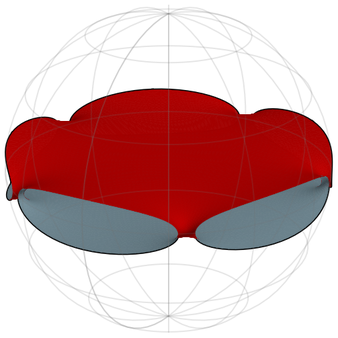} & \includegraphics[width=\panelw]{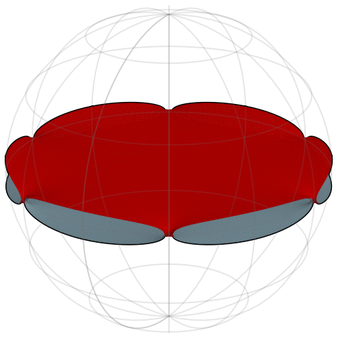}\\[-1pt]
\includegraphics[width=\panelw]{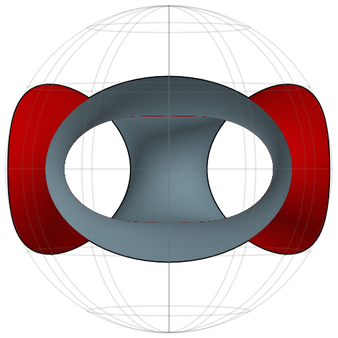} & \includegraphics[width=\panelw]{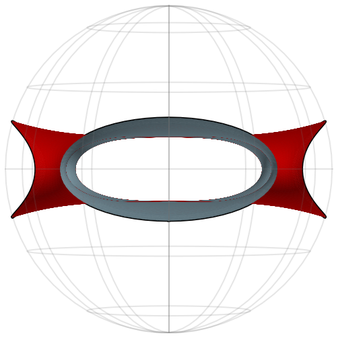} & \includegraphics[width=\panelw]{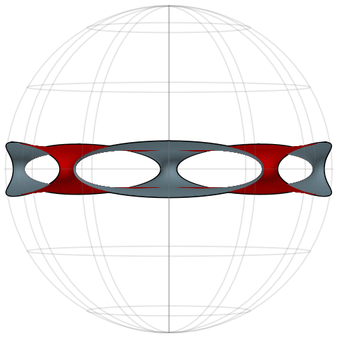} & \includegraphics[width=\panelw]{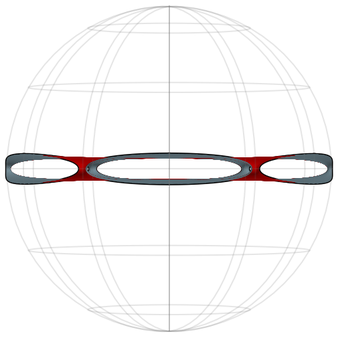}\\[10pt]
\paneltitle{$b=7$, $\bar\sigma_1^\star=12.56$\\\zoo{b0007/equatorial_7_gon}{heptagon}} & \paneltitle{$b=8$, $\bar\sigma_1^\star=12.57$\\\zoo{b0008/equatorial_8_gon}{octagon}} & \paneltitle{$b=9$, $\bar\sigma_1^\star=12.57$\\\zoo{b0009/equatorial_9_gon}{9-gon}} & \paneltitle{$b=10$, $\bar\sigma_1^\star=12.57$\\\zoo{b0010/equatorial_10_gon}{10-gon}}\\[1pt]
\includegraphics[width=\panelw]{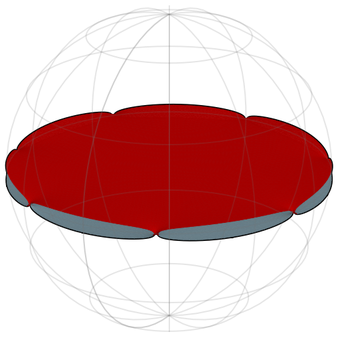} & \includegraphics[width=\panelw]{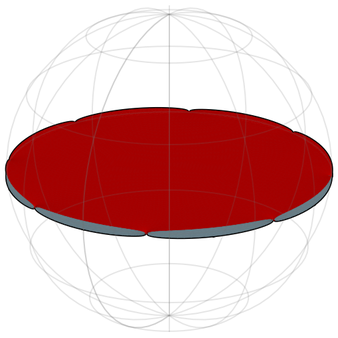} & \includegraphics[width=\panelw]{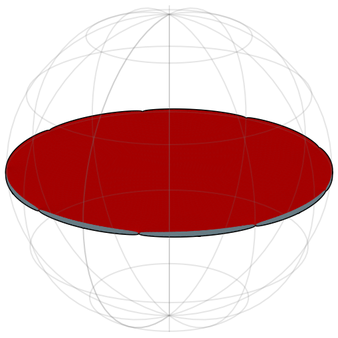} & \includegraphics[width=\panelw]{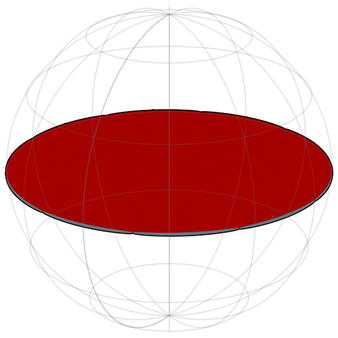}\\[-1pt]
\includegraphics[width=\panelw]{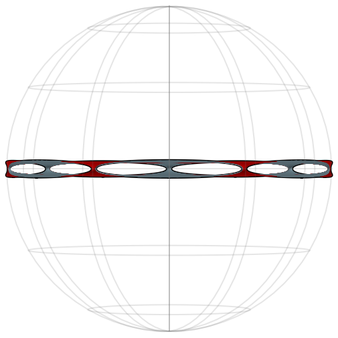} & \includegraphics[width=\panelw]{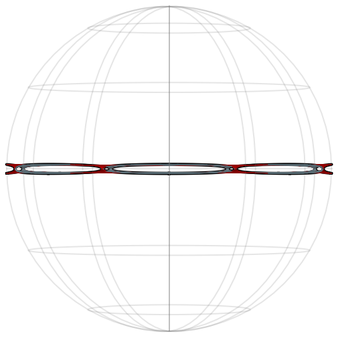} & \includegraphics[width=\panelw]{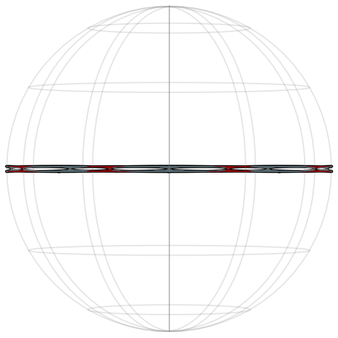} & \includegraphics[width=\panelw]{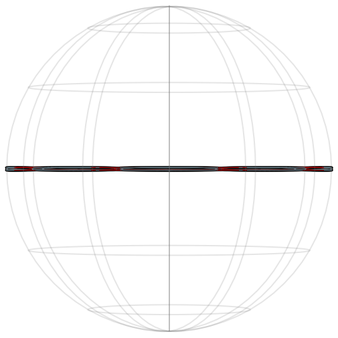}
\end{tabular}

\caption{The equatorial family (Section~\ref{ssec:equatorial}): two sheets joined by $b$ equatorial necks, converging to the doubly covered disk as $b\to\infty$. Values in Table~\ref{tab:catalog}.}\label{fig:type-equatorial}
\end{figure}

\begin{figure}[t]\centering
\setlength{\panelw}{\dimexpr(\textwidth-6\tabcolsep)*100/400\relax}
\setlength{\tabcolsep}{1pt}
\begin{tabular}{@{}cccc@{}}
\paneltitle{$b=5$, $\bar\sigma_1^\star=14.43$\\\zoo{b0005/square_pyramid}{square pyramid}} & \paneltitle{$b=6$, $\bar\sigma_1^\star=14.82$\\\zoo{b0006/pentagonal_pyramid}{pentagonal pyramid}} & \paneltitle{$b=7$, $\bar\sigma_1^\star=15.04$\\\zoo{b0007/hexagonal_pyramid}{hexagonal pyramid}} & \paneltitle{$b=8$, $\bar\sigma_1^\star=15.18$\\\zoo{b0008/heptagonal_pyramid}{heptagonal pyramid}}\\[1pt]
\includegraphics[width=\panelw]{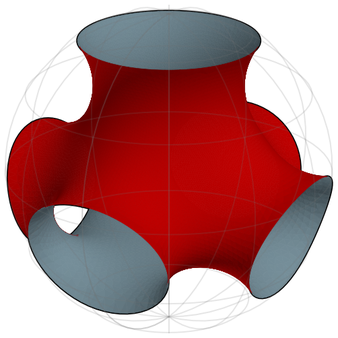} & \includegraphics[width=\panelw]{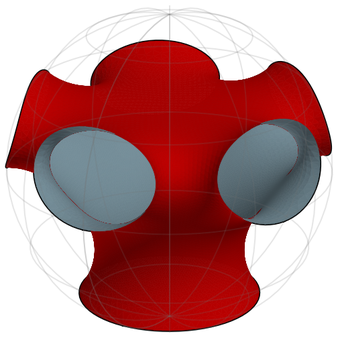} & \includegraphics[width=\panelw]{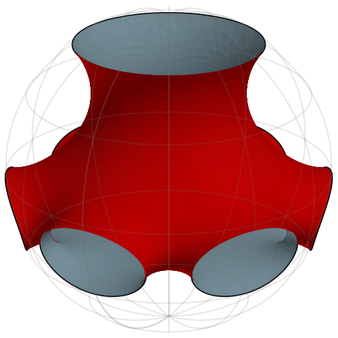} & \includegraphics[width=\panelw]{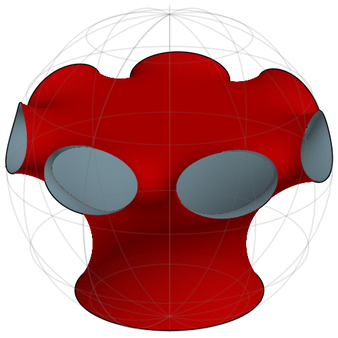}\\[-1pt]
\includegraphics[width=\panelw]{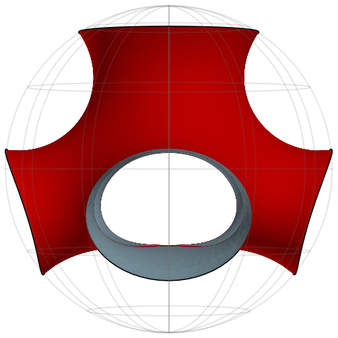} & \includegraphics[width=\panelw]{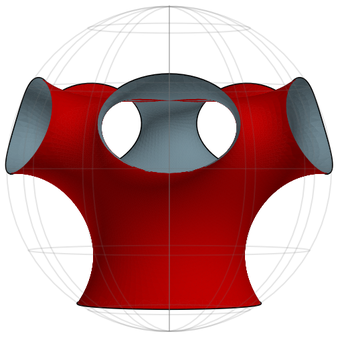} & \includegraphics[width=\panelw]{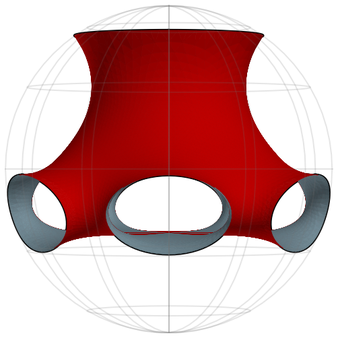} & \includegraphics[width=\panelw]{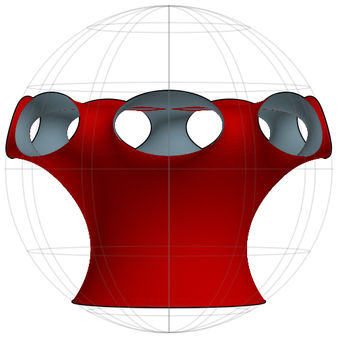}\\[10pt]
\paneltitle{$b=9$, $\bar\sigma_1^\star=15.28$\\\zoo{b0009/octagonal_pyramid}{octagonal pyramid}} & \paneltitle{$b=13$, $\bar\sigma_1^\star=15.48$\\\zoo{b0013/12_gonal_pyramid}{12-gonal pyramid}} & \paneltitle{$b=17$, $\bar\sigma_1^\star=15.57$\\\zoo{b0017/16_gonal_pyramid}{16-gonal pyramid}} & \paneltitle{$b=21$, $\bar\sigma_1^\star=15.61$\\\zoo{b0021/20_gonal_pyramid}{20-gonal pyramid}}\\[1pt]
\includegraphics[width=\panelw]{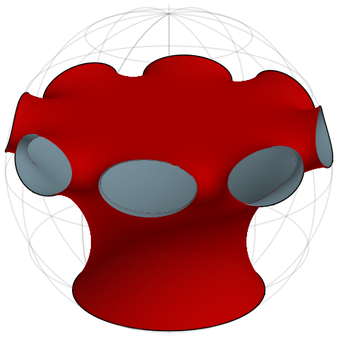} & \includegraphics[width=\panelw]{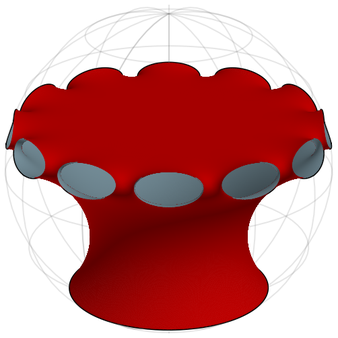} & \includegraphics[width=\panelw]{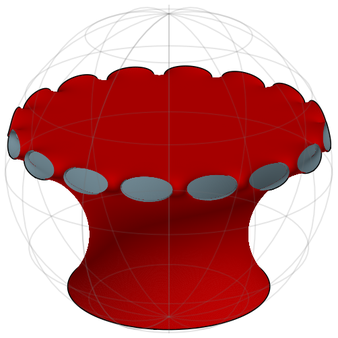} & \includegraphics[width=\panelw]{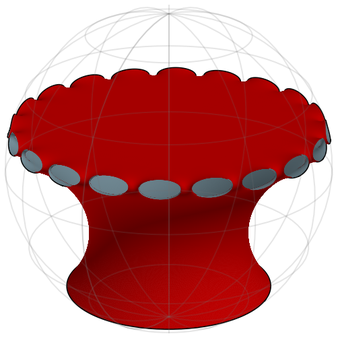}\\[-1pt]
\includegraphics[width=\panelw]{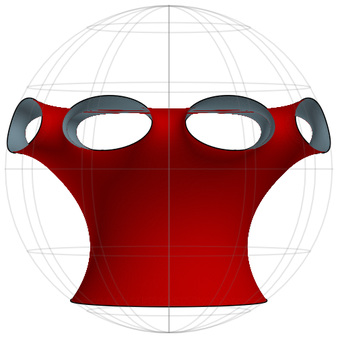} & \includegraphics[width=\panelw]{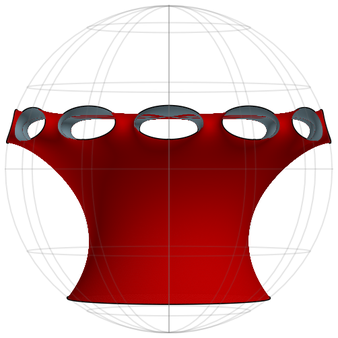} & \includegraphics[width=\panelw]{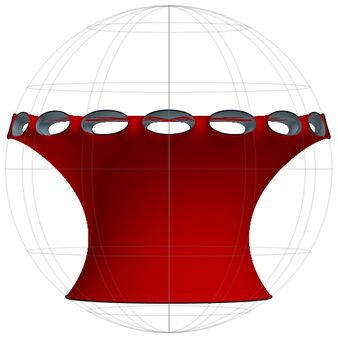} & \includegraphics[width=\panelw]{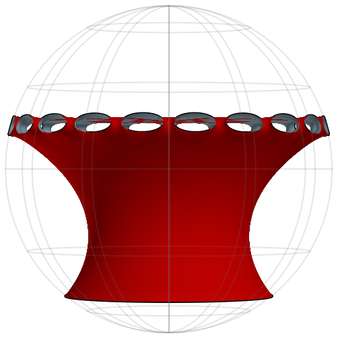}
\end{tabular}

\caption{Pyramids, a ring of $j$ caps and one polar cap, for $j=4,\dots,8$, $12$, $16$, $20$ (Section~\ref{ssec:pyramids}). As $j\to\infty$ they appear to converge to a cake stand, a horn with a round plate (Section~\ref{sec:fcms}). Values in Table~\ref{tab:catalog}.}\label{fig:type-pyramid}
\end{figure}

\begin{figure}[t]\centering
\setlength{\panelw}{\dimexpr(\textwidth-6\tabcolsep)*100/400\relax}
\setlength{\tabcolsep}{1pt}
\begin{tabular}{@{}cccc@{}}
\paneltitle{$b=5$, $\bar\sigma_1^\star=14.47$\\\zoo{b0005/triangular_bipyramid}{triangular bipyramid}} & \paneltitle{$b=7$, $\bar\sigma_1^\star=15.95$\\\zoo{b0007/pentagonal_bipyramid}{pentagonal bipyramid}} & \paneltitle{$b=8$, $\bar\sigma_1^\star=16.28$\\\zoo{b0008/hexagonal_bipyramid}{hexagonal bipyramid}} & \paneltitle{$b=10$, $\bar\sigma_1^\star=16.65$\\\zoo{b0010/octagonal_bipyramid}{octagonal bipyramid}}\\[1pt]
\includegraphics[width=\panelw]{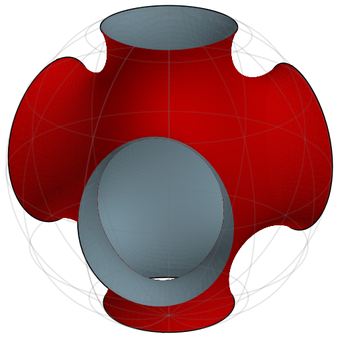} & \includegraphics[width=\panelw]{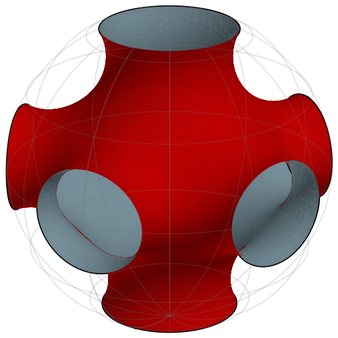} & \includegraphics[width=\panelw]{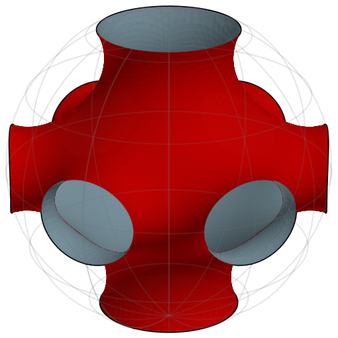} & \includegraphics[width=\panelw]{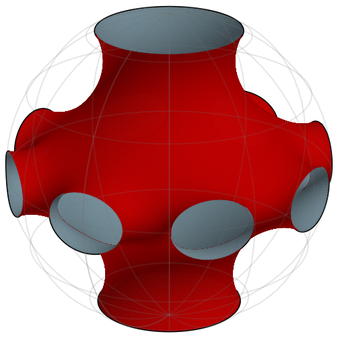}\\[-1pt]
\includegraphics[width=\panelw]{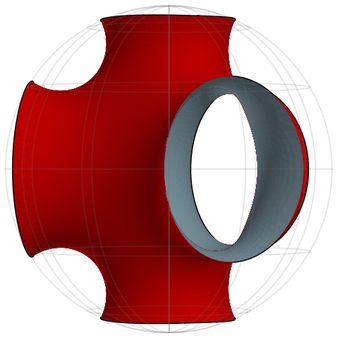} & \includegraphics[width=\panelw]{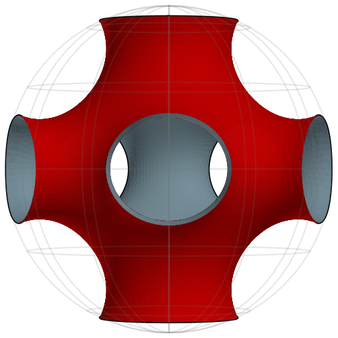} & \includegraphics[width=\panelw]{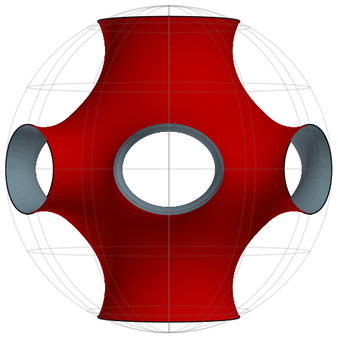} & \includegraphics[width=\panelw]{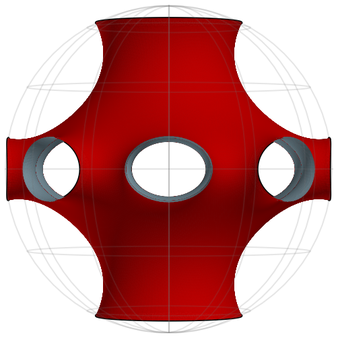}\\[10pt]
\paneltitle{$b=18$, $\bar\sigma_1^\star=17.13$\\\zoo{b0018/16_gonal_bipyramid}{16-gonal bipyramid}} & \paneltitle{$b=20$, $\bar\sigma_1^\star=17.18$\\\zoo{b0020/18_gonal_bipyramid}{18-gonal bipyramid}} & \paneltitle{$b=26$, $\bar\sigma_1^\star=17.28$\\\zoo{b0026/24_gonal_bipyramid}{24-gonal bipyramid}} & \paneltitle{$b=32$, $\bar\sigma_1^\star=17.33$\\\zoo{b0032/30_gonal_bipyramid}{30-gonal bipyramid}}\\[1pt]
\includegraphics[width=\panelw]{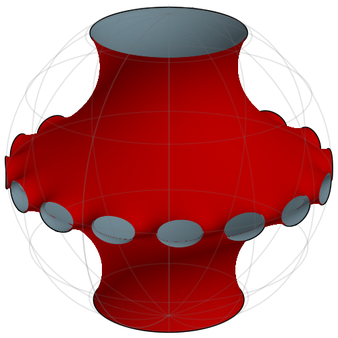} & \includegraphics[width=\panelw]{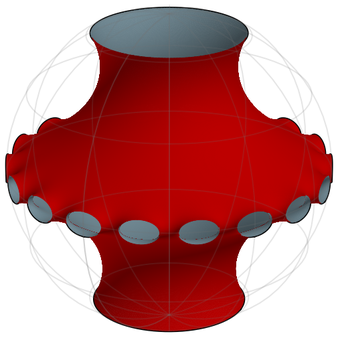} & \includegraphics[width=\panelw]{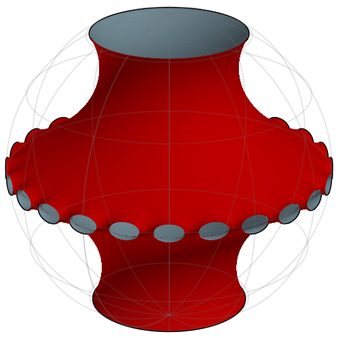} & \includegraphics[width=\panelw]{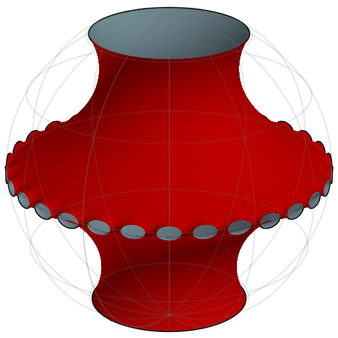}\\[-1pt]
\includegraphics[width=\panelw]{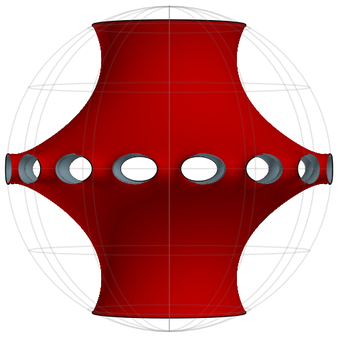} & \includegraphics[width=\panelw]{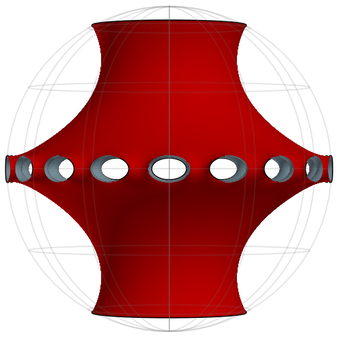} & \includegraphics[width=\panelw]{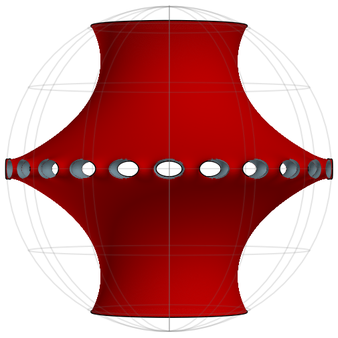} & \includegraphics[width=\panelw]{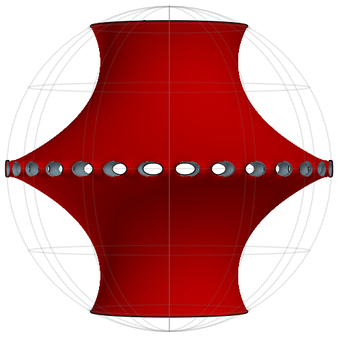}
\end{tabular}

\caption{Bipyramids, an equatorial ring of $j$ caps and two polar caps, for $j=3$, $5$, $6$, $8$, $16$, $18$, $24$, $30$ (Section~\ref{ssec:bipyramids}). As $j\to\infty$ they appear to converge to a wheel, two horns glued along the equator (Section~\ref{sec:fcms}). Values in Table~\ref{tab:catalog}.}\label{fig:type-bipyramid}
\end{figure}

\begin{figure}[t]\centering
\setlength{\panelw}{\dimexpr(\textwidth-6\tabcolsep)*100/400\relax}
\setlength{\tabcolsep}{1pt}
\begin{tabular}{@{}cccc@{}}
\paneltitle{$b=6$, $\bar\sigma_1^\star=15.13$\\\zoo{b0006/triangular_prism}{triangular prism}} & \paneltitle{$b=7$, $\bar\sigma_1^\star=15.45$\\\zoo{b0007/triangular_prism_plus_1pc}{triangular prism + 1PC}} & \paneltitle{$b=8$, $\bar\sigma_1^\star=15.68$\\\zoo{b0008/triangular_prism_plus_2pc}{triangular prism + 2PC}} & \paneltitle{$b=14$, $\bar\sigma_1^\star=17.95$\\\zoo{b0014/heptagonal_prism}{heptagonal prism}}\\[1pt]
\includegraphics[width=\panelw]{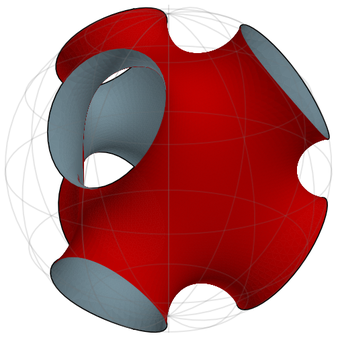} & \includegraphics[width=\panelw]{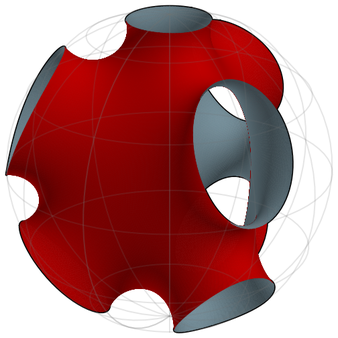} & \includegraphics[width=\panelw]{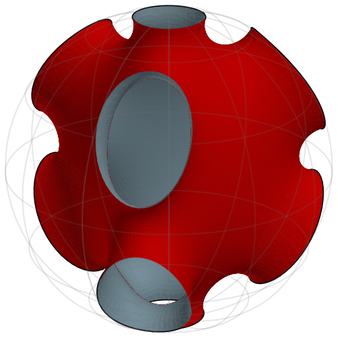} & \includegraphics[width=\panelw]{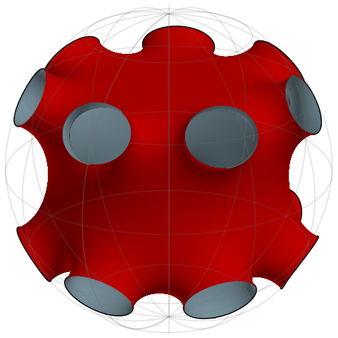}\\[-1pt]
\includegraphics[width=\panelw]{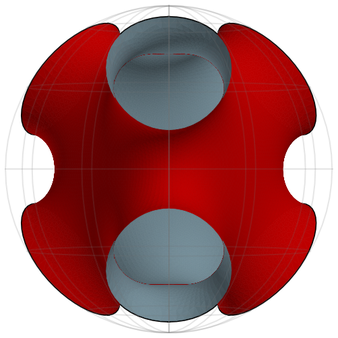} & \includegraphics[width=\panelw]{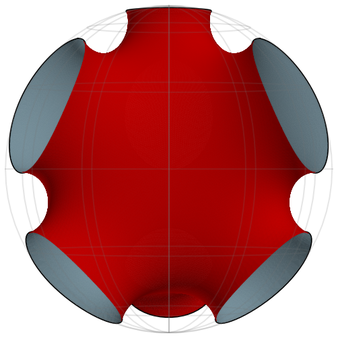} & \includegraphics[width=\panelw]{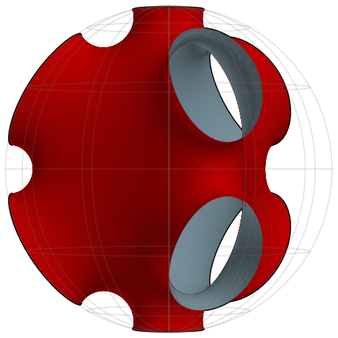} & \includegraphics[width=\panelw]{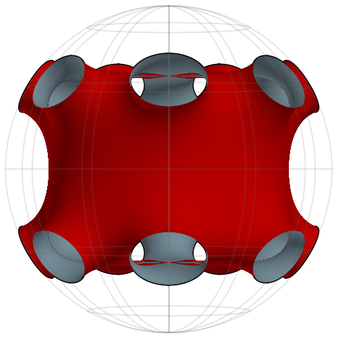}\\[10pt]
\paneltitle{$b=10$, $\bar\sigma_1^\star=17.16$\\\zoo{b0010/pentagonal_prism}{pentagonal prism}} & \paneltitle{$b=18$, $\bar\sigma_1^\star=18.37$\\\zoo{b0018/9_gonal_prism}{9-gonal prism}} & \paneltitle{$b=20$, $\bar\sigma_1^\star=19.42$\\\zoo{b0020/9_gonal_prism_plus_2pc}{9-gonal prism + 2PC}} & \paneltitle{$b=35$, $\bar\sigma_1^\star=19.71$\\\zoo{b0035/17_gonal_antiprism_plus_1pc}{17-gonal prism + 1PC}}\\[1pt]
\includegraphics[width=\panelw]{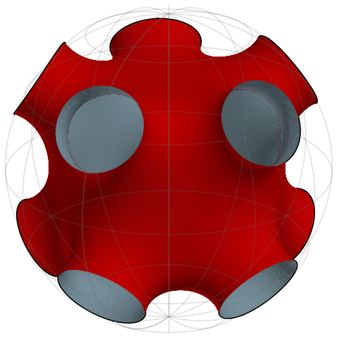} & \includegraphics[width=\panelw]{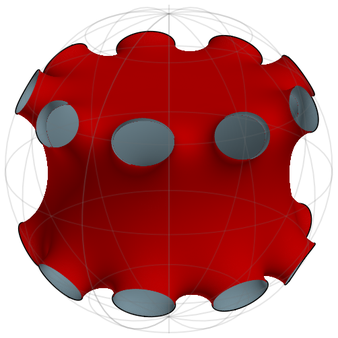} & \includegraphics[width=\panelw]{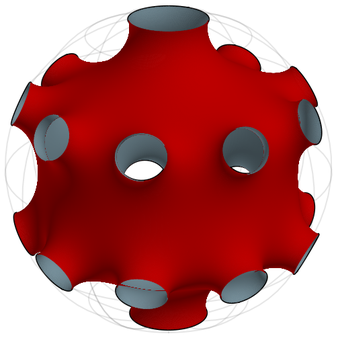} & \includegraphics[width=\panelw]{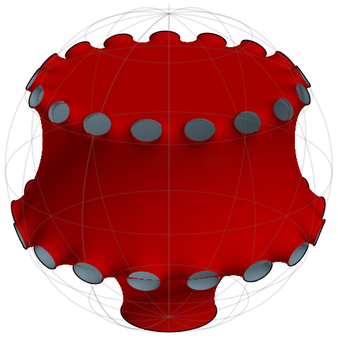}\\[-1pt]
\includegraphics[width=\panelw]{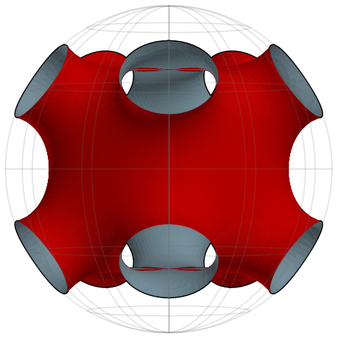} & \includegraphics[width=\panelw]{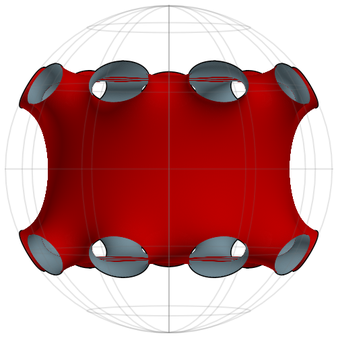} & \includegraphics[width=\panelw]{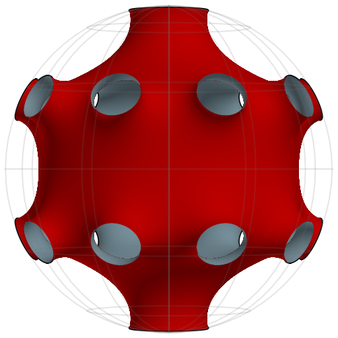} & \includegraphics[width=\panelw]{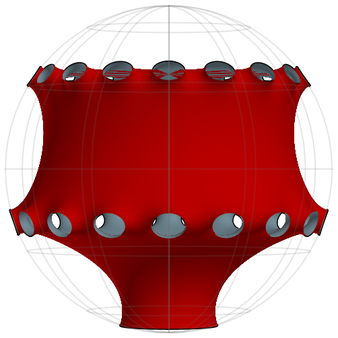}
\end{tabular}

\caption{Prisms with zero, one or two polar caps (PC) (Section~\ref{ssec:prisms}). As $j\to\infty$ they appear to converge to the critical drum, and with one or two polar caps to the drum with one or both disks replaced by horns (Section~\ref{sec:fcms}). Values in Table~\ref{tab:catalog}.}\label{fig:type-prism}
\end{figure}

\begin{figure}[t]\centering
\setlength{\panelw}{\dimexpr(\textwidth-6\tabcolsep)*100/400\relax}
\setlength{\tabcolsep}{1pt}
\begin{tabular}{@{}cccc@{}}
\paneltitle{$b=8$, $\bar\sigma_1^\star=16.50$\\\zoo{b0008/square_antiprism}{square antiprism}} & \paneltitle{$b=10$, $\bar\sigma_1^\star=17.18$\\\zoo{b0010/pentagonal_antiprism}{pentagonal antiprism}} & \paneltitle{$b=7$, $\bar\sigma_1^\star=15.92$\\\zoo{b0007/triangular_antiprism_plus_1pc}{triangular antiprism + 1PC}} & \paneltitle{$b=19$, $\bar\sigma_1^\star=18.97$\\\zoo{b0019/9_gonal_antiprism_plus_1pc}{9-gonal antiprism + 1PC}}\\[1pt]
\includegraphics[width=\panelw]{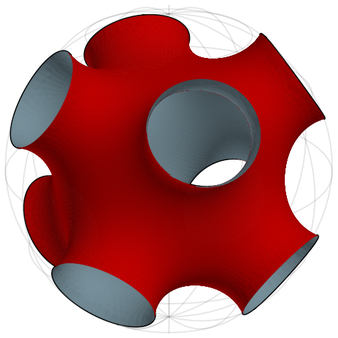} & \includegraphics[width=\panelw]{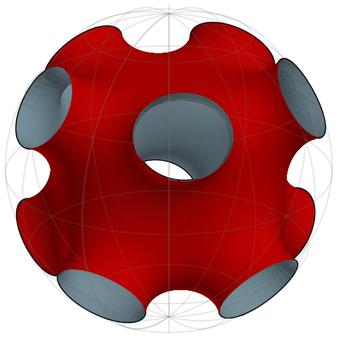} & \includegraphics[width=\panelw]{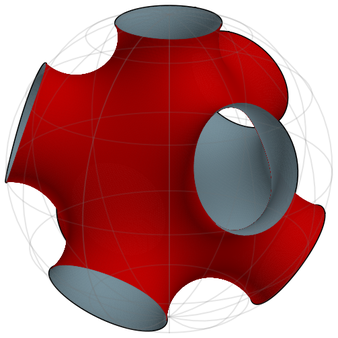} & \includegraphics[width=\panelw]{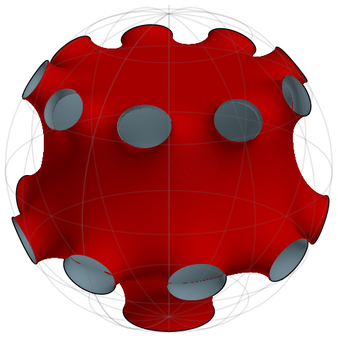}\\[-1pt]
\includegraphics[width=\panelw]{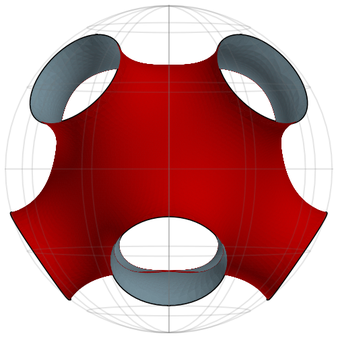} & \includegraphics[width=\panelw]{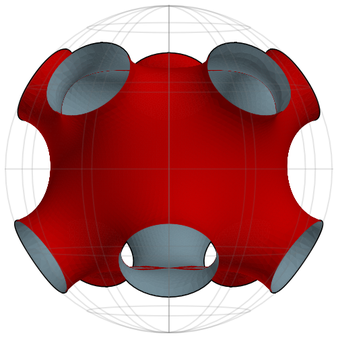} & \includegraphics[width=\panelw]{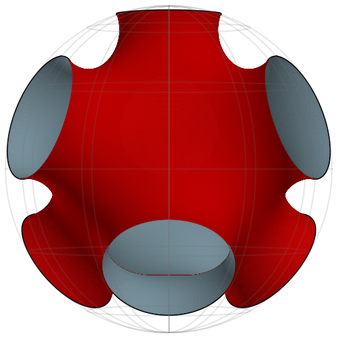} & \includegraphics[width=\panelw]{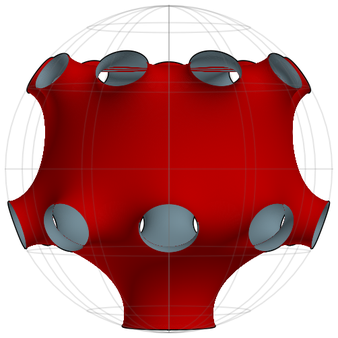}\\[10pt]
\paneltitle{$b=10$, $\bar\sigma_1^\star=17.37$\\\zoo{b0010/square_antiprism_plus_2pc}{square antiprism + 2PC}} & \paneltitle{$b=14$, $\bar\sigma_1^\star=18.57$\\\zoo{b0014/hexagonal_antiprism_plus_2pc}{hexagonal antiprism + 2PC}} & \paneltitle{$b=40$, $\bar\sigma_1^\star=19.11$\\\zoo{b0040/20_gonal_antiprism}{20-gonal antiprism}} & \paneltitle{$b=34$, $\bar\sigma_1^\star=20.15$\\\zoo{b0034/16_gonal_antiprism_plus_2pc}{16-gonal antiprism + 2PC}}\\[1pt]
\includegraphics[width=\panelw]{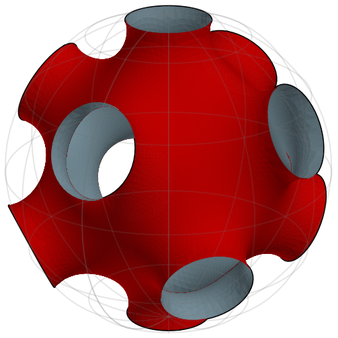} & \includegraphics[width=\panelw]{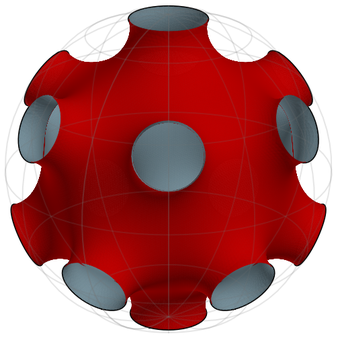} & \includegraphics[width=\panelw]{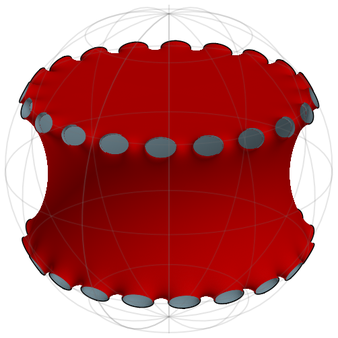} & \includegraphics[width=\panelw]{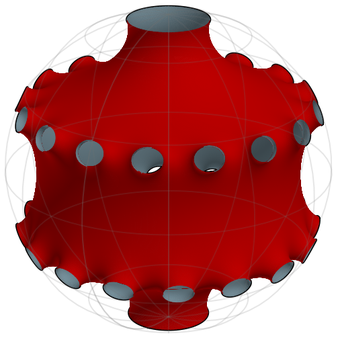}\\[-1pt]
\includegraphics[width=\panelw]{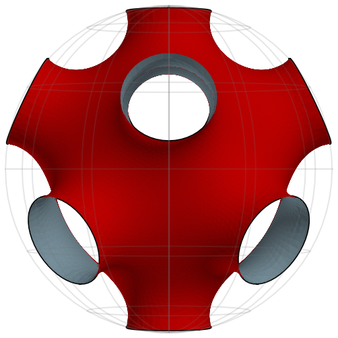} & \includegraphics[width=\panelw]{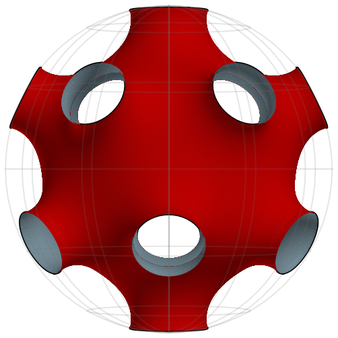} & \includegraphics[width=\panelw]{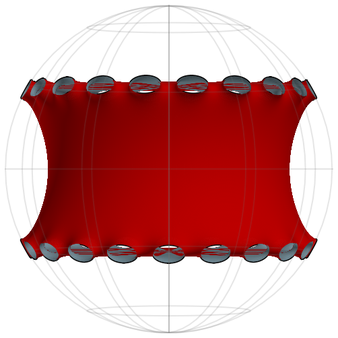} & \includegraphics[width=\panelw]{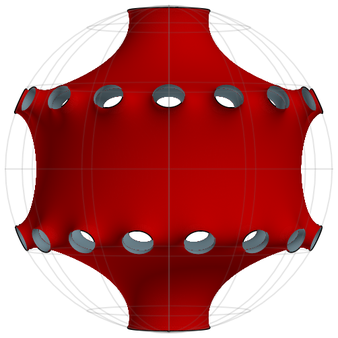}
\end{tabular}

\caption{Antiprisms with zero, one or two polar caps (PC) (Section~\ref{ssec:antiprisms}). As $j\to\infty$ they appear to converge to the critical drum, and with one or two polar caps to the drum with one or both disks replaced by horns (Section~\ref{sec:fcms}). Values in Table~\ref{tab:catalog}.}\label{fig:type-antiprism}
\end{figure}

\begin{figure}[p]\centering
\setlength{\panelw}{\dimexpr(\textwidth-6\tabcolsep)*100/400\relax}
\setlength{\tabcolsep}{1pt}
\begin{tabular}{@{}cccc@{}}
\paneltitle{$b=8$, $\bar\sigma_1^\star=16.49$\\\zoo{b0008/four_stacked_pairs}{four stacked pairs}} & \paneltitle{$b=8$, $\bar\sigma_1^\star=16.29$\\\zoo{b0008/three_pairs_plus_2pc}{three pairs + 2PC}} & \paneltitle{$b=9$, $\bar\sigma_1^\star=16.97$\\\zoo{b0009/three_stacked_triangles}{three stacked triangles}} & \paneltitle{$b=15$, $\bar\sigma_1^\star=18.79$\\\zoo{b0015/five_stacked_triangles}{five stacked triangles}}\\[1pt]
\includegraphics[width=\panelw]{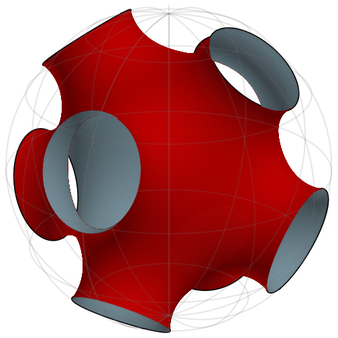} & \includegraphics[width=\panelw]{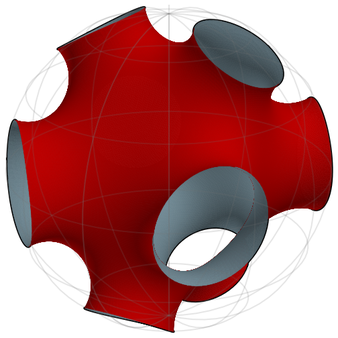} & \includegraphics[width=\panelw]{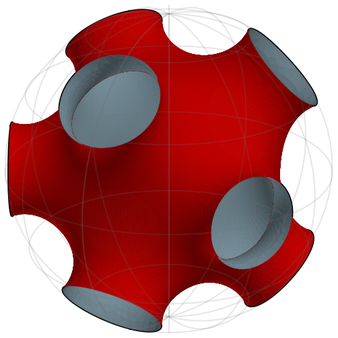} & \includegraphics[width=\panelw]{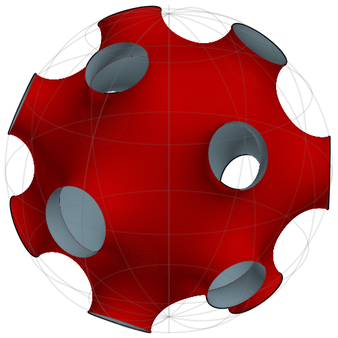}\\[-1pt]
\includegraphics[width=\panelw]{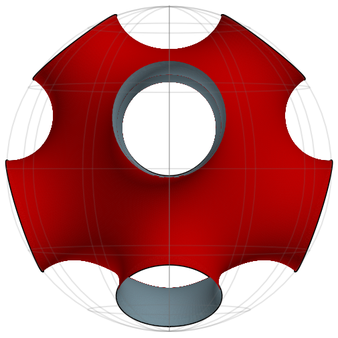} & \includegraphics[width=\panelw]{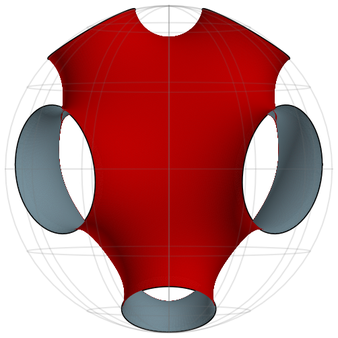} & \includegraphics[width=\panelw]{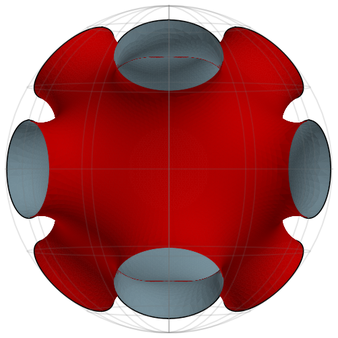} & \includegraphics[width=\panelw]{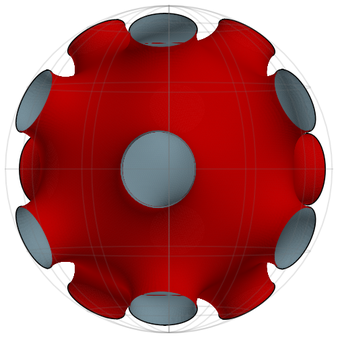}\\[10pt]
\paneltitle{$b=16$, $\bar\sigma_1^\star=19.00$\\\zoo{b0016/four_stacked_squares}{four stacked squares}} & \paneltitle{$b=20$, $\bar\sigma_1^\star=19.51$\\\zoo{b0020/six_stacked_triangles_plus_2pc}{six stacked triangles + 2PC}} & \paneltitle{$b=20$, $\bar\sigma_1^\star=19.70$\\\zoo{b0020/five_stacked_squares}{five stacked squares}} & \paneltitle{$b=13$, $\bar\sigma_1^\star=18.28$\\\zoo{b0013/three_stacked_squares_plus_1pc}{three stacked squares + 1PC}}\\[1pt]
\includegraphics[width=\panelw]{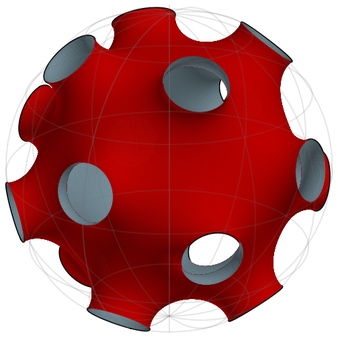} & \includegraphics[width=\panelw]{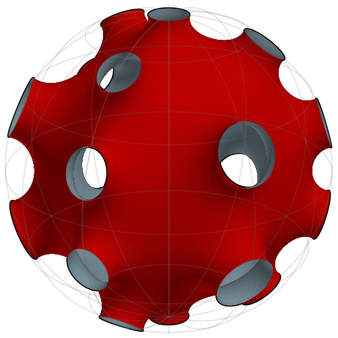} & \includegraphics[width=\panelw]{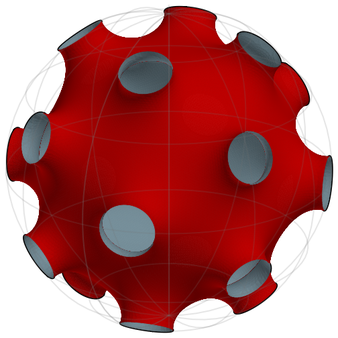} & \includegraphics[width=\panelw]{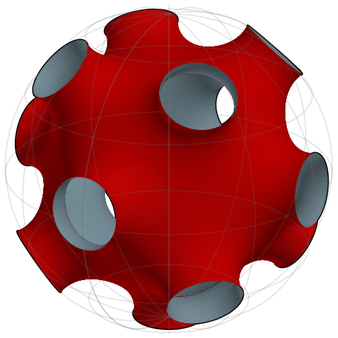}\\[-1pt]
\includegraphics[width=\panelw]{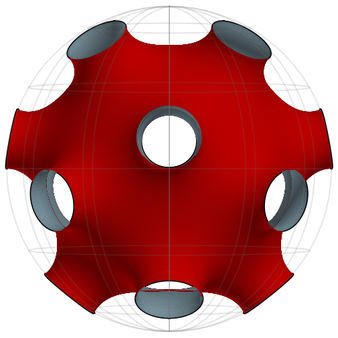} & \includegraphics[width=\panelw]{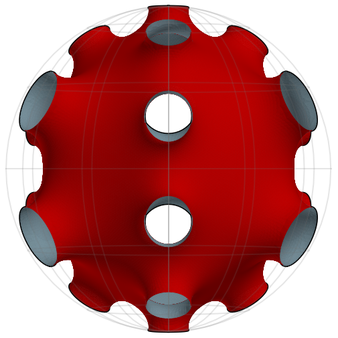} & \includegraphics[width=\panelw]{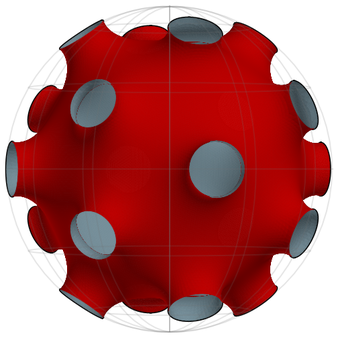} & \includegraphics[width=\panelw]{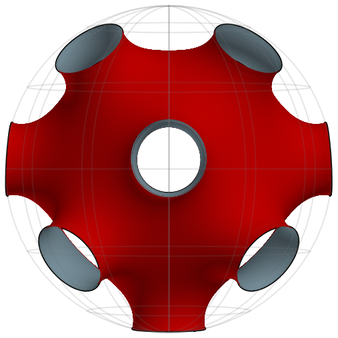}
\end{tabular}

\caption{Stacked latitudinal rings, with or without polar caps (PC) (Section~\ref{ssec:stacked}). As $j\to\infty$, $m$ rings appear to converge to a drum stack, $m-1$ catenoidal annuli closed at each end by a disk, or by a horn where there is a polar cap (Section~\ref{sec:fcms}). Values in Table~\ref{tab:catalog}.}\label{fig:type-stacked}
\end{figure}

\begin{figure}[p]\centering
\setlength{\panelw}{\dimexpr(\textwidth-6\tabcolsep)*100/400\relax}
\setlength{\tabcolsep}{1pt}
\begin{tabular}{@{}cccc@{}}
\paneltitle{$b=2$, $\bar\sigma_1^\star=10.47$\\\zoo{b0002/antipodal_pair}{critical catenoid}} & \paneltitle{$b=3$, $\bar\sigma_1^\star=12.01$\\\zoo{b0003/equatorial_3_gon}{equatorial 3-gon}} & \paneltitle{$b=4$, $\bar\sigma_1^\star=13.67$\\\zoo{b0004/tetrahedron}{tetrahedron}} & \paneltitle{$b=5$, $\bar\sigma_1^\star=14.47$\\\zoo{b0005/triangular_bipyramid}{triangular bipyramid}}\\[1pt]
\includegraphics[width=\panelw]{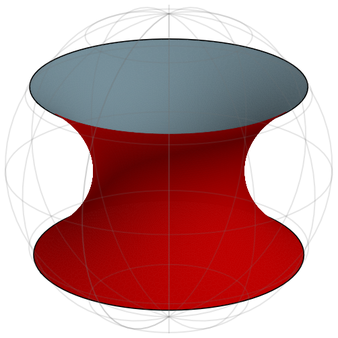} & \includegraphics[width=\panelw]{surfaces_small/b03_triangle_oblique.png} & \includegraphics[width=\panelw]{surfaces_small/b04_tetrahedron_oblique.png} & \includegraphics[width=\panelw]{surfaces_small/b05_triangular_bipyramid_oblique.png}\\[10pt]
\paneltitle{$b=6$, $\bar\sigma_1^\star=15.42$\\\zoo{b0006/octahedron}{octahedron}} & \paneltitle{$b=7$, $\bar\sigma_1^\star=15.95$\\\zoo{b0007/pentagonal_bipyramid}{pentagonal bipyramid}} & \paneltitle{$b=8$, $\bar\sigma_1^\star=16.50$\\\zoo{b0008/square_antiprism}{square antiprism}} & \paneltitle{$b=9$, $\bar\sigma_1^\star=16.97$\\\zoo{b0009/three_stacked_triangles}{three stacked triangles}}\\[1pt]
\includegraphics[width=\panelw]{surfaces_small/b06_octahedron_oblique.png} & \includegraphics[width=\panelw]{surfaces_small/b07_pentagonal_bipyramid_oblique.png} & \includegraphics[width=\panelw]{surfaces_small/b08_square_antiprism_oblique.png} & \includegraphics[width=\panelw]{surfaces_small/b09_triaugmented_prism_oblique.png}\\[10pt]
\paneltitle{$b=10$, $\bar\sigma_1^\star=17.37$\\\zoo{b0010/square_antiprism_plus_2pc}{square antiprism + 2PC}} & \paneltitle{$b=11$, $\bar\sigma_1^\star=17.69$\\\zoo{b0011/c2v_17p6912}{$C_{2v}$: $4+3\times2+1$}} & \paneltitle{$b=12$, $\bar\sigma_1^\star=18.07$\\\zoo{b0012/icosahedron}{icosahedron}} & \paneltitle{$b=13$, $\bar\sigma_1^\star=18.31$\\\zoo{b0013/c2v_18p3064}{$C_{2v}$: $4+4\times2+1$}}\\[1pt]
\includegraphics[width=\panelw]{surfaces_small/b10_square_antiprism_plus_2pc_oblique.png} & \includegraphics[width=\panelw]{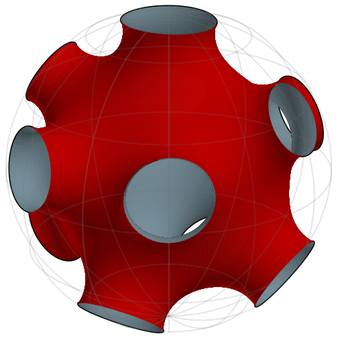} & \includegraphics[width=\panelw]{surfaces_small/b12_icosahedron_oblique.png} & \includegraphics[width=\panelw]{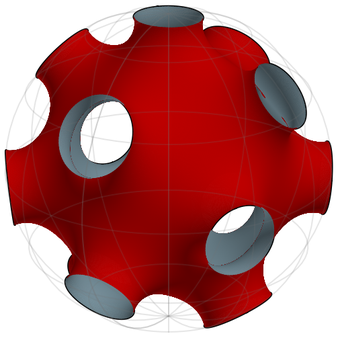}\\[10pt]
\paneltitle{$b=14$, $\bar\sigma_1^\star=18.57$\\\zoo{b0014/hexagonal_antiprism_plus_2pc}{hexagonal antiprism + 2PC}} & \paneltitle{$b=15$, $\bar\sigma_1^\star=18.79$\\\zoo{b0015/five_stacked_triangles}{five stacked triangles}} & \paneltitle{$b=16$, $\bar\sigma_1^\star=19.01$\\\zoo{b0016/td_19p0089}{$T_d$: $12+4$}} & \paneltitle{$b=17$, $\bar\sigma_1^\star=19.20$\\\zoo{b0017/c2v_19p2012}{$C_{2v}$: $2\times4+4\times2+1$}}\\[1pt]
\includegraphics[width=\panelw]{surfaces_small/b14_hexagonal_antiprism_plus_2pc_oblique.png} & \includegraphics[width=\panelw]{surfaces_small/b15_5_triangles_oblique.png} & \includegraphics[width=\panelw]{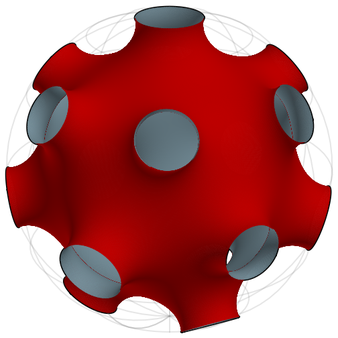} & \includegraphics[width=\panelw]{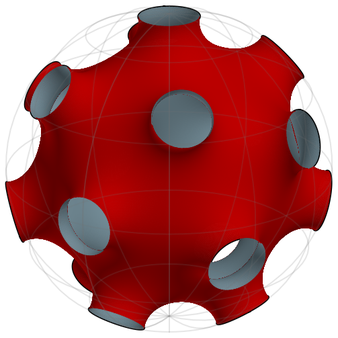}
\end{tabular}
\caption{The best configuration found for each $b=2,\dots,17$ (Section~\ref{ssec:best}, Table~\ref{tab:best}): oblique view of the surface $X^\star(\overline\Omega)\subset\ball$, boundary circles in black. Each name links to the surface in the zoo.}\label{fig:best-a}\label{fig:best-b2}\label{fig:best-b3}\label{fig:best-b4}\label{fig:best-b5}\label{fig:best-b6}\label{fig:best-b7}\label{fig:best-b8}\label{fig:best-b9}\label{fig:best-b10}\label{fig:best-b11}\label{fig:best-b12}\label{fig:best-b13}\label{fig:best-b14}\label{fig:best-b15}\label{fig:best-b16}\label{fig:best-b17}
\end{figure}

\begin{figure}[p]\centering
\setlength{\panelw}{\dimexpr(\textwidth-6\tabcolsep)*100/400\relax}
\setlength{\tabcolsep}{1pt}
\begin{tabular}{@{}cccc@{}}
\paneltitle{$b=18$, $\bar\sigma_1^\star=19.39$\\\zoo{b0018/four_stacked_squares_plus_2pc}{four stacked squares + 2PC}} & \paneltitle{$b=19$, $\bar\sigma_1^\star=19.54$\\\zoo{b0019/fekete_class_19p5446}{Fekete class, $C_{2v}$}} & \paneltitle{$b=20$, $\bar\sigma_1^\star=19.71$\\\zoo{b0020/fekete_class}{Fekete class, $D_{3h}$}} & \paneltitle{$b=21$, $\bar\sigma_1^\star=19.85$\\\zoo{b0021/fekete_class}{Fekete class, $C_{2v}$}}\\[1pt]
\includegraphics[width=\panelw]{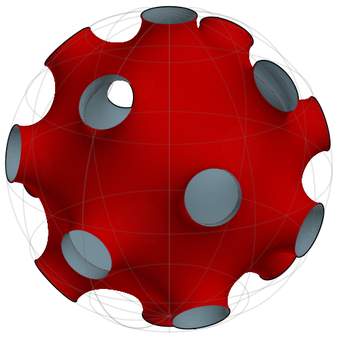} & \includegraphics[width=\panelw]{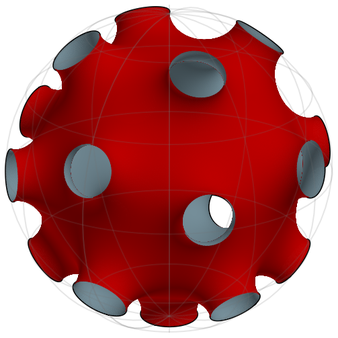} & \includegraphics[width=\panelw]{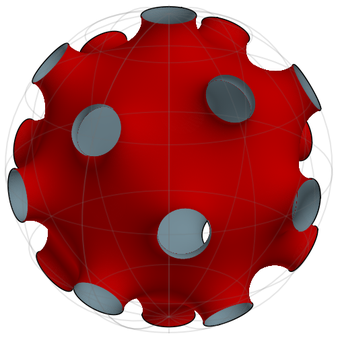} & \includegraphics[width=\panelw]{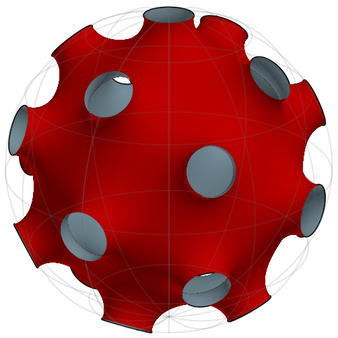}\\[10pt]
\paneltitle{$b=22$, $\bar\sigma_1^\star=19.99$\\\zoo{b0022/fekete_class_19p9935}{Fekete class, $T_d$}} & \paneltitle{$b=23$, $\bar\sigma_1^\star=20.12$\\\zoo{b0023/fekete_class}{Fekete class, $D_{3}$}} & \paneltitle{$b=24$, $\bar\sigma_1^\star=20.24$\\\zoo{b0024/fekete_class}{Fekete class, $O$}} & \paneltitle{$b=25$, $\bar\sigma_1^\star=20.35$\\\zoo{b0025/c1_20p3531}{$C_1$ (20.3531)}}\\[1pt]
\includegraphics[width=\panelw]{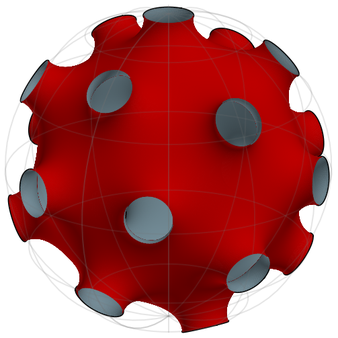} & \includegraphics[width=\panelw]{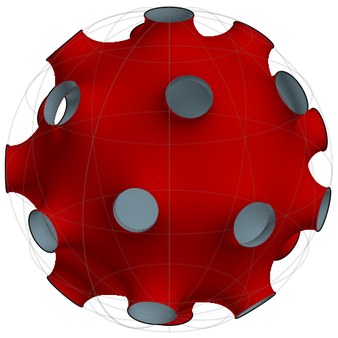} & \includegraphics[width=\panelw]{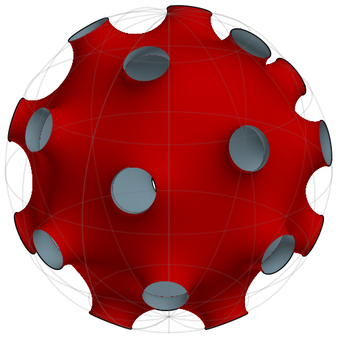} & \includegraphics[width=\panelw]{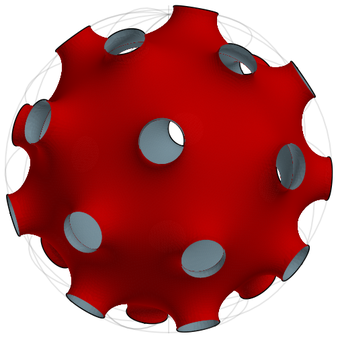}\\[10pt]
\paneltitle{$b=26$, $\bar\sigma_1^\star=20.46$\\\zoo{b0026/fekete_class}{Fekete class, $C_{2}$}} & \paneltitle{$b=27$, $\bar\sigma_1^\star=20.57$\\\zoo{b0027/fekete_class}{Fekete class, $D_{5h}$}} & \paneltitle{$b=28$, $\bar\sigma_1^\star=20.66$\\\zoo{b0028/fekete_class}{Fekete class, $T$}} & \paneltitle{$b=29$, $\bar\sigma_1^\star=20.75$\\\zoo{b0029/fekete_class}{Fekete class, $D_{3}$}}\\[1pt]
\includegraphics[width=\panelw]{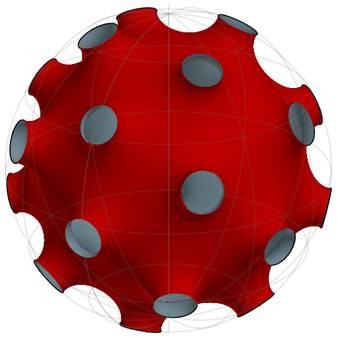} & \includegraphics[width=\panelw]{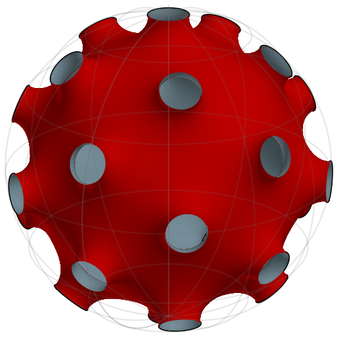} & \includegraphics[width=\panelw]{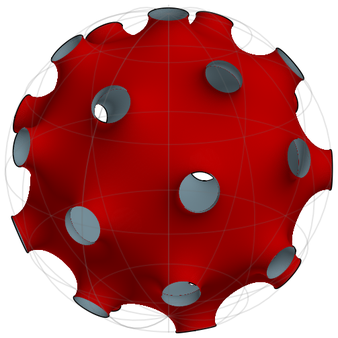} & \includegraphics[width=\panelw]{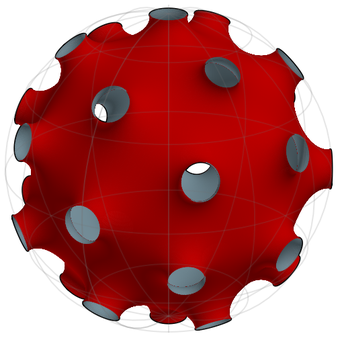}\\[10pt]
\paneltitle{$b=30$, $\bar\sigma_1^\star=20.84$\\\zoo{b0030/fekete_class}{Fekete class, $D_{2}$}} & \paneltitle{$b=31$, $\bar\sigma_1^\star=20.93$\\\zoo{b0031/fekete_class_20p9278}{Fekete class, $C_{3v}$}} & \paneltitle{$b=32$, $\bar\sigma_1^\star=21.01$\\\zoo{b0032/icosahedron_plus_dodecahedron}{icosahedron + dodecahedron}} & \paneltitle{$b=33$, $\bar\sigma_1^\star=21.08$\\\zoo{b0033/c1_21p0835}{$C_1$ (21.0835)}}\\[1pt]
\includegraphics[width=\panelw]{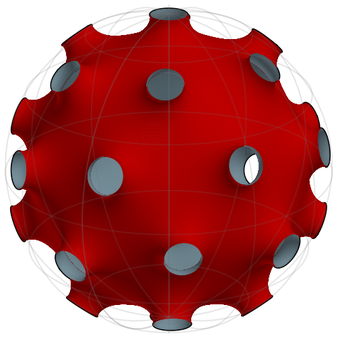} & \includegraphics[width=\panelw]{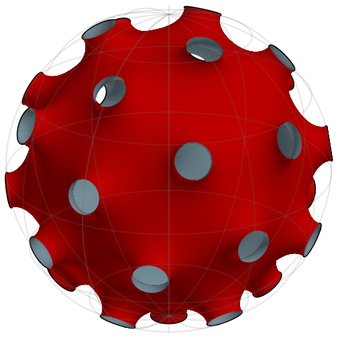} & \includegraphics[width=\panelw]{surfaces_small/b32_icosahedronplusdodecahedron_oblique.png} & \includegraphics[width=\panelw]{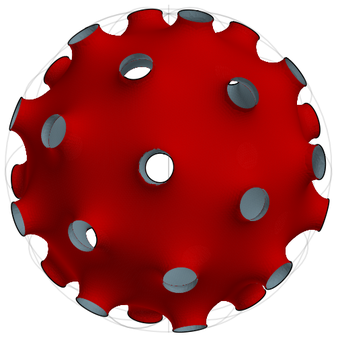}
\end{tabular}
\caption{The best configuration found for each $b=18,\dots,33$ (Section~\ref{ssec:best}, Table~\ref{tab:best}): oblique view of the surface $X^\star(\overline\Omega)\subset\ball$, boundary circles in black. Each name links to the surface in the zoo.}\label{fig:best-b}\label{fig:best-b18}\label{fig:best-b19}\label{fig:best-b20}\label{fig:best-b21}\label{fig:best-b22}\label{fig:best-b23}\label{fig:best-b24}\label{fig:best-b25}\label{fig:best-b26}\label{fig:best-b27}\label{fig:best-b28}\label{fig:best-b29}\label{fig:best-b30}\label{fig:best-b31}\label{fig:best-b32}\label{fig:best-b33}
\end{figure}

\begin{figure}[p]\centering
\setlength{\panelw}{\dimexpr(\textwidth-6\tabcolsep)*100/400\relax}
\setlength{\tabcolsep}{1pt}
\begin{tabular}{@{}cccc@{}}
\paneltitle{$b=34$, $\bar\sigma_1^\star=21.16$\\\zoo{b0034/fekete_class}{Fekete class, $D_{2}$}} & \paneltitle{$b=35$, $\bar\sigma_1^\star=21.23$\\\zoo{b0035/fekete_class}{Fekete class, $C_{2}$}} & \paneltitle{$b=36$, $\bar\sigma_1^\star=21.29$\\\zoo{b0036/fekete_class}{Fekete class, $D_{2}$}} & \paneltitle{$b=37$, $\bar\sigma_1^\star=21.36$\\\zoo{b0037/fekete_class_21p3600}{Fekete class, $D_{5h}$}}\\[1pt]
\includegraphics[width=\panelw]{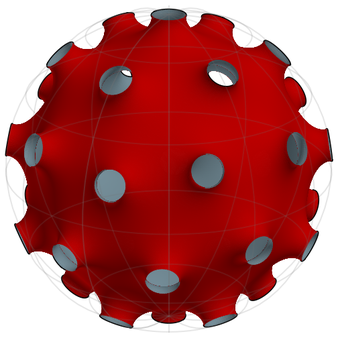} & \includegraphics[width=\panelw]{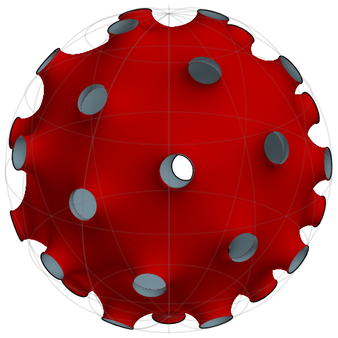} & \includegraphics[width=\panelw]{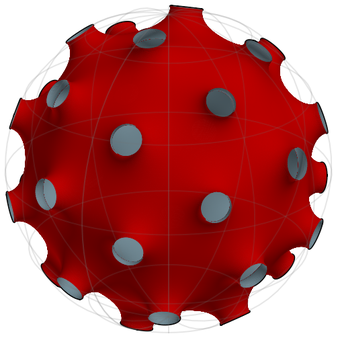} & \includegraphics[width=\panelw]{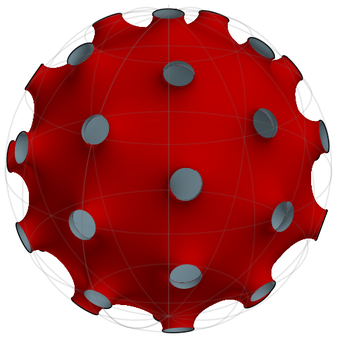}\\[10pt]
\paneltitle{$b=38$, $\bar\sigma_1^\star=21.42$\\\zoo{b0038/six_stacked_hexagons_plus_2pc_21p4221}{six stacked hexagons + 2PC}} & \paneltitle{$b=39$, $\bar\sigma_1^\star=21.48$\\\zoo{b0039/fekete_class_21p4831}{Fekete class, $D_{3h}$}} & \paneltitle{$b=40$, $\bar\sigma_1^\star=21.54$\\\zoo{b0040/fekete_class_21p5411}{Fekete class, $T_d$}} & \paneltitle{$b=42$, $\bar\sigma_1^\star=21.65$\\\zoo{b0042/geodesic_icosahedron_nu2}{geodesic icosahedron $\nu=2$}}\\[1pt]
\includegraphics[width=\panelw]{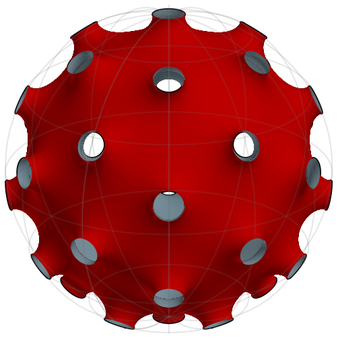} & \includegraphics[width=\panelw]{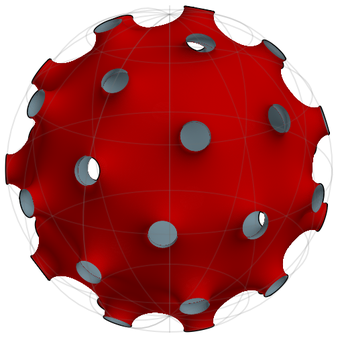} & \includegraphics[width=\panelw]{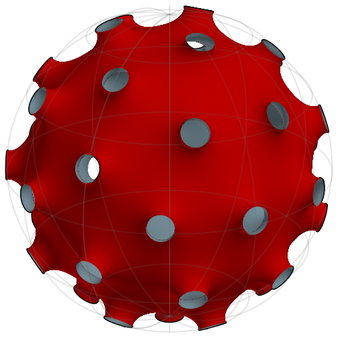} & \includegraphics[width=\panelw]{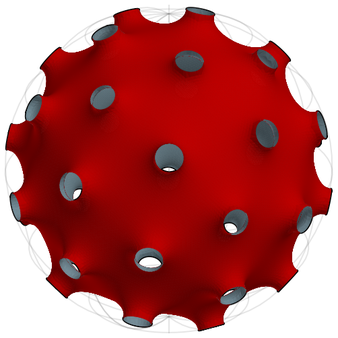}\\[10pt]
\paneltitle{$b=60$, $\bar\sigma_1^\star=22.28$\\\zoo{b0060/60_orbit_on_the_vertex_to_vertex_arcs}{60-orbit on the vertex-to-vertex arcs}} & \paneltitle{$b=72$, $\bar\sigma_1^\star=22.70$\\\zoo{b0072/icosahedron_plus_generic_60_orbit}{icosahedron + generic 60-orbit}} & \paneltitle{$b=92$, $\bar\sigma_1^\star=23.08$\\\zoo{b0092/geodesic_icosahedron_nu3}{geodesic icosahedron $\nu=3$}} & \paneltitle{$b=162$, $\bar\sigma_1^\star=23.77$\\\zoo{b0162/geodesic_icosahedron_nu4}{geodesic icosahedron $\nu=4$}}\\[1pt]
\includegraphics[width=\panelw]{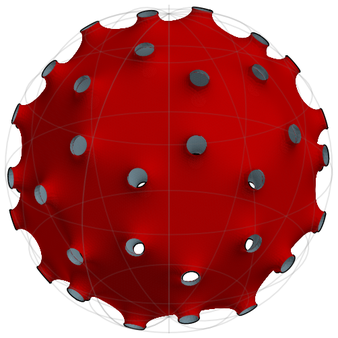} & \includegraphics[width=\panelw]{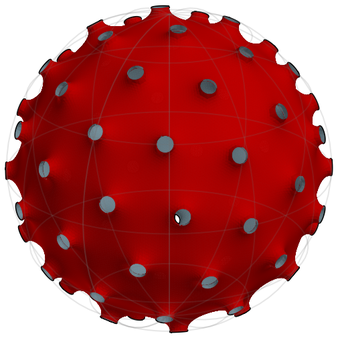} & \includegraphics[width=\panelw]{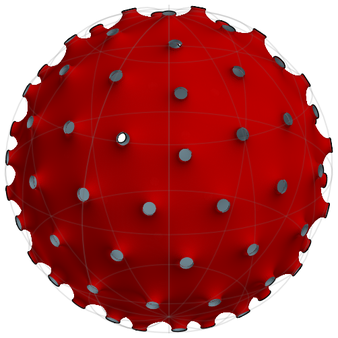} & \includegraphics[width=\panelw]{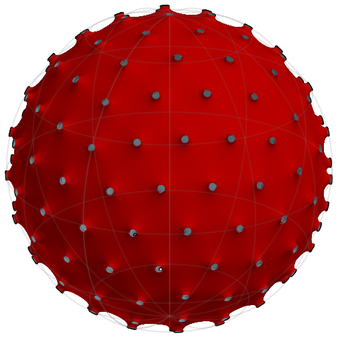}\\[10pt]
\paneltitle{$b=252$, $\bar\sigma_1^\star=24.16$\\\zoo{b0252/geodesic_icosahedron_nu5}{geodesic icosahedron $\nu=5$}} & \paneltitle{$b=362$, $\bar\sigma_1^\star=24.40$\\\zoo{b0362/geodesic_icosahedron_nu6}{geodesic icosahedron $\nu=6$}} & \paneltitle{$b=492$, $\bar\sigma_1^\star=24.56$\\\zoo{b0492/geodesic_icosahedron_nu7}{geodesic icosahedron $\nu=7$}} & \paneltitle{$b=642$, $\bar\sigma_1^\star=24.67$\\\zoo{b0642/geodesic_icosahedron_nu8}{geodesic icosahedron $\nu=8$}}\\[1pt]
\includegraphics[width=\panelw]{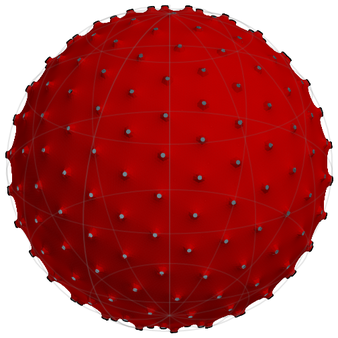} & \includegraphics[width=\panelw]{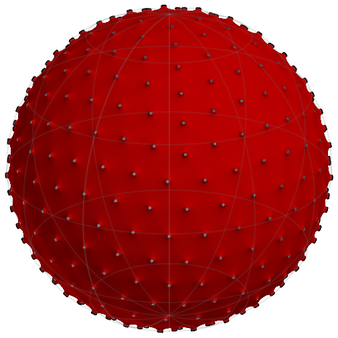} & \includegraphics[width=\panelw]{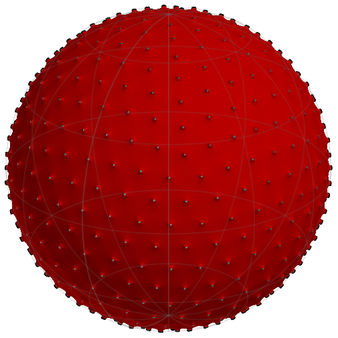} & \includegraphics[width=\panelw]{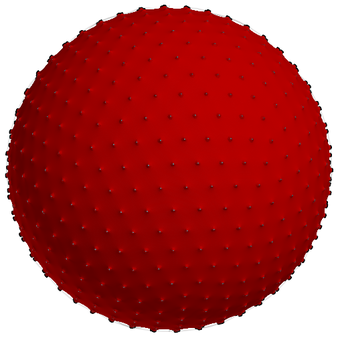}
\end{tabular}
\caption{The best configuration found for each $b$ computed with $34\le b\le642$ (Section~\ref{ssec:best}, Table~\ref{tab:best}): oblique view of the surface $X^\star(\overline\Omega)\subset\ball$, boundary circles in black. Each name links to the surface in the zoo.}\label{fig:best-c}\label{fig:best-b34}\label{fig:best-b35}\label{fig:best-b36}\label{fig:best-b37}\label{fig:best-b38}\label{fig:best-b39}\label{fig:best-b40}\label{fig:best-b42}\label{fig:best-b60}\label{fig:best-b72}\label{fig:best-b92}\label{fig:best-b162}\label{fig:best-b252}\label{fig:best-b362}\label{fig:best-b492}\label{fig:best-b642}
\end{figure}

\begin{figure}[p]\centering
\setlength{\panelw}{\dimexpr(\textwidth-6\tabcolsep)*100/400\relax}
\setlength{\tabcolsep}{1pt}
\begin{tabular}{@{}cccc@{}}
\paneltitle{$b=812$, $\bar\sigma_1^\star=24.75$\\\zoo{b0812/geodesic_icosahedron_nu9}{geodesic icosahedron $\nu=9$}} & \paneltitle{$b=1002$, $\bar\sigma_1^\star=24.81$\\\zoo{b1002/geodesic_icosahedron_nu10}{geodesic icosahedron $\nu=10$}} & \paneltitle{$b=1212$, $\bar\sigma_1^\star=24.86$\\\zoo{b1212/geodesic_icosahedron_nu11}{geodesic icosahedron $\nu=11$}} & \paneltitle{$b=1442$, $\bar\sigma_1^\star=24.90$\\\zoo{b1442/geodesic_icosahedron_nu12}{geodesic icosahedron $\nu=12$}}\\[1pt]
\includegraphics[width=\panelw]{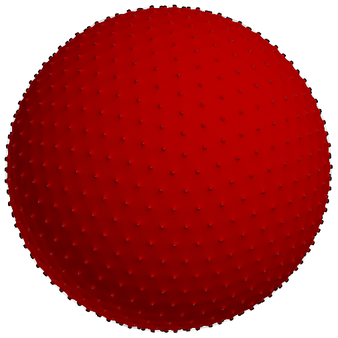} & \includegraphics[width=\panelw]{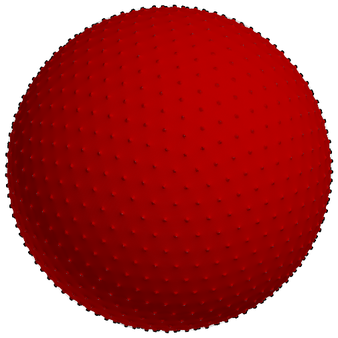} & \includegraphics[width=\panelw]{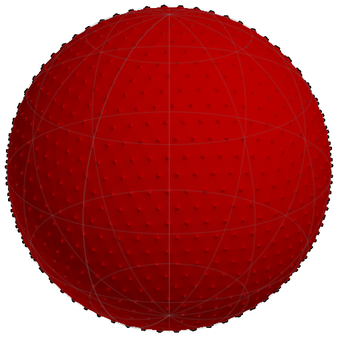} & \includegraphics[width=\panelw]{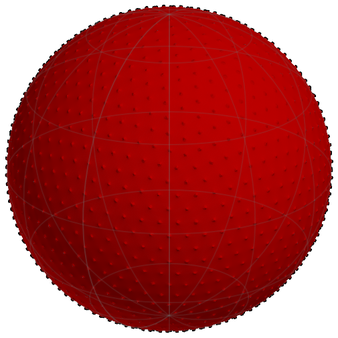}\\[10pt]
\paneltitle{$b=1692$, $\bar\sigma_1^\star=24.93$\\\zoo{b1692/geodesic_icosahedron_nu13}{geodesic icosahedron $\nu=13$}} & \paneltitle{$b=1962$, $\bar\sigma_1^\star=24.95$\\\zoo{b1962/geodesic_icosahedron_nu14}{geodesic icosahedron $\nu=14$}} & \paneltitle{$b=2252$, $\bar\sigma_1^\star=24.97$\\\zoo{b2252/geodesic_icosahedron_nu15}{geodesic icosahedron $\nu=15$}} & \\[1pt]
\includegraphics[width=\panelw]{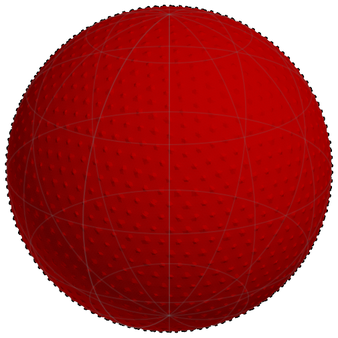} & \includegraphics[width=\panelw]{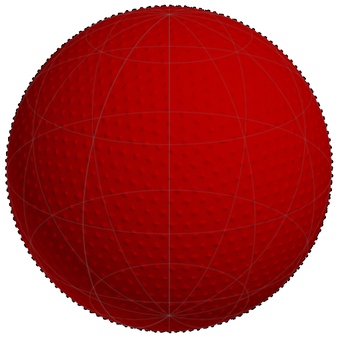} & \includegraphics[width=\panelw]{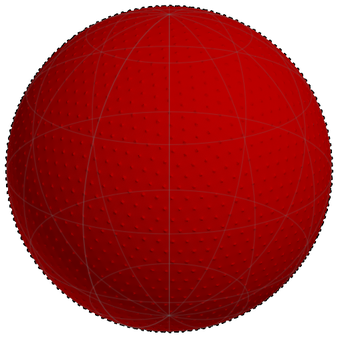} & 
\end{tabular}
\caption{The best configuration found for each $b$ computed with $812\le b\le2252$ (Section~\ref{ssec:best}, Table~\ref{tab:best}): oblique view of the surface $X^\star(\overline\Omega)\subset\ball$, boundary circles in black. Each name links to the surface in the zoo. A larger view of the surface with $b=2252$ is given in Figure~\ref{fig:headline}.}\label{fig:best-d}\label{fig:best-b812}\label{fig:best-b1002}\label{fig:best-b1212}\label{fig:best-b1442}\label{fig:best-b1692}\label{fig:best-b1962}\label{fig:best-b2252}\label{fig:best-last}
\end{figure}

\clearpage
\clearpage
\twocolumn[{\captionof{table}{Catalog of the zoo: the certified critical points by family (Sections~\ref{ssec:platonic}--\ref{ssec:thomson}), and within a family by $b$ and decreasing value. The columns give the symmetry group of the result, $\bar\sigma_1^\star$ (italics: best for its $b$), the Morse index (Section~\ref{ssec:catalog}; a dash where not computed), the order of magnitude of the residual $\mathcal K$ of \eqref{eq:K}, as $-\log_{10}\mathcal K$ rounded down, so that an entry $n$ means $\mathcal K\le10^{-n}$, and the figure of the surface. Each surface is listed once, under the most restricted class that reaches it. In the class names, $\nu$ is the frequency of the geodesic icosahedron, and `3 triangles' means three stacked triangles. Each class name links to its data folder in the repository \cite{Chen_Osting_genus0_code}.}\label{tab:catalog}\vspace{4pt}}]
\begingroup\scriptsize\setlength{\tabcolsep}{2pt}
\tablehead{\toprule $b$ & class & sym & $\bar\sigma_1^\star$ & ind & $-\log_{10}\!\mathcal K$ & fig.\\\midrule}
\tabletail{\bottomrule}
\begin{supertabular}{@{}rC{2.45cm}lrrrr@{}}
\multicolumn{6}{@{}l}{\textbf{Polyhedral configurations}}\\[2pt]
4 & \zoo{b0004/tetrahedron}{tetrahedron} & $T_d$ & \textit{13.66761484} & 0 & 13 & \ref{fig:best-a}\\
6 & \zoo{b0006/octahedron}{octahedron} & $O_h$ & \textit{15.42492045} & 0 & 9 & \ref{fig:best-a}\\
8 & \zoo{b0008/cube}{cube} & $O_h$ & 16.42485423 & 2 & 12 & \ref{fig:type-polyhedral}\\
12 & \zoo{b0012/icosahedron}{icosahedron} & $I_h$ & \textit{18.06871939} & 0 & 12 & \ref{fig:best-a}\\
20 & \zoo{b0020/dodecahedron}{dodecahedron} & $I_h$ & 19.57267499 & 8 & 12 & \ref{fig:type-polyhedral}\\
30 & \zoo{b0030/icosidodecahedron}{icosidodecahedron} & $I_h$ & 20.73171822 & 11 & 12 & \ref{fig:type-polyhedral}\\
32 & \zoo{b0032/icosahedron_plus_dodecahedron}{icosahedron + dodecahedron} & $I_h$ & \textit{21.01400764} & 0 & 12 & \ref{fig:best-b}\\
60 & \zoo{b0060/60_orbit_on_the_vertex_to_vertex_arcs}{60-orbit on edge arcs} & $I_h$ & \textit{22.28356339} & 28 & 12 & \ref{fig:best-c}\\
72 & \zoo{b0072/icosahedron_plus_generic_60_orbit}{icosahedron + generic 60-orbit} & $I$ & \textit{22.69609123} & 0 & 11 & \ref{fig:best-c}\\
72 & \zoo{b0072/icosahedron_plus_60_orbit_on_the_mirror_arcs}{icosahedron + 60-orbit on mirror arcs} & $I_h$ & 22.57754557 & -- & 9 & --\\
\addlinespace[4pt]
\multicolumn{6}{@{}l}{\textbf{Geodesic icosahedra}}\\[2pt]
42 & \zoo{b0042/geodesic_icosahedron_nu2}{$\nu=2$} & $I_h$ & \textit{21.64565725} & 0 & 11 & \ref{fig:best-c}\\
92 & \zoo{b0092/geodesic_icosahedron_nu3}{$\nu=3$} & $I_h$ & \textit{23.07686419} & 0 & 10 & \ref{fig:best-c}\\
92 & \zoo{b0092/icosahedral_nu3}{$\nu=3$, saddle} & $I_h$ & 22.99638099 & 50 & 11 & --\\
162 & \zoo{b0162/geodesic_icosahedron_nu4}{$\nu=4$} & $I_h$ & \textit{23.77221789} & 0 & 9 & \ref{fig:best-c}\\
162 & \zoo{b0162/icosahedral_nu4}{$\nu=4$, saddle} & $I_h$ & 23.75163737 & 52 & 10 & --\\
252 & \zoo{b0252/geodesic_icosahedron_nu5}{$\nu=5$} & $I_h$ & \textit{24.16173249} & 0 & 8 & \ref{fig:best-c}\\
362 & \zoo{b0362/geodesic_icosahedron_nu6}{$\nu=6$} & $I_h$ & \textit{24.40228582} & 0 & 9 & \ref{fig:best-c}\\
492 & \zoo{b0492/geodesic_icosahedron_nu7}{$\nu=7$} & $I_h$ & \textit{24.56163111} & 0 & 8 & \ref{fig:best-c}\\
642 & \zoo{b0642/geodesic_icosahedron_nu8}{$\nu=8$} & $I_h$ & \textit{24.67287741} & 0 & 8 & \ref{fig:best-c}\\
812 & \zoo{b0812/geodesic_icosahedron_nu9}{$\nu=9$} & $I_h$ & \textit{24.75377466} & -- & 9 & \ref{fig:best-d}\\
1002 & \zoo{b1002/geodesic_icosahedron_nu10}{$\nu=10$} & $I_h$ & \textit{24.81454590} & -- & 9 & \ref{fig:best-d}\\
1212 & \zoo{b1212/geodesic_icosahedron_nu11}{$\nu=11$} & $I_h$ & \textit{24.86142450} & -- & 9 & \ref{fig:best-d}\\
1442 & \zoo{b1442/geodesic_icosahedron_nu12}{$\nu=12$} & $I_h$ & \textit{24.89839157} & -- & 9 & \ref{fig:best-d}\\
1692 & \zoo{b1692/geodesic_icosahedron_nu13}{$\nu=13$} & $I_h$ & \textit{24.92808940} & -- & 9 & \ref{fig:best-d}\\
1962 & \zoo{b1962/geodesic_icosahedron_nu14}{$\nu=14$} & $I_h$ & \textit{24.95232944} & -- & 9 & \ref{fig:best-d}\\
2252 & \zoo{b2252/geodesic_icosahedron_nu15}{$\nu=15$} & $I_h$ & \textit{24.97238832} & -- & 8 & \ref{fig:best-d}\\
\addlinespace[4pt]
\multicolumn{6}{@{}l}{\textbf{The equatorial family}}\\[2pt]
3 & \zoo{b0003/equatorial_3_gon}{equatorial 3-gon} & $D_{3h}$ & \textit{12.01198922} & 0 & 9 & \ref{fig:best-a}\\
4 & \zoo{b0004/equatorial_4_gon}{equatorial 4-gon} & $D_{4h}$ & 12.42990413 & 1 & 13 & \ref{fig:type-equatorial}\\
5 & \zoo{b0005/equatorial_5_gon}{equatorial 5-gon} & $D_{5h}$ & 12.53139727 & 2 & 12 & \ref{fig:type-equatorial}\\
6 & \zoo{b0006/equatorial_6_gon}{equatorial 6-gon} & $D_{6h}$ & 12.55681968 & 3 & 12 & \ref{fig:type-equatorial}\\
7 & \zoo{b0007/equatorial_7_gon}{equatorial 7-gon} & $D_{7h}$ & 12.56361680 & 4 & 12 & \ref{fig:type-equatorial}\\
8 & \zoo{b0008/equatorial_8_gon}{equatorial 8-gon} & $D_{8h}$ & 12.56554312 & 5 & 11 & \ref{fig:type-equatorial}\\
9 & \zoo{b0009/equatorial_9_gon}{equatorial 9-gon} & $D_{9h}$ & 12.56611415 & 6 & 12 & \ref{fig:type-equatorial}\\
10 & \zoo{b0010/equatorial_10_gon}{equatorial 10-gon} & $D_{10h}$ & 12.56628924 & 7 & 10 & \ref{fig:type-equatorial}\\
11 & \zoo{b0011/equatorial_11_gon}{equatorial 11-gon} & $D_{11h}$ & 12.56634432 & 8 & 10 & --\\
12 & \zoo{b0012/equatorial_12_gon}{equatorial 12-gon} & $D_{12h}$ & 12.56636200 & 9 & 11 & --\\
13 & \zoo{b0013/equatorial_13_gon}{equatorial 13-gon} & $D_{13h}$ & 12.56636776 & -- & 10 & --\\
14 & \zoo{b0014/equatorial_14_gon}{equatorial 14-gon} & $D_{14h}$ & 12.56636966 & -- & 9 & --\\
15 & \zoo{b0015/equatorial_15_gon}{equatorial 15-gon} & $D_{15h}$ & 12.56637029 & -- & 10 & --\\
16 & \zoo{b0016/equatorial_16_gon}{equatorial 16-gon} & $D_{16h}$ & 12.56637050 & -- & 9 & --\\
18 & \zoo{b0018/equatorial_18_gon}{equatorial 18-gon} & $D_{18h}$ & 12.56637060 & -- & 9 & --\\
20 & \zoo{b0020/equatorial_20_gon}{equatorial 20-gon} & $D_{20h}$ & 12.56637061 & -- & 9 & --\\
\addlinespace[4pt]
\multicolumn{6}{@{}l}{\textbf{Pyramids}}\\[2pt]
2 & \zoo{b0002/antipodal_pair}{antipodal pair} & $D_{\infty h}$ & \textit{10.47478066} & 0 & 13 & \ref{fig:best-a}\\
5 & \zoo{b0005/square_pyramid}{square pyramid} & $C_{4v}$ & 14.43163401 & 1 & 13 & \ref{fig:type-pyramid}\\
6 & \zoo{b0006/pentagonal_pyramid}{pentagonal pyramid} & $C_{5v}$ & 14.81894243 & 2 & 12 & \ref{fig:type-pyramid}\\
7 & \zoo{b0007/hexagonal_pyramid}{hexagonal pyramid} & $C_{6v}$ & 15.04160665 & 3 & 13 & \ref{fig:type-pyramid}\\
8 & \zoo{b0008/heptagonal_pyramid}{heptagonal pyramid} & $C_{7v}$ & 15.18338502 & 4 & 12 & \ref{fig:type-pyramid}\\
9 & \zoo{b0009/octagonal_pyramid}{octagonal pyramid} & $C_{8v}$ & 15.28070739 & 5 & 12 & \ref{fig:type-pyramid}\\
10 & \zoo{b0010/9_gonal_pyramid}{9-gonal pyramid} & $C_{9v}$ & 15.35132872 & 6 & 12 & --\\
11 & \zoo{b0011/10_gonal_pyramid}{10-gonal pyramid} & $C_{10v}$ & 15.40477635 & 7 & 12 & --\\
12 & \zoo{b0012/11_gonal_pyramid}{11-gonal pyramid} & $C_{11v}$ & 15.44657189 & 8 & 12 & --\\
13 & \zoo{b0013/12_gonal_pyramid}{12-gonal pyramid} & $C_{12v}$ & 15.48011866 & 9 & 10 & \ref{fig:type-pyramid}\\
14 & \zoo{b0014/13_gonal_pyramid}{13-gonal pyramid} & $C_{13v}$ & 15.50762152 & 10 & 11 & --\\
15 & \zoo{b0015/14_gonal_pyramid}{14-gonal pyramid} & $C_{14v}$ & 15.53056857 & 11 & 11 & --\\
16 & \zoo{b0016/15_gonal_pyramid}{15-gonal pyramid} & $C_{15v}$ & 15.54999901 & 12 & 12 & --\\
17 & \zoo{b0017/16_gonal_pyramid}{16-gonal pyramid} & $C_{16v}$ & 15.56665978 & 13 & 11 & \ref{fig:type-pyramid}\\
21 & \zoo{b0021/20_gonal_pyramid}{20-gonal pyramid} & $C_{20v}$ & 15.61479310 & 17 & 10 & \ref{fig:type-pyramid}\\
\addlinespace[4pt]
\multicolumn{6}{@{}l}{\textbf{Bipyramids}}\\[2pt]
5 & \zoo{b0005/triangular_bipyramid}{triangular bipyramid} & $D_{3h}$ & \textit{14.46866597} & 0 & 12 & \ref{fig:best-a}\\
7 & \zoo{b0007/pentagonal_bipyramid}{pentagonal bipyramid} & $D_{5h}$ & \textit{15.95213199} & 0 & 12 & \ref{fig:best-a}\\
8 & \zoo{b0008/hexagonal_bipyramid}{hexagonal bipyramid} & $D_{6h}$ & 16.27617583 & 3 & 12 & \ref{fig:type-bipyramid}\\
9 & \zoo{b0009/heptagonal_bipyramid}{heptagonal bipyramid} & $D_{7h}$ & 16.49293784 & 4 & 10 & --\\
10 & \zoo{b0010/octagonal_bipyramid}{octagonal bipyramid} & $D_{8h}$ & 16.64732271 & 5 & 12 & \ref{fig:type-bipyramid}\\
11 & \zoo{b0011/9_gonal_bipyramid}{9-gonal bipyramid} & $D_{9h}$ & 16.76257176 & 6 & 9 & --\\
12 & \zoo{b0012/10_gonal_bipyramid}{10-gonal bipyramid} & $D_{10h}$ & 16.85176711 & 7 & 12 & --\\
18 & \zoo{b0018/16_gonal_bipyramid}{16-gonal bipyramid} & $D_{16h}$ & 17.13368319 & 13 & 12 & \ref{fig:type-bipyramid}\\
20 & \zoo{b0020/18_gonal_bipyramid}{18-gonal bipyramid} & $D_{18h}$ & 17.18276531 & 15 & 11 & \ref{fig:type-bipyramid}\\
26 & \zoo{b0026/24_gonal_bipyramid}{24-gonal bipyramid} & $D_{24h}$ & 17.27812713 & 21 & 12 & \ref{fig:type-bipyramid}\\
32 & \zoo{b0032/30_gonal_bipyramid}{30-gonal bipyramid} & $D_{30h}$ & 17.33359467 & -- & 11 & \ref{fig:type-bipyramid}\\
\addlinespace[4pt]
\multicolumn{6}{@{}l}{\textbf{Prisms}}\\[2pt]
6 & \zoo{b0006/triangular_prism}{triangular prism} & $D_{3h}$ & 15.12598742 & 1 & 13 & \ref{fig:type-prism}\\
7 & \zoo{b0007/triangular_prism_plus_1pc}{triangular prism + 1PC} & $C_{3v}$ & 15.45115152 & 3 & 12 & \ref{fig:type-prism}\\
8 & \zoo{b0008/triangular_prism_plus_2pc}{triangular prism + 2PC} & $D_{3h}$ & 15.68321560 & 3 & 12 & \ref{fig:type-prism}\\
9 & \zoo{b0009/square_prism_plus_1pc}{square prism + 1PC} & $C_{4v}$ & 16.84599277 & 1 & 11 & --\\
10 & \zoo{b0010/square_prism_plus_2pc}{square prism + 2PC} & $D_{4h}$ & 17.16125175 & 3 & 12 & --\\
10 & \zoo{b0010/pentagonal_prism}{pentagonal prism} & $D_{5h}$ & 17.15927759 & 3 & 13 & \ref{fig:type-prism}\\
11 & \zoo{b0011/pentagonal_prism_plus_1pc}{pentagonal prism + 1PC} & $C_{5v}$ & 17.64159745 & 3 & 13 & --\\
12 & \zoo{b0012/pentagonal_prism_plus_2pc}{pentagonal prism + 2PC} & $D_{5h}$ & 18.00475200 & 1 & 12 & --\\
12 & \zoo{b0012/hexagonal_prism}{hexagonal prism} & $D_{6h}$ & 17.62680627 & 5 & 12 & --\\
13 & \zoo{b0013/hexagonal_prism_plus_1pc}{hexagonal prism + 1PC} & $C_{6v}$ & 18.15187720 & 4 & 12 & --\\
14 & \zoo{b0014/hexagonal_prism_plus_2pc}{hexagonal prism + 2PC} & $D_{6h}$ & 18.54537110 & 2 & 12 & --\\
14 & \zoo{b0014/heptagonal_prism}{heptagonal prism} & $D_{7h}$ & 17.94983940 & 9 & 12 & \ref{fig:type-prism}\\
15 & \zoo{b0015/heptagonal_prism_plus_1pc}{heptagonal prism + 1PC} & $C_{7v}$ & 18.50736333 & 7 & 12 & --\\
16 & \zoo{b0016/heptagonal_prism_plus_2pc}{heptagonal prism + 2PC} & $D_{7h}$ & 18.92222890 & 5 & 12 & --\\
16 & \zoo{b0016/octagonal_prism}{octagonal prism} & $D_{8h}$ & 18.18597888 & 11 & 12 & --\\
17 & \zoo{b0017/octagonal_prism_plus_1pc}{octagonal prism + 1PC} & $C_{8v}$ & 18.76949374 & 9 & 12 & --\\
18 & \zoo{b0018/octagonal_prism_plus_2pc}{octagonal prism + 2PC} & $D_{8h}$ & 19.20077284 & 7 & 12 & --\\
18 & \zoo{b0018/9_gonal_prism}{9-gonal prism} & $D_{9h}$ & 18.36581718 & 13 & 12 & \ref{fig:type-prism}\\
19 & \zoo{b0019/9_gonal_prism_plus_1pc}{9-gonal prism + 1PC} & $C_{9v}$ & 18.97080398 & 11 & 12 & --\\
20 & \zoo{b0020/9_gonal_prism_plus_2pc}{9-gonal prism + 2PC} & $D_{9h}$ & 19.41539505 & 11 & 12 & \ref{fig:type-prism}\\
21 & \zoo{b0021/10_gonal_prism_plus_1pc}{10-gonal prism + 1PC} & $C_{10v}$ & 19.13021129 & 13 & 12 & --\\
23 & \zoo{b0023/11_gonal_prism_plus_1pc}{11-gonal prism + 1PC} & $C_{11v}$ & 19.25951094 & 15 & 12 & --\\
\addlinespace[4pt]
\multicolumn{6}{@{}l}{\textbf{Antiprisms}}\\[2pt]
7 & \zoo{b0007/triangular_antiprism_plus_1pc}{triangular antiprism + 1PC} & $C_{3v}$ & 15.91577614 & 2 & 12 & \ref{fig:type-antiprism}\\
8 & \zoo{b0008/square_antiprism}{square antiprism} & $D_{4d}$ & \textit{16.49540952} & 0 & 12 & \ref{fig:best-a}\\
9 & \zoo{b0009/square_antiprism_plus_1pc}{square antiprism + 1PC} & $C_{4v}$ & 16.96874135 & 1 & 12 & --\\
10 & \zoo{b0010/square_antiprism_plus_2pc}{square antiprism + 2PC} & $D_{4d}$ & \textit{17.37203030} & 0 & 12 & \ref{fig:best-a}\\
10 & \zoo{b0010/pentagonal_antiprism}{pentagonal antiprism} & $D_{5d}$ & 17.17571147 & 2 & 13 & \ref{fig:type-antiprism}\\
11 & \zoo{b0011/pentagonal_antiprism_plus_1pc}{pentagonal antiprism + 1PC} & $C_{5v}$ & 17.67481751 & 2 & 12 & --\\
12 & \zoo{b0012/hexagonal_antiprism}{hexagonal antiprism} & $D_{6d}$ & 17.63066847 & 6 & 12 & --\\
13 & \zoo{b0013/hexagonal_antiprism_plus_1pc}{hexagonal antiprism + 1PC} & $C_{6v}$ & 18.16106926 & 3 & 12 & --\\
14 & \zoo{b0014/hexagonal_antiprism_plus_2pc}{hexagonal antiprism + 2PC} & $D_{6d}$ & \textit{18.56544315} & 0 & 12 & \ref{fig:best-a}\\
14 & \zoo{b0014/heptagonal_antiprism}{heptagonal antiprism} & $D_{7d}$ & 17.95076062 & 8 & 12 & --\\
15 & \zoo{b0015/heptagonal_antiprism_plus_1pc}{heptagonal antiprism + 1PC} & $C_{7v}$ & 18.50995876 & 6 & 12 & --\\
16 & \zoo{b0016/heptagonal_antiprism_plus_2pc}{heptagonal antiprism + 2PC} & $D_{7d}$ & 18.92868081 & 4 & 12 & --\\
16 & \zoo{b0016/octagonal_antiprism}{octagonal antiprism} & $D_{8d}$ & 18.18620192 & 10 & 12 & --\\
17 & \zoo{b0017/octagonal_antiprism_plus_1pc}{octagonal antiprism + 1PC} & $C_{8v}$ & 18.77023943 & 8 & 12 & --\\
18 & \zoo{b0018/octagonal_antiprism_plus_2pc}{octagonal antiprism + 2PC} & $D_{8d}$ & 19.20288658 & 8 & 12 & --\\
18 & \zoo{b0018/9_gonal_antiprism}{9-gonal antiprism} & $D_{9d}$ & 18.36587192 & 12 & 12 & --\\
19 & \zoo{b0019/9_gonal_antiprism_plus_1pc}{9-gonal antiprism + 1PC} & $C_{9v}$ & 18.97102138 & 10 & 11 & \ref{fig:type-antiprism}\\
20 & \zoo{b0020/9_gonal_antiprism_plus_2pc}{9-gonal antiprism + 2PC} & $D_{9d}$ & 19.41609831 & 10 & 12 & --\\
20 & \zoo{b0020/10_gonal_antiprism}{10-gonal antiprism} & $D_{10d}$ & 18.50715415 & 14 & 12 & --\\
21 & \zoo{b0021/10_gonal_antiprism_plus_1pc}{10-gonal antiprism + 1PC} & $C_{10v}$ & 19.13027544 & 12 & 12 & --\\
22 & \zoo{b0022/10_gonal_antiprism_plus_2pc}{10-gonal antiprism + 2PC} & $D_{10d}$ & 19.58618325 & 12 & 12 & --\\
22 & \zoo{b0022/11_gonal_antiprism}{11-gonal antiprism} & $D_{11d}$ & 18.62100292 & 16 & 12 & --\\
23 & \zoo{b0023/11_gonal_antiprism_plus_1pc}{11-gonal antiprism + 1PC} & $C_{11v}$ & 19.25953007 & 14 & 12 & --\\
24 & \zoo{b0024/11_gonal_antiprism_plus_2pc}{11-gonal antiprism + 2PC} & $D_{11d}$ & 19.72483665 & 14 & 11 & --\\
24 & \zoo{b0024/12_gonal_antiprism}{12-gonal antiprism} & $D_{12d}$ & 18.71461596 & 18 & 11 & --\\
26 & \zoo{b0026/12_gonal_antiprism_plus_2pc}{12-gonal antiprism + 2PC} & $D_{12d}$ & 19.83995001 & 16 & 11 & --\\
26 & \zoo{b0026/13_gonal_antiprism}{13-gonal antiprism} & $D_{13d}$ & 18.79289896 & 20 & 11 & --\\
27 & \zoo{b0027/13_gonal_antiprism_plus_1pc}{13-gonal antiprism + 1PC} & $C_{13v}$ & 19.45635380 & 18 & 12 & --\\
28 & \zoo{b0028/13_gonal_antiprism_plus_2pc}{13-gonal antiprism + 2PC} & $D_{13d}$ & 19.93700822 & 18 & 11 & --\\
28 & \zoo{b0028/14_gonal_antiprism}{14-gonal antiprism} & $D_{14d}$ & 18.85930278 & -- & 12 & --\\
29 & \zoo{b0029/14_gonal_antiprism_plus_1pc}{14-gonal antiprism + 1PC} & $C_{14v}$ & 19.53295856 & 20 & 12 & --\\
30 & \zoo{b0030/14_gonal_antiprism_plus_2pc}{14-gonal antiprism + 2PC} & $D_{14d}$ & 20.01992869 & 20 & 11 & --\\
30 & \zoo{b0030/15_gonal_antiprism}{15-gonal antiprism} & $D_{15d}$ & 18.91632175 & -- & 11 & --\\
31 & \zoo{b0031/15_gonal_antiprism_plus_1pc}{15-gonal antiprism + 1PC} & $C_{15v}$ & 19.59900546 & 22 & 12 & --\\
32 & \zoo{b0032/15_gonal_antiprism_plus_2pc}{15-gonal antiprism + 2PC} & $D_{15d}$ & 20.09157846 & -- & 11 & --\\
32 & \zoo{b0032/16_gonal_antiprism}{16-gonal antiprism} & $D_{16d}$ & 18.96580166 & -- & 12 & --\\
33 & \zoo{b0033/16_gonal_antiprism_plus_1pc}{16-gonal antiprism + 1PC} & $C_{16v}$ & 19.65652786 & -- & 11 & --\\
34 & \zoo{b0034/16_gonal_antiprism_plus_2pc}{16-gonal antiprism + 2PC} & $D_{16d}$ & 20.15410172 & -- & 11 & \ref{fig:type-antiprism}\\
35 & \zoo{b0035/17_gonal_antiprism_plus_1pc}{17-gonal antiprism + 1PC} & $C_{17v}$ & 19.70707031 & -- & 11 & \ref{fig:type-prism}\\
37 & \zoo{b0037/18_gonal_antiprism_plus_1pc}{18-gonal antiprism + 1PC} & $C_{18v}$ & 19.75182665 & -- & 11 & --\\
39 & \zoo{b0039/19_gonal_antiprism_plus_1pc}{19-gonal antiprism + 1PC} & $C_{19v}$ & 19.79173390 & -- & 11 & --\\
40 & \zoo{b0040/19_gonal_antiprism_plus_2pc}{19-gonal antiprism + 2PC} & $D_{19d}$ & 20.30151860 & -- & 11 & --\\
40 & \zoo{b0040/20_gonal_antiprism}{20-gonal antiprism} & $D_{20d}$ & 19.11187927 & -- & 11 & \ref{fig:type-antiprism}\\
42 & \zoo{b0042/20_gonal_antiprism_plus_2pc}{20-gonal antiprism + 2PC} & $D_{20d}$ & 20.34066584 & 32 & 11 & --\\
\addlinespace[4pt]
\multicolumn{6}{@{}l}{\textbf{Stacked latitudinal rings}}\\[2pt]
7 & \zoo{b0007/quadruple_plus_pair_plus_1pc}{quadruple + pair + 1PC} & $C_{2v}$ & 15.91819689 & 1 & 10 & --\\
7 & \zoo{b0007/quadruple_plus_pair_plus_1pc_pair_near_the_pole}{quadruple + polar pair + 1PC} & $C_{2v}$ & 15.46032046 & 2 & 13 & --\\
8 & \zoo{b0008/four_stacked_pairs}{4 pairs} & $D_{2d}$ & 16.49051351 & 1 & 12 & \ref{fig:type-stacked}\\
8 & \zoo{b0008/three_pairs_plus_2pc}{three pairs + 2PC} & $D_{2d}$ & 16.29009761 & 1 & 12 & \ref{fig:type-stacked}\\
8 & \zoo{b0008/quadruple_plus_pair_plus_2pc}{quadruple + pair + 2PC} & $C_{2v}$ & 16.28575637 & 2 & 9 & --\\
8 & \zoo{b0008/quadruple_plus_two_pairs_staggered}{quadruple + two pairs, staggered} & $C_{2v}$ & 15.71150575 & 3 & 12 & --\\
8 & \zoo{b0008/quadruple_plus_two_pairs}{quadruple + two pairs} & $C_{2v}$ & 15.68321410 & 4 & 9 & --\\
9 & \zoo{b0009/three_stacked_triangles}{3 triangles} & $D_{3h}$ & \textit{16.97069066} & 0 & 13 & \ref{fig:best-a}\\
10 & \zoo{b0010/three_stacked_triangles_plus_1pc}{3 triangles + 1PC} & $C_{3v}$ & 17.36763805 & 2 & 13 & --\\
11 & \zoo{b0011/three_stacked_triangles_plus_2pc}{3 triangles + 2PC} & $D_{3h}$ & 17.67498004 & 2 & 12 & --\\
12 & \zoo{b0012/three_stacked_squares}{3 squares} & $O_h$ & 17.97506201 & 4 & 12 & \ref{fig:type-polyhedral}\\
13 & \zoo{b0013/three_stacked_squares_plus_1pc}{3 squares + 1PC} & $C_{4v}$ & 18.28105516 & 3 & 12 & \ref{fig:type-stacked}\\
14 & \zoo{b0014/three_stacked_squares_plus_2pc}{3 squares + 2PC} & $O_h$ & 18.55117279 & 3 & 12 & --\\
15 & \zoo{b0015/five_stacked_triangles}{5 triangles} & $D_{3h}$ & \textit{18.79339513} & 0 & 12 & \ref{fig:best-a}\\
15 & \zoo{b0015/three_stacked_pentagons}{3 pentagons} & $D_{5h}$ & 18.70485004 & 4 & 12 & --\\
16 & \zoo{b0016/four_stacked_squares}{4 squares} & $D_{4d}$ & 19.00297008 & 0 & 9 & \ref{fig:type-stacked}\\
16 & \zoo{b0016/three_stacked_pentagons_plus_1pc}{3 pentagons + 1PC} & $C_{5v}$ & 18.96676603 & 2 & 12 & --\\
17 & \zoo{b0017/three_stacked_pentagons_plus_2pc}{3 pentagons + 2PC} & $D_{5h}$ & 19.20095197 & 2 & 12 & --\\
18 & \zoo{b0018/four_stacked_squares_plus_2pc}{4 squares + 2PC} & $D_{4d}$ & \textit{19.38539453} & 0 & 9 & \ref{fig:best-b}\\
18 & \zoo{b0018/six_stacked_triangles}{6 triangles} & $D_{3d}$ & 19.33709715 & 1 & 12 & --\\
18 & \zoo{b0018/three_stacked_hexagons}{3 hexagons} & $D_{6h}$ & 19.24666524 & 6 & 12 & --\\
19 & \zoo{b0019/three_stacked_hexagons_plus_1pc}{3 hexagons + 1PC} & $C_{6v}$ & 19.49029702 & 3 & 11 & --\\
20 & \zoo{b0020/three_stacked_hexagons_plus_2pc}{3 hexagons + 2PC} & $D_{6h}$ & 19.70418876 & 1 & 12 & --\\
20 & \zoo{b0020/five_stacked_squares}{5 squares} & $D_{4h}$ & 19.70374025 & 2 & 10 & \ref{fig:type-stacked}\\
20 & \zoo{b0020/four_stacked_pentagons}{4 pentagons} & $D_{5d}$ & 19.65221836 & 6 & 12 & --\\
20 & \zoo{b0020/six_stacked_triangles_plus_2pc_alternating}{6 triangles + 2PC, alt.} & $D_{6h}$ & 19.60976866 & 4 & 12 & --\\
20 & \zoo{b0020/six_stacked_triangles_plus_2pc}{6 triangles + 2PC} & $D_{3h}$ & 19.50838431 & 9 & 11 & \ref{fig:type-stacked}\\
21 & \zoo{b0021/three_stacked_heptagons}{3 heptagons} & $D_{7h}$ & 19.65480708 & 8 & 12 & --\\
22 & \zoo{b0022/four_stacked_pentagons_plus_2pc}{4 pentagons + 2PC} & $D_{5d}$ & 19.99205479 & 0 & 12 & --\\
22 & \zoo{b0022/three_stacked_7_gons_plus_1pc}{3 7-gons + 1PC} & $C_{7v}$ & 19.89228042 & 4 & 12 & --\\
23 & \zoo{b0023/three_stacked_heptagons_plus_2pc}{3 heptagons + 2PC} & $D_{7h}$ & 20.09724473 & 2 & 11 & --\\
24 & \zoo{b0024/six_stacked_squares}{6 squares} & $D_{4d}$ & 20.24291718 & 1 & 9 & --\\
24 & \zoo{b0024/four_stacked_hexagons}{4 hexagons} & $D_{6d}$ & 20.15142930 & 8 & 12 & --\\
25 & \zoo{b0025/five_stacked_pentagons}{5 pentagons} & $D_{5h}$ & 20.32141112 & 6 & 12 & --\\
25 & \zoo{b0025/three_stacked_8_gons_plus_1pc}{3 8-gons + 1PC} & $C_{8v}$ & 20.20538352 & 8 & 12 & --\\
26 & \zoo{b0026/four_stacked_hexagons_plus_2pc}{4 hexagons + 2PC} & $D_{6d}$ & 20.45390415 & 6 & 12 & --\\
26 & \zoo{b0026/six_stacked_squares_plus_2pc_20p4521}{6 squares + 2PC} & $D_{4d}$ & 20.45211659 & 3 & 11 & --\\
26 & \zoo{b0026/six_stacked_squares_plus_2pc_20p4122}{6 squares + 2PC} & $D_{4h}$ & 20.41219517 & 6 & 12 & --\\
27 & \zoo{b0027/three_stacked_9_gons}{3 9-gons} & $D_{9h}$ & 20.21668512 & 12 & 12 & --\\
28 & \zoo{b0028/four_stacked_heptagons_20p5476}{4 heptagons} & $D_{7d}$ & 20.54757463 & 8 & 12 & --\\
28 & \zoo{b0028/four_stacked_heptagons_20p4947}{4 heptagons} & $D_{7h}$ & 20.49466892 & 11 & 11 & --\\
28 & \zoo{b0028/three_stacked_9_gons_plus_1pc}{3 9-gons + 1PC} & $C_{9v}$ & 20.45393289 & 10 & 12 & --\\
29 & \zoo{b0029/three_stacked_9_gons_plus_2pc}{3 9-gons + 2PC} & $D_{9h}$ & 20.65478913 & 10 & 11 & --\\
30 & \zoo{b0030/six_stacked_pentagons}{6 pentagons} & $D_{5d}$ & 20.82752606 & 4 & 11 & --\\
30 & \zoo{b0030/five_stacked_hexagons}{5 hexagons} & $D_{6h}$ & 20.51696912 & 15 & 12 & --\\
32 & \zoo{b0032/six_stacked_pentagons_plus_2pc}{6 pentagons + 2PC} & $D_{5h}$ & 20.98446690 & 3 & 10 & --\\
33 & \zoo{b0033/three_stacked_11_gons}{3 11-gons} & $D_{11h}$ & 20.58014530 & 20 & 10 & --\\
34 & \zoo{b0034/four_stacked_octagons_plus_2pc}{4 octagons + 2PC} & $D_{8d}$ & 21.12308897 & -- & 11 & --\\
34 & \zoo{b0034/three_stacked_11_gons_plus_1pc}{3 11-gons + 1PC} & $C_{11v}$ & 20.82083077 & -- & 11 & --\\
35 & \zoo{b0035/three_stacked_11_gons_plus_2pc}{3 11-gons + 2PC} & $D_{11h}$ & 21.02271261 & -- & 11 & --\\
35 & \zoo{b0035/five_stacked_heptagons}{5 heptagons} & $D_{7h}$ & 20.99936713 & -- & 12 & --\\
36 & \zoo{b0036/six_stacked_hexagons}{6 hexagons} & $D_{6d}$ & 21.25608219 & -- & 10 & --\\
38 & \zoo{b0038/six_stacked_hexagons_plus_2pc_21p4221}{6 hexagons + 2PC} & $D_{6d}$ & \textit{21.42213538} & 0 & 11 & \ref{fig:best-c}\\
38 & \zoo{b0038/six_stacked_hexagons_plus_2pc_20p9935}{6 hexagons + 2PC} & $D_{6h}$ & 20.99346879 & -- & 11 & --\\
39 & \zoo{b0039/three_stacked_13_gons}{3 13-gons} & $D_{13h}$ & 20.83311527 & -- & 10 & --\\
40 & \zoo{b0040/five_stacked_octagons}{5 octagons} & $D_{8h}$ & 21.36190951 & -- & 11 & --\\
40 & \zoo{b0040/four_stacked_10_gons_21p3336}{4 10-gons} & $D_{10d}$ & 21.33356007 & -- & 11 & --\\
40 & \zoo{b0040/four_stacked_10_gons_21p3274}{4 10-gons} & $D_{10h}$ & 21.32737717 & -- & 11 & --\\
40 & \zoo{b0040/three_stacked_13_gons_plus_1pc}{3 13-gons + 1PC} & $C_{13v}$ & 21.07732753 & -- & 11 & --\\
49 & \zoo{b0049/three_stacked_16_gons_plus_1pc}{3 16-gons + 1PC} & $C_{16v}$ & 21.34327300 & -- & 11 & --\\
61 & \zoo{b0061/three_stacked_20_gons_plus_1pc}{3 20-gons + 1PC} & $C_{20v}$ & 21.57466863 & -- & 11 & --\\
\addlinespace[4pt]
\multicolumn{6}{@{}l}{\textbf{Classes of the Fekete configurations}}\\[2pt]
19 & \zoo{b0019/fekete_class_19p5446}{Fekete class} & $C_{2v}$ & \textit{19.54459459} & 0 & 11 & \ref{fig:best-b}\\
19 & \zoo{b0019/fekete_class_19p5362}{Fekete class} & $C_{2v}$ & 19.53616762 & 2 & 11 & --\\
20 & \zoo{b0020/fekete_class}{Fekete class} & $D_{3h}$ & \textit{19.70770682} & 0 & 12 & \ref{fig:best-b}\\
21 & \zoo{b0021/fekete_class}{Fekete class} & $C_{2v}$ & \textit{19.85473178} & 0 & 12 & \ref{fig:best-b}\\
22 & \zoo{b0022/fekete_class_19p9935}{Fekete class} & $T_d$ & \textit{19.99354385} & 0 & 12 & \ref{fig:best-b}\\
22 & \zoo{b0022/fekete_class_19p7018}{Fekete class} & $T_d$ & 19.70178848 & 10 & 12 & --\\
23 & \zoo{b0023/fekete_class}{Fekete class} & $D_{3}$ & \textit{20.11635049} & 0 & 12 & \ref{fig:best-b}\\
24 & \zoo{b0024/fekete_class}{Fekete class} & $O$ & \textit{20.24322425} & 0 & 12 & \ref{fig:best-b}\\
25 & \zoo{b0025/fekete_class}{Fekete class} & $C_s$ & 20.35311125 & 1 & 11 & --\\
26 & \zoo{b0026/fekete_class}{Fekete class} & $C_{2}$ & \textit{20.46172826} & 0 & 12 & \ref{fig:best-b}\\
27 & \zoo{b0027/fekete_class}{Fekete class} & $D_{5h}$ & \textit{20.56911733} & 0 & 12 & \ref{fig:best-b}\\
28 & \zoo{b0028/fekete_class}{Fekete class} & $T$ & \textit{20.66487632} & 0 & 12 & \ref{fig:best-b}\\
29 & \zoo{b0029/fekete_class}{Fekete class} & $D_{3}$ & \textit{20.75364318} & 0 & 11 & \ref{fig:best-b}\\
30 & \zoo{b0030/fekete_class}{Fekete class} & $D_{2}$ & \textit{20.84340804} & 0 & 12 & \ref{fig:best-b}\\
31 & \zoo{b0031/fekete_class_20p9278}{Fekete class} & $C_{3v}$ & \textit{20.92777889} & 0 & 11 & \ref{fig:best-b}\\
31 & \zoo{b0031/fekete_class_20p9163}{Fekete class} & $C_{3v}$ & 20.91628299 & 4 & 11 & --\\
31 & \zoo{b0031/fekete_class_20p9024}{Fekete class} & $C_{3v}$ & 20.90244106 & 8 & 11 & --\\
33 & \zoo{b0033/fekete_class}{Fekete class} & $C_s$ & 21.08353676 & 1 & 11 & --\\
34 & \zoo{b0034/fekete_class}{Fekete class} & $D_{2}$ & \textit{21.15725596} & 0 & 12 & \ref{fig:best-c}\\
35 & \zoo{b0035/fekete_class}{Fekete class} & $C_{2}$ & \textit{21.22645244} & 0 & 11 & \ref{fig:best-c}\\
36 & \zoo{b0036/fekete_class}{Fekete class} & $D_{2}$ & \textit{21.29451111} & 1 & 12 & \ref{fig:best-c}\\
37 & \zoo{b0037/fekete_class_21p3600}{Fekete class} & $D_{5h}$ & \textit{21.35998448} & 0 & 12 & \ref{fig:best-c}\\
37 & \zoo{b0037/fekete_class_21p1917}{Fekete class} & $D_{5h}$ & 21.19166865 & -- & 11 & --\\
39 & \zoo{b0039/fekete_class_21p4831}{Fekete class} & $D_{3h}$ & \textit{21.48308018} & 0 & 11 & \ref{fig:best-c}\\
39 & \zoo{b0039/fekete_class_21p4727}{Fekete class} & $D_{3h}$ & 21.47271751 & -- & 12 & --\\
39 & \zoo{b0039/fekete_class_21p4257}{Fekete class} & $D_{3h}$ & 21.42573429 & -- & 11 & --\\
39 & \zoo{b0039/fekete_class_21p3318}{Fekete class} & $D_{3h}$ & 21.33182995 & -- & 11 & --\\
40 & \zoo{b0040/fekete_class_21p5411}{Fekete class} & $T_d$ & \textit{21.54113373} & 0 & 11 & \ref{fig:best-c}\\
40 & \zoo{b0040/fekete_class_21p4941}{Fekete class} & $T_d$ & 21.49406074 & -- & 12 & --\\
40 & \zoo{b0040/fekete_class_21p4122}{Fekete class} & $T_d$ & 21.41215845 & -- & 11 & --\\
40 & \zoo{b0040/fekete_class_21p2717}{Fekete class} & $T_d$ & 21.27170505 & -- & 11 & --\\
\addlinespace[4pt]
\multicolumn{6}{@{}l}{\textbf{Other configurations}}\\[2pt]
8 & \zoo{b0008/cs_16p2814}{$C_s$: $3\times2+2\times1$} & $C_s$ & 16.28141689 & 2 & 12 & --\\
8 & \zoo{b0008/cs_16p2812}{$C_s$: $3\times2+2\times1$} & $C_s$ & 16.28123654 & 3 & 12 & --\\
9 & \zoo{b0009/c2v_16p8408}{$C_{2v}$: $4+2\times2+1$} & $C_{2v}$ & 16.84081989 & 2 & 12 & --\\
9 & \zoo{b0009/cs_16p8224}{$C_s$: $3\times2+3\times1$} & $C_s$ & 16.82240540 & 2 & 12 & --\\
9 & \zoo{b0009/cs_16p82239}{$C_s$: $4\times2+1$} & $C_s$ & 16.82238669 & 3 & 12 & --\\
9 & \zoo{b0009/c3v_16p6111}{$C_{3v}$: $3\times3$} & $C_{3v}$ & 16.61110182 & 3 & 13 & --\\
9 & \zoo{b0009/c3v_16p5361}{$C_{3v}$: $6+3$} & $C_{3v}$ & 16.53610663 & 3 & 12 & --\\
9 & \zoo{b0009/c3v_16p5230}{$C_{3v}$: $3\times3$} & $C_{3v}$ & 16.52295783 & 4 & 12 & --\\
9 & \zoo{b0009/d3h_15p9547}{$D_{3h}$: $6+3$} & $D_{3h}$ & 15.95465009 & 4 & 12 & --\\
9 & \zoo{b0009/cs_15p8633}{$C_s$: $4\times2+1$} & $C_s$ & 15.86332538 & 4 & 12 & --\\
9 & \zoo{b0009/cs_15p8629}{$C_s$: $3\times2+3\times1$} & $C_s$ & 15.86294017 & 5 & 12 & --\\
10 & \zoo{b0010/c2v_17p3680}{$C_{2v}$: $4+3\times2$} & $C_{2v}$ & 17.36797545 & 1 & 11 & --\\
10 & \zoo{b0010/c3v_17p2779}{$C_{3v}$: $6+3+1$} & $C_{3v}$ & 17.27786230 & 1 & 11 & --\\
10 & \zoo{b0010/c2v_17p2769}{$C_{2v}$: $4+3\times2$} & $C_{2v}$ & 17.27686572 & 2 & 12 & --\\
10 & \zoo{b0010/cs_17p23495}{$C_s$: $4\times2+2\times1$} & $C_s$ & 17.23495171 & 2 & 12 & --\\
10 & \zoo{b0010/cs_17p2350}{$C_s$: $3\times2+4\times1$} & $C_s$ & 17.23495080 & 3 & 11 & --\\
10 & \zoo{b0010/d2h_17p1680}{$D_{2h}$: $2\times4+2$} & $D_{2h}$ & 17.16802835 & 3 & 12 & --\\
10 & \zoo{b0010/c2v_17p0583}{$C_{2v}$: $4+3\times2$} & $C_{2v}$ & 17.05834704 & 3 & 12 & --\\
10 & \zoo{b0010/c2v_17p0578}{$C_{2v}$: $2\times4+2$} & $C_{2v}$ & 17.05777516 & 4 & 12 & --\\
10 & \zoo{b0010/c3v_16p8513}{$C_{3v}$: $6+3+1$} & $C_{3v}$ & 16.85130687 & 3 & 12 & --\\
10 & \zoo{b0010/cs_16p8512}{$C_s$: $4\times2+2\times1$} & $C_s$ & 16.85123562 & 4 & 12 & --\\
10 & \zoo{b0010/c3v_16p8041}{$C_{3v}$: $3\times3+1$} & $C_{3v}$ & 16.80405847 & 5 & 12 & --\\
10 & \zoo{b0010/td_16p8035}{$T_d$: $6+4$} & $T_d$ & 16.80346515 & 6 & 12 & --\\
10 & \zoo{b0010/cs_16p69742}{$C_s$: $4\times2+2\times1$} & $C_s$ & 16.69741899 & 4 & 12 & --\\
10 & \zoo{b0010/cs_16p6974}{$C_s$: $4\times2+2\times1$} & $C_s$ & 16.69741644 & 5 & 11 & --\\
10 & \zoo{b0010/c2v_16p0880}{$C_{2v}$: $4+3\times2$} & $C_{2v}$ & 16.08800759 & 5 & 12 & --\\
10 & \zoo{b0010/c3v_16p0791}{$C_{3v}$: $3\times3+1$} & $C_{3v}$ & 16.07906572 & 6 & 10 & --\\
10 & \zoo{b0010/c2v_15p9769}{$C_{2v}$: $2\times4+2$} & $C_{2v}$ & 15.97690925 & 5 & 12 & --\\
10 & \zoo{b0010/c2v_15p9742}{$C_{2v}$: $4+3\times2$} & $C_{2v}$ & 15.97416080 & 6 & 12 & --\\
11 & \zoo{b0011/c2v_17p6912}{$C_{2v}$: $4+3\times2+1$} & $C_{2v}$ & \textit{17.69115333} & 0 & 12 & \ref{fig:best-a}\\
11 & \zoo{b0011/cs_17p6897}{$C_s$: $4\times2+3\times1$} & $C_s$ & 17.68970529 & 1 & 11 & --\\
12 & \zoo{b0012/d3h_17p9944}{$D_{3h}$: $2\times6$} & $D_{3h}$ & 17.99436595 & 1 & 11 & --\\
12 & \zoo{b0012/c2v_17p9751}{$C_{2v}$: $2\times4+2+2\times1$} & $C_{2v}$ & 17.97511342 & 2 & 12 & --\\
12 & \zoo{b0012/cs_17p9518}{$C_s$: $5\times2+2\times1$} & $C_s$ & 17.95180696 & 2 & 10 & --\\
12 & \zoo{b0012/cs_17p9263}{$C_s$: $4\times2+4\times1$} & $C_s$ & 17.92628595 & 3 & 10 & --\\
13 & \zoo{b0013/c2v_18p3064}{$C_{2v}$: $4+4\times2+1$} & $C_{2v}$ & \textit{18.30638945} & 0 & 12 & \ref{fig:best-a}\\
16 & \zoo{b0016/td_19p0089}{$T_d$: $12+4$} & $T_d$ & \textit{19.00885253} & 0 & 12 & \ref{fig:best-a}\\
17 & \zoo{b0017/c2v_19p2012}{$C_{2v}$: $2\times4+4\times2+1$} & $C_{2v}$ & \textit{19.20121924} & 0 & 12 & \ref{fig:best-a}\\
19 & \zoo{b0019/c2_19p5433}{$C_2$: $9\times2+1$} & $C_{2}$ & 19.54328806 & 1 & 11 & --\\
19 & \zoo{b0019/c2v_19p5422}{$C_{2v}$: $2\times4+5\times2+1$} & $C_{2v}$ & 19.54218732 & 1 & 11 & --\\
25 & \zoo{b0025/c1_20p3531}{$C_1$ (20.3531)} & $C_1$ & \textit{20.35311176} & 0 & 11 & \ref{fig:best-b}\\
33 & \zoo{b0033/c1_21p0835}{$C_1$ (21.0835)} & $C_1$ & \textit{21.08353796} & 0 & 11 & \ref{fig:best-b}\\
\end{supertabular}
\endgroup
\onecolumn

\end{document}